\documentclass[12pt,a4paper]{article}
\usepackage[latin1]{inputenc}
\usepackage[french]{babel}
\usepackage{amsfonts}
\usepackage{amssymb}

\newtheorem{lemme}{Lemme}
\newtheorem{proposition}[lemme]{Proposition}
\newtheorem{theoreme}[lemme]{Th\'eor\`eme}
\newtheorem{corollaire}[lemme]{Corollaire}
\newtheorem{definition}[lemme]{Définition}

\newcommand\dem{\noindent{\it D\'emonstration.}\ }
\newcommand\findem{\hfill$\square$}
\newcommand\hfld[2]{\smash{\mathop{\hbox to 10mm{\rightarrowfill}}
     \limits^{\scriptstyle#1}_{\scriptstyle#2}}}
\newcommand\hflg[2]{\smash{\mathop{\hbox to 10mm{\leftarrowfill}}
     \limits^{\scriptstyle#1}_{\scriptstyle#2}}}
\newcommand\vfld[2]{\llap{$\scriptstyle#1$}
     \left\downarrow\vbox to 4mm{}\right.\rlap{$\scriptstyle #2$}}

\newcommand\rta{\rightarrow}
\newcommand\hrta{\hookrightarrow}
\newcommand\sv{{\scriptscriptstyle\vee}}
\newcommand\la{\langle}
\newcommand\ra{\rangle}
\newcommand\bab{\overline}

\newcommand\calA{{\mathcal A}}
\newcommand\calF{{\mathcal F}}

\newcommand\calL{{\mathcal L}}
\newcommand\calM{{\mathcal M}}
\newcommand\calO{{\mathcal O}}

\newcommand\calE{{\mathcal E}}
\newcommand\calD{{\mathcal D}}

\newcommand\calH{{\mathcal H}}
\newcommand\calV{{\mathcal V}}
\newcommand\calU{{\mathcal U}}
\newcommand\calW{{\mathcal W}}
\newcommand\calX{{\mathcal X}}
\newcommand\calQ{{\mathcal Q}}
\newcommand\calS{{\mathcal S}}
\newcommand\calG{{\mathcal G}}
\newcommand\CC{{\mathbb C}}

\newcommand\FF{{\mathbb F}}

\newcommand\ZZ{{\mathbb Z}}
\newcommand\NN{{\mathbb N}}
\newcommand\QQ{{\mathbb Q}}
\newcommand\bbA{{\mathbb A}}
\newcommand\rmH{{\rm H}}
\newcommand\rmR{{\rm R}}
\newcommand\frakS{\mathfrak S}
\newcommand\frakC{\mathfrak C}

\newcommand\gl{\mathfrak{gl}}
\newcommand\GL{{\rm GL}}
\newcommand\Tr{{\rm Tr}}
\newcommand\Gal{{\rm Gal}}
\newcommand\pr{{\rm pr}}
\newcommand\DFib{\calD\mbox{-}{\rm Fib}}
\newcommand\DHck{\calD\mbox{-}{\rm Hck}}
\newcommand\DCht{\calD\mbox{-}{\rm Cht}}
\newcommand\DICht{\calD^I\mbox{-}{\rm Cht}}
\newcommand\DITrv{\calD_I\mbox{-}{\rm Trv}}
\newcommand\smv{\sv}
\newcommand\Frob{{\rm Frob}}
\newcommand\id{{\rm id}}

\newcommand\Spec{{\rm Spec}}
\newcommand\hookr\hookrightarrow
\newcommand\hookl\hookleftarrow
\newcommand\und{\underline}
\newcommand\inv{{\rm inv}}
\newcommand\isom{\,\smash{\mathop{\hbox to 5mm{\rightarrowfill}}
    \limits^\sim}\,}
\newcommand\ul{{\und\lambda}}

\newcommand\Utau{U_{\langle\tau\rangle}}

\begin{document}

\title{$\calD$-chtoucas de Drinfeld à modifications symétriques et
       identité de changement de base}
\author{Ng\^o Bao Ch\^au \footnote{Département de mathématiques, Université
Paris-Sud 91405 Orsay France}}
\date{\small\em à Gérard Laumon pour son cinquantième anniversaire}
\maketitle

\section*{Introduction}

Soit $F$ un corps local non-archimédien, $E$ une extension
non-ramifiée de degré $r$ de $F$ et $\sigma$ un générateur de
$\Gal(E/F)$. Soit $G=\GL_d$. Pour tous $\delta\in G(E)$, la classe
de conjugaison de $\delta \sigma(\delta)\ldots
\sigma^{r-1}(\delta)$ est $\sigma$-stable et définit donc une
classe de conjugaison dans $G(F)$. Un élément $\gamma\in G(F)$
dans cette classe de conjugaison est appelé une norme de $\delta$.
Supposons que $\gamma$ est semi-simple régulier.

Soit $\calH_F$ l'algèbre des fonctions sphériques c'est-à-dire les
fonctions complexes à support compacts sur $G(F)$ qui sont
invariantes à gauche et à droite par $G(\calO_F)$, $\calO_F$ étant
l'anneau des entiers de $F$. D'après Satake, il existe un
isomorphisme d'algèbres
$$\calH_F\isom \CC[t_1^{\pm 1},\ldots,t_d^{\pm 1}]^{\mathfrak S_d}$$
où l'algèbre polynômiale $\CC[t_1,\ldots,t_d]$ est l'algèbre des
fonctions à support compact sur $T(F)$ qui sont invariantes par
$T(\calO_F)$, $T$ étant le tore diagonal de $G$. De même, on a un
isomorphisme
$$\calH_E\isom \CC[u_1^{\pm 1},\ldots,u_d^{\pm 1}]^{\mathfrak S_d}.$$
L'homomorphisme de changement de base ${\mathbf
b}:\calH_E\rta\calH_F$ est alors défini par transport de
l'homomorphisme $u_i\mapsto t_i^r$ du côté des polynômes
symétriques. Cet homomorphisme est compatible aux changements de
base de la série principale non-ramifiée et est caractérisé par
cette propriété.

Le lemme fondamental pour le changement de base affirme que
l'intégrale d'une fonction sphérique $f$ sur une classe de
conjugaison tordue $\delta$ est égale à l'intégrale de ${\mathbf
b}(f)$ sur la norme $\gamma$ de $\delta$
$${\mathbf{TO}}_\delta(f)={\mathbf O}_\gamma({\mathbf b}(f)).$$

Ce lemme a été établi pour la fonction unité $f=1_{G(\calO_E)}$
par Kottwitz, puis pour toutes les fonctions sphériques par
Arthur-Clozel dans le cas $\GL_n$. Pour un groupe réductif
quelconque, le lemme fondamental pour le changement de base stable
a été démontré par Clozel et Labesse du moins dans le cas
$p$-adique. Leur méthode est basée sur l'usage de la formule des
traces et le cas d'unité établi par Kottwitz.

Ce lemme fondamental pour le changement de base est un outil
essentiel dans le comptage de points des variétés de Shimura ou
des variétés modulaires de Drinfeld et la comparaison avec la
formule des traces d'Arthur-Selberg. Le comptage naturel fait
apparaître des intégrales orbitales tordues et le lemme est
nécessaire pour les convertir en intégrales orbitales non tordues.

Dans ce travail, nous adoptons une stratégie inverse en partant du
comptage des points dans les variétés modulaires de Drinfeld pour
arriver au lemme fondamental. L'idée du départ est la suivante.
Dans les variétés de Shimura avec mauvaise réduction, le comptage
fait intervenir des fonctions de pondération provenant du modèle
local des singularités qui est une sorte de sous-variétés de
Schubert de variétés de drapeaux affines ou d'une Grassmannienne
affine \cite{Ra}. On a par ailleurs une interprétation géométrique
de l'homomorphisme de changement base ${\mathbf b}$ à l'aide des
faisceaux pervers sur la Grassmannienne affine \cite{N}. Il s'agit
donc de fabriquer un problème de module global tel que la
construction \cite{N} apparaisse comme son modèle local. On est
conduit ainsi à étudier un espace de module des chtoucas à
modifications multiples et l'action du groupe symétrique sur
celles-ci. Une variante globale du lemme fondamental pour le
changement de base prend la forme d'une comparaison du problème de
module de chtoucas à modifications multiples et de la restriction
à la Weil du problème de module des chtoucas à modifications
simples. Nous ne faisons pas appel au cas d'unité établi par
Kottwitz.

Signalons que E. Lau a indépendamment obtenu une démonstration de
ce lemme fondamental, également à l'aide des chtoucas à
modifications multiples mais par un argument différent.

Nous pensons que ce type de comparaison existe pour des chtoucas
pour d'autres groupes réductifs, pourrait produire de nouvelles
identités et établir certains cas de fonctorialité. Il faudra
toutefois résoudre avant de nombreux problèmes notamment le
comptage de $G$-chtoucas et la compactification des espaces de
modules de $G$-chtoucas. Dans ce travail, on se restreint à des
$\calD$-chtoucas où $\calD$ est une algèbre à division ayant
beaucoup de places totalement ramifiées ce qui assure la compacité
du problème de module. De plus, puisque le groupe $G=\calD^\times$
est une forme intérieure de $\GL_n$, le problème d'endoscopie
disparaît ce qui simplifie considérablement le comptage des
points.

Nous pensons aussi que ces constructions sont peut-être possibles
même pour les variétés de Shimura car les fonctions de changement
de bases ont aussi apparus sur les modèles locaux de
Pappas-Rapoport, des variétés de Shimura pour des groupes
unitaires en une place où le groupe unitaire est ramifié, voir
\cite{PaRa}.

Le résultat démontré dans cet article n'est pas nouveau, seule la
démarche l'est. Vu son caractère un peu expérimental, nous avons
voulu privilégier la clarté de l'exposé à la généralité du
résultat. On ne démontrera dans ce texte, le lemme fondamental que
pour les fonctions sphériques $f$ dont les transformés de Satake
sont des fractions symétriques de degré $0$. \footnote{On référera
à une version précédente de ce papier, disponible sur {\tt arXiv},
où ce lemme est démontré pour toutes les fonctions de Hecke. E.
Lau nous a informé d'une erreur dans cette ancienne version : le
problème de module de $\calD$-chtoucas n'est pas nécessairement
compact si on ne suppose pas qu'il y a des places totalement
ramifiées en grande quantité par rapport aux modifications. En
augmentant le nombre de places totalement ramifiées, la
généralisation à des fonctions sphériques de degré non nul reste
valable. Cette généralisation est abandonnée ici à cause de sa
longueur.}

Passons en revue les différents chapitres de ce texte. On rappelle
dans le chapitre 1 la définition des $\calD$-chtoucas à
modifications multiples, leur problème de module et son modèle
local des singularités. La notion de chtoucas avec modifications
multiples, est due à Drinfeld, a été reprise récemment par
Varshavsky. On donnera également un critère nécessaire pour que le
morphisme caracté\-ristique soit propre. Dans ce chapitre, on suit
de près le livre de Lafforgue \cite{LL}.

Sous l'hypothèse que le morphisme caractéristique est propre, les
images directes du complexe d'intersection sont des systèmes
locaux. On rappelle dans le chapitre 2 la définition de l'action
des correspondances de Hecke sur ces systèmes locaux. On énonce
dans 2.2 la description conjecturale de la somme alternée de ces
systèmes locaux en suivant les conjectures de Langlands. Nous nous
en servons comme un guide heuristique.

On énonce dans le paragraphe 3.3 le théorème principal de
l'article, suivi d'un argument heuristique basé sur la description
conjecturale de la somme alternée.

La démonstration de ce théorème est basée sur le comptage qui est
le sujet du chapitre 4. On y présente le comptage des chtoucas à
modifications multiples, en suivant le comptage de Kottwitz des
points des variétés de Shimura. Ce chapitre contient quelques
innovations techniques et devrait se généraliser aux $G$-chtoucas
pour un groupe réductif $G$ arbitraire.

On termine la démonstration du théorème 1 de 3.3 en 5.4 en
comptant les points au-dessus des points cycliques et en évoquant
le théorème de densité de Chebotarev. Le point un peu subtil de
cette démonstration est la description de ces points cycliques en
5.2.

En spécialisant à la petite diagonale, on obtient une égalité où
apparaît l'opérateur de changement de base ${\mathbf b}$. On en
déduit de cette égalité le lemme fondamental pour le changement de
base pour les fonctions sphériques de degré $0$.

Ce travail a été réalisé en grande partie durant une visite à
l'IHES de 2001 à 2003. Je remercie A. Genestier, R. Kottwitz, L.
Lafforgue, V. Lafforgue et G. Laumon pour leurs conseils et leurs
encouragements, tout particulièrement E. Lau pour sa vigilance.
J'ai été très influencé par certaines idées de Rapoport sur les
modèles locaux des variétés de Shimura. Je suis particulièrement
reconnaissant au referee de cet article pour sa lecture attentive
et pour sa grande patience.

\section*{Notation}
\setcounter{lemme}{0}

Soit $X$ une courbe projective lisse géométriquement connexe sur
un corps fini $\FF_q$. On note $F$ son corps des fonctions
rationnelles. Pour tout schéma $S$, on note $|S|$ l'ensemble des
points fermés de $S$. En particulier, pour $X$, l'ensemble $|X|$
est l'ensemble des valuations de $F$. Pour tout $x\in|X|$, on note
$F_x$ le complété de $F$, $\calO_x$ son anneau des entiers,
$\kappa(x)$ le corps résiduel de $\calO_x$. On note $\deg(x)$ le
degré de l'extension $[\kappa(x):\FF_q]$.

Soit $\bbA$ l'anneau des adèles de $F$ et $\calO_\bbA$ son anneau
des entiers. On a un homomorphisme degré $\bbA^\times\rta\ZZ$
défini par $$\deg(a_x)_{x\in|X|}=\sum_x \deg(x)x(a_x).$$ On fixera
une idèle $a\in\bbA^\times$ de degré $1$.

On va fixer une clôture algébrique $k$ de $\FF_q$. Si $S$ est un
sch\'ema sur $\FF_q$, on va noter $\bab S=S\otimes_{\FF_q} k$. En
particulier, on a $\bab X=X\otimes_{\FF_q}k$.

On fixera $D$ une algèbre à division centrale sur $F$, de
dimension $d^2$. Soit $\calD$ une $\calO_X$-Algèbre de fibre
générique $D$ telle que $\calD_x$ soit un ordre maximal de $D_x$
pour tout $x\in |X|$. Soit $X'$ le plus grand ouvert de $X$ constitué
des points $x\in|X|$ tels que $\calD_x$ est isomorphe à
$\gl_d(\calO_x)$.

Soit $G$ le schéma en groupes sur $X$ qui associe à tout schéma
affine $\Spec(R)$ au-dessus de $X$ le groupe
$(\calD\otimes_{\calO_X} R)^\times$. Puisque la fibre générique
$D$ de $\calD$ est une algèbre à division centrale sur $F$, le
groupe $G_F$ est anisotropique c'est-à-dire qu'il ne contient
aucun $F$-tore déployé non trivial. Au-dessus de $X'$, le schéma
en groupes $G$ est isomorphe à $X'\times\GL_d$ localement pour la
topologie de Zariski de $X'$. Sur $X$ tout entier, il est
isomorphe à $X\times \GL_d$ localement pour la topologie étale.

On utilisera souvent l'équivalence de Morita entre modules sur
l'algèbre des matrices $M_d(K)$ et espace vectoriel sur un corps
$K$. Ceci se fait comme habitude à l'aide d'un idempotent $e\in
M_d$ qu'on va choisir une fois pour toutes. Ce choix est b\'enin
dans la mesure o\`u il n'intervient que dans les discussions sur
les classes d'isomorphisme de modifications, qui finalement, ne
d\'ependent pas de ce choix.

\section{$\calD$-chtoucas à modifications multiples}

\subsection{Modifications et Grassmannienne affine}
\setcounter{lemme}{0}

Considérons le champ $\DFib$ des $\calD$-fibrés \emph{de rang
$1$}. Il associe à tout $\FF_q$-schéma la catégorie en groupoïdes
des $\calD\boxtimes\calO_S$-Modules à droite, localement libre de
rang $1$. Ceux-ci sont en particulier des $\calO_{X\times
S}$-Modules localement libres de rang $d^2$. Il est clair que
$\DFib$ est aussi le champ des $G$-torseurs localement triviaux
pour la topologie étale de $X$, parce que $G$ est le schéma en
groupes des automorphismes du $\calD$-module trivial $\calD$ et
que tous les autres objets de la catégorie $\DFib(k)$ sont
localement isomorphes à ce module trivial.

Soit $\bab T$ un sous-schéma fini de $\bab X'$. Soient $\calV$ et
$\calV'$ deux $\calD$-fibrés de rang $1$ sur $\bab X$. Une {\em
$\bab T$-modification} de $\calV'$ dans $\calV$ est un
isomorphisme entre les restrictions de $\calV$ et de $\calV'$ à
$\bab X -\bab T$
$$
t:\calV'|_{\bab X -\bab T} \isom\calV|_{\bab X -\bab T}.$$ Pour
alléger les notations, notons $\calV^{\bab T}$ et ${\calV'}^{\bab
T}$ les restrictions de $\calV$ et $\calV'$ à l'ouvert $\bab
X-\bab T$. Une $\bab T$-modification sera alors un isomorphisme
$t:{\calV'}^{\bab T}\isom \calV^{\bab T}$.

Rappelons la définition de {\em l'invariant} d'une $\bab
T$-modification. Soit $\bab x\in |\bab T|$ un point de $\bab T$.
Notons $\calO_{\bab x}$ le complété de $\calO_{\bab X}$ en
$\overline x$, et notons $F_{\bab x}$ son corps des fractions.
Notons $V$ et $V'$ les fibres génériques de $\calV$ et $\calV'$
respectivement. La $\bab T$-modification $t$ permet d'identifier
$V\isom V'$ ainsi que leurs complétés en $\bab x$ qu'on note
$V_{\bab x}\isom V'_{\bab x}$. Le complété en $\bab x$ de $\calV$
définit un $\calD_{\bab x}$-réseau de $V_{\bab x}$, celui de
$\calV'$ définit un $\calD_{\bab x}$-réseau de $V'_{\bab x}$.
Puisque $\bab T$ est contenu dans $\bab X'$ l'ouvert de
non-ramification de $D$, il existe un isomorphisme  $\calD_{\bab
x}\simeq M_d(\calO_{\bab x})$ qui induit un isomorphisme $D_{\bab
x}\simeq M_d(F_{\bab x})$. En utilisant l'idempotent standard
$e\in M_d(\calO_{\bab x})$, on obtient un $F_{\bab x}$-espace
vectoriel $e(V_{\bab x})=e(V'_{\bab x})$ de dimension $d$ avec
deux $\calO_{\bab x}$-réseaux $e(\calV_{\bab x})$ et
$e(\calV'_{\bab x})$. La position relative de ces deux réseaux est
donnée par une suite de $d$ entiers décroissants
$$
\lambda_{\bab x}=(\lambda^1\geq \lambda^2 \geq\cdots\geq
\lambda^d)$$ d'après le théorème des diviseurs élémentaires. On
pose ${\rm inv}_{\bab x}(t)=\lambda_{\bab x}$.

On peut donc associer à toute $\bab T$-modification $t:\calV^{\bab
T}\rta {\calV'}^{\bab T'}$ une fonction ${\rm inv}(t):|T|\rta
\ZZ^d_+$ à valeurs dans
$$\ZZ^d_+=\{\lambda=(\lambda^1,\ldots,\lambda^d)\in\ZZ^d\mid
\lambda^1\geq\cdots\geq\lambda^d\},$$ définie par $\bab x\mapsto
{\rm inv}_{\bab x}(t)$. Il sera assez commode d'écrire
formellement $\inv(t)$ sous la forme
$$\inv(t)=\sum_{\bab x\in\bab T}\inv_{\bab x}(t)\bab x.$$

En termes des $G$-torseurs, une $\bab T$-modification est un
isomorphisme entre les restrictions de deux $G$-torseurs $\calV$
et $\calV'$ à $\bab X-\bab T$. Dans un voisinage assez petit de
$\bab T$ contenu dans $\bab X'$, $G$ est isomorphe à $\GL_d$. Un
tel isomorphisme étant choisi, on peut associer à une $\bab
T$-modification de $G$-torseurs, une $\bab T$-modification de
fibrés vectoriels de rang $d$. L'invariant défini ci-dessus
coïncide avec l'invariant habituel associé à une $\bab
T$-modification entre deux fibrés vectoriels de rang $d$ définis
sur un voisinage arbitrairement petit de $\bab T$.

Nous allons noter $\mu$ et $\mu^\smv\in\ZZ^d_+$ les suites
d'entiers $$ \mu=(1,0,\ldots,0)\ \mbox{ et }\
\mu^\smv=(0,\ldots,0,-1). $$ Elles sont parmi les éléments
minimaux de $\ZZ^d_+$ par rapport à l'ordre partiel usuel, défini
comme suit : pour tous $\lambda,\lambda'\in\ZZ^d_+$, $\lambda\geq
\lambda'$ si et seulement si
$$
\begin{array}{rcl}
\lambda^1 & \geq &{\lambda'}^1 \\
& \ldots & \\
\lambda^1+\cdots+\lambda^{d-1} & \geq & {\lambda'}^1+\cdots+
{\lambda'}^{d-1} \\
\lambda^1+\cdots+\lambda^{d} & = & {\lambda'}^1+\cdots+
{\lambda'}^{d}.\\
\end{array}
$$
Nous aurons aussi besoin de la fonction croissante $\rho:\ZZ^d_+\rta
\mathbb N$ définie par
$$ \rho(\lambda)={1\over 2}[(d-1)\lambda^1+(d-3)
\lambda^2+ \cdots+(1-d)\lambda^d]
$$
pour les calculs de dimension. On a, en particulier,
$2\rho(\mu)=2\rho(\mu^\smv)=d-1$.

Les modifications de position $\mu$ et $\mu^\smv$ correspondent
bien aux notions usuelles de modifications supérieures et
inférieures de fibrés vectoriels de rang $d$. Soit $\bab u\in
X'(k)$ et prenons $\bab T= \{\bab u\}$. Une $\bab T$-modification
de $\calD$-fibrés $t:\calV^{\bab T} \rta {\calV'}^{\bab T}$
d'invariant $\mu$ se prolonge un une injection
$\calD\otimes_{\FF_q}k$-linéaire $\calV\rta \calV'$ dont le
conoyau est un $k$-espace vectoriel de rang $d$, supporté par
$\bab u$. Si de plus, on a choisi un voisinage assez petit de
$\bab u$ sur lequel il existe un isomorphisme entre $G$ et
$\GL_d$, une $\bab T$-modification de $\calD$-fibrés
$t:\calV^{\bab T} \rta {\calV'}^{\bab T}$ de position $\mu$
induit, sur ce voisinage, une injection de fibrés vectoriels de
rang $d$ dont le conoyau est un $k$-espace vectoriel de dimension
$1$, supporté par $\bab u$. Idem, les $\bab T$-modifications
d'invariant $\mu^\smv$ correspondent aux modifications inférieures
de fibrés vectoriels de rang $d$.

\begin{definition}
Soit $\lambda\in\ZZ^d_+$. Le champ $\DHck_\lambda$ associe à tout
$\FF_q$-schéma $S$ la catégorie en groupoïdes des données
$(u,\calV,\calV')$ comprenant :
\begin{itemize}

\item un point $u$ de $X'$ à valeurs dans $S$,

\item deux points $\calV,\calV'$ de $\DFib$ à valeurs dans $S$,

\item  un isomorphisme
$$t: \calV'|_{X\times S-\Gamma(u)}\isom
\calV|_{X\times S-\Gamma(u)}$$
entre les restrictions de $\calV$ et $\calV'$ au complémentaire
dans $X\times S$ du graphe $\Gamma(u)$ de $u:S\rta X$,
\end{itemize}
telles que pour tout $\bab s\in S(k)$, pour $\bab u=u(\bab s)$,
au-dessus de $\bab s$, la modification $t_{\bab s}:{\calV'}_{\bab
s}^{\{\bab u\}}\isom \calV_{\bab s}^{\{\bab u\}}$ a un invariant
$\inv_{\bab u}(t_{\bab s})\leq\lambda$.
\end{definition}

On va aussi considérer le sous-champ $\DHck'_{\lambda}$ de
$\DHck_{\lambda}$ classifiant les mêmes données, mais en imposant
une condition plus forte sur la modification $t$ : on demande
qu'en tout point géométrique $\bab s\in S(k)$, la $\{u(\bab
s)\}$-modification $t_{\bab s}$ a l'invariant ${\rm inv}_{u(\bab
s)}(t)=\lambda$.

Si $\lambda$ est un élément minimal de $\ZZ^d_+$, comme $\mu$ ou
$\mu^\smv$, il est clair que $\DHck'_\lambda=\DHck_\lambda$.

\begin{proposition} Le morphisme
$$\DHck_{\lambda}\rta X'\times \DFib$$ qui associe
$(u;\calV,\calV';t) \mapsto (u,\calV)$ est un morphisme
représentable et projectif. De plus, $\DHck'_{\lambda}$ est un
ouvert dense de $\DHck_{\lambda}$, lisse et de purement de
dimension relative $2\rho(\lambda)$ sur ${X'}\times \DFib$.

En particulier, si $\lambda$ est égal à $\mu$ où $\mu^\smv$, le
morphisme $\DHck_{\lambda}\rta X'\times \DFib$ est projectif,
lisse et purement de dimension relative $d-1$.
\end{proposition}

Pour un $\lambda$ non minuscule, $\DHck_{\lambda}$ n'est pas
nécessairement lisse au-dessus de $X'\times \DFib$. Toutefois, il
est lisse au-dessus de $\calQ_\lambda\times\DFib$ où
$\calQ_\lambda$ est une variante globale de la Grassmannienne
affine, introduite par Beilinson et Drinfeld \cite{BD}. Rappelons
la définition de $\calQ_\lambda$.

Soit $Q$ un $\calO_{\bab X}$-Module supporté par un point $u\in
\bab X$. Comme $\calO_u$-module, il est isomorphe \`a un module de
la forme
$$Q\simeq\calO_u/\varpi_u^{\lambda^1}\oplus\cdots\oplus\calO_u/
\varpi_u^{\lambda^{d'}}$$ où la suite d\'ecroissante d'entiers
naturels $\lambda^1\geq\cdots\geq\lambda^{d'}>0$ est bien défini.
Ici, $\varpi_u$ est un uniformisant de $\bab X$ en $u$. On ne va
considérer que les suites de longueur $d'\leq d$ lesquelles seront
complétées par nombre nécessaire de $0$ pour devenir un élément
$\lambda\in\NN^d_+$ c'est-à-dire
$$\lambda=(\lambda^1\geq\lambda^2\geq\cdots\geq\lambda^d\geq 0).$$
On notera alors ${\rm inv}(Q)=\lambda$.

\begin{definition}
Pour tout $\lambda\in\NN^d_+$, le champ $\calQ_\lambda$ qui
associe \`a tout $\FF_q$-sch\'ema $S$ la cat\'egorie en groupoïdes
des donn\'ees suivantes :
\begin{itemize}
\item un morphisme $u: S\rta X'$,

\item un $\calD_X\otimes\calO_S$-Module $Q$, $\calO_S$-plat,
support\'e par le graphe de $u$ tel que pour tout point
g\'eom\'etrique $s\in S(k)$, pour $u=u(s)\in \bab X'$, pour tout
isomorphisme $\calD_{u}\isom M_d(\calO_u)$, au-dessus de $s$,
$e(Q)$ est un module de torsion support\'e par $u$, d'invariant
$\inv(Q)\leq\lambda$. Ici, $e$ est l'idempotent standard de
$M_d(\calO_u)$ qu'on a transf\'er\'e \`a $\calD_u$.
\end{itemize}
\end{definition}

\begin{proposition}
Soit $\lambda\in\NN^d_+$. Soit
$$(u;\calV,\calV';t:\calV^{\{u\}}\isom{\calV'}^{\{u\}})$$ un
objet de $\DHck_\lambda(S)$. Alors $t$ se prolonge de fa\c con
unique en une injection $\calD_X\boxtimes \calO_S$-lin\'eaire
$t:\calV\hookrightarrow \calV'$. De plus, si $Q$ d\'esigne le
conoyau de cette injection, la paire $(u,Q)$ est un objet de
$\calQ_\lambda(S)$.

Le morphisme $$\DHck_\lambda\rta \calQ_\lambda\times\DFib$$ qui
associe \`a $(u;\calV,\calV';t)$ la paire $(u,Q)\in\calQ_\lambda$
et le $\calD$-fibr\'e $\calV\in\DFib$, est un morphisme lisse. Il
en est de même du morphisme $(u;\calV,\calV';t)\mapsto
(u;Q,\calV')$.
\end{proposition}

\dem Pour démontrer que le morphisme
$$f:(u;\calV,\calV';t)\mapsto (u;Q,\calV')$$
est lisse, il suffit de démontrer qu'il est le composé d'une
immersion ouverte suivie de la projection d'un fibré vectoriel sur
sa base. La modification étant locale dans un voisinage de $u$, on
peut remplacer un $\calD$-module par un fibré vectoriel de rang
$n$. Étant donné un fibré vectoriel de rang $n$ et un
$\calO$-module de torsion $Q$ de type $\lambda$, la fibre du
morphisme $f$ consiste en des application $\calO$-linéaire
surjective $\calV'\rta Q$. En omettant la surjectivité, la fibre
est clairement un espace vectoriel de dimension fixe. Quant \`a la
condition de surjectivité, c'est une condition ouverte. \hfill
$\square$

\bigskip
Il n'est pas difficile d'étendre la définition de $\calQ_\lambda$
et la proposition ci-dessus à un élément $\lambda\in\ZZ^d_+$
arbitraire. Soit $\lambda\in\ZZ^d_+$ et soit $N\in\NN$ un entier
assez grand tel que
$\lambda+N=(\lambda^1+N\geq\cdots\geq\lambda^d+N\geq 0)$. On va
poser formellement $\calQ_{\lambda}=\calQ_{\lambda+N}$. Ceci ne
dépend pas de l'entier $N$ pourvu que $\lambda_d+N\geq 0$. Si $t$
est une modification en $u$ de $\calV'$ dans $\calV$ de position
$\lambda$, $t$ induit une modification en $u$ de $\calV$ dans
$\calV'[N]$ de type $\lambda+N$. La proposition s'\'etend donc
sans difficult\'es au cas $\lambda\in\ZZ^d_+$.

L'énoncé suivant utilise le lien entre $\calQ_\lambda$ et les
cellules de la Grassmannienne affine introduite par Lusztig
\cite{Lus}, voir \cite{G}, \cite{MV} et \cite{NP}.

\begin{proposition}
Il existe une variété projective $\tilde\calQ_\lambda$ telle que
localement pour la topologie étale de $X'$, il existe de
$X'$-morphismes lisses surjectifs de $X'\times\tilde\calQ_\lambda$ sur
$\calQ_\lambda$.
\end{proposition}

\dem Suivant Mirkovic et Vilonen \cite{MV}, on considère le foncteur
$\tilde\calQ_{\lambda,X'}$ qui associe à tout $\FF_q$-schéma $S$
l'ensemble des données $(u,\calV,t)$
\begin{itemize}
\item $u$ est un point de $X'$ à valeur dans $S$, \item $\calV$
est un point de $\DFib$ à valeur dans $S$, \item
$t:\calD^{\{u\}}\rta \calV^{\{u\}}$ est une modification en $u$
d'invariant $\inv_u(t)\leq\lambda$, du $\calD$-fibré trivial
$\calD$ dans $\calV$.
\end{itemize}
Ce foncteur est représentable par un schéma projectif au-dessus
$X'$. D'après la proposition précédente, on sait que le morphisme
$\tilde\calQ_{\lambda,X'}\rta \calQ_\lambda$ est un morphisme
lisse. Il n'est pas difficile de vérifier qu'il est surjectif.

Si $Y'\rta X'$ est un revêtement étale, il est clair que
$$\tilde\calQ_{\lambda,Y'}=\tilde\calQ_{\lambda,X'}\times_{X'}
Y'.$$ Localement pour la topologie étale, on peut supposer que
$Y'$ est isomorphe à la droite affine. Puisque $Y'\simeq \bbA^1$,
on a une action par translation de $\bbA^1$ sur
$\tilde\calQ_{\lambda,X}$ qui permet d'obtenir un isomorphisme
$\tilde\calQ_{\lambda,Y'}\simeq Y'\times \tilde\calQ_\lambda$ où
$\tilde\calQ_\lambda$ est une fibre quelconque de
$\tilde\calQ_{\lambda,Y'}$ au-dessus de $Y'$. \findem

\begin{corollaire} Soit $\calA_\lambda$ le complexe d'intersection
de $\calQ_\lambda$. Il satisfait deux propriétés :
\begin{enumerate}
\item
 Il est localement acyclique par rapport au morphisme
$\calQ_\lambda\rta X'$.

\item Ses restrictions aux fibres de $\calQ_\lambda$ sur $X'$
sont, à décalage près, des faisceaux pervers irréductibles.
\end{enumerate}
\end{corollaire}

\dem Localement pour la topologie étale de $X'$, il existe un
$X'$-morphisme lisse surjectif $X'\times \tilde \calQ_\lambda \rta
\calQ_\lambda$. On se ramène donc à démontrer l'énoncé pour
$X'\times \tilde\calQ_\lambda$ où il se déduit de la formule de
Künneth. \findem

\bigskip
D'après un théorème de Lusztig, les faisceaux pervers
$\calA_\lambda$ sur $\calQ_\lambda$ réalisent via le dictionnaire
faisceaux-fonctions une base de l'algèbre de Hecke ayant des
propriétés remarquables. Pour préciser cet énoncé dont on va s'en
servir, rappelons la définition des algèbres de Hecke.

Soit $x$ un point fermé de $X'$. On définit
$\calH_x$ comme l'algèbre des fonctions à valeurs dans $\QQ_\ell$, à support
compact dans $G(F_x)$ qui sont invariantes à gauche et à droite par
$G(\calO_x)$. La structure d'algèbre de $\calH_x$ est définie par
le produit de convolution usuelle
$$(f*f')(g)=\int_{G(F_x)}f(h)f'(h^{-1}g)\,dg$$
en intégrant par rapport à la mesure de Haar $dg$ sur $G(F_x)$, normalisée
par ${\rm vol}(G(\calO_x))=1$.

Pour tout $x\in |X'|$, $G$ est isomorphe à $\GL_d$ dans un
voisinage de Zariski de $x$ de sorte que $\calH_x$ est isomorphe à
l'algèbre de Hecke de $\GL_d$. L'algèbre de Hecke de $\GL_d$ admet
donc une base évidente $(\phi_{\lambda})$, paramétrée par les
$\lambda\in\ZZ^d_+$ ; la fonction $\phi_{\lambda}$ étant la
fonction caractéristique de la double-classe de $\lambda(x)$
$$\GL_d(\calO_x){\rm diag}(\varpi_x^{\lambda_1},\ldots,
\varpi_x^{\lambda_n}) \GL_d(\calO_x)$$
où $\varpi_x$ est une uniformisante de $F_x$. Transférée à $G$, on
obtient une base de $\calH_x$. Il n'est pas difficile de montrer
que la base ainsi obtenue, ne dépend ni du choix de
l'isomorphisme entre $G_x$ et $\GL_d$, ni du choix de l'uniformisante
$\varpi_x$.

L'isomorphisme de Satake permet d'identifier canoniquement
l'algèbre $\calH_x$ avec l'algèbre $\QQ_\ell[\hat G]^{{\rm
Ad}(\hat G)}$ des fonctions polynômiales
$$\hat f:\hat G(\QQ_\ell)\rta\QQ_\ell$$
qui sont invariantes par l'action adjointe de $\hat G(\QQ_\ell)$.
Dans la situation présente, le groupe dual de Langlands $\hat G$
est aussi égal à $\GL_d$. Pour tout $\lambda\in\ZZ^d_+$, notons
$\hat\psi_\lambda$ la fonction définie par la trace de tout
élément $\hat g\in\hat G(\QQ_\ell)$ sur la représentation
irréductible de $\hat G$ de plus haut poids $\lambda$. Ces
fonctions forment une base de $\QQ_\ell[\hat G]^{{\rm Ad}(\hat
G)}$ lorsque $\lambda$ parcourt $\ZZ^d_+$. Le transformé de Satake
$\hat\phi_\lambda$ de la fonction double-classe $\phi_\lambda$
diffère en général de $\hat\psi_\lambda$. Le faisceau pervers
$\calA_\lambda$ tient compte précisément de cette différence.

Supposons $\deg(x)=s$. Les $k$-points du champ $\calQ_{\lambda}$
au-dessus de $x$, en nombre fini, tous définis sur $\FF_{q^s}$,
correspondent bijectivement aux éléments $\alpha\in\ZZ^d_+$ tels
que $\alpha\leq\lambda$. On les note $\alpha(x)\in
\calQ_\lambda(x)$. Soit $P_{\lambda,\alpha}(q^s)$ la trace de
$\Frob_{q^s}$ sur la fibre de $\calA_\lambda$ au-dessus de
$\alpha(x)$.

\begin{theoreme}[Lusztig-Kato]
La fonction $$\psi_\lambda=\sum_{\alpha\leq \lambda}
P_{\lambda,\alpha} (q^s)\phi_{\alpha}$$ dans $\calH_x$, a pour
transformé de Satake la fonction $\hat \psi_\lambda$ définie par
la trace sur la représentation irréductible de plus haut poids
$\lambda$ de $\hat G(\QQ_\ell)$.
\end{theoreme}

On renvoie à \cite{Lus} pour la démonstration de ce théorème et à
\cite{G} et \cite{MV} pour son interprétation géométrique.

\subsection{Modifications it\'er\'ees et propriété de factorisation}
\setcounter{lemme}{0}

On aura aussi besoin de considérer des suites de modifications
avec les invariants bornés. Soient $n$ un entier naturel et
$\und\lambda=(\lambda_1,\ldots, \lambda_n)$ un $n$-uple d'éléments
de $\ZZ^d_+$.

\begin{definition}
Le champ $\DHck_{\und\lambda}$ associe à tout $\FF_q$-schéma $S$
la catégorie en groupoïdes des données
$((u_j)_{j=1}^n,(\calV_j)_{j=0}^n,(t_j)_{j=1}^n)$ comprenant
\begin{itemize}

\item $n$ points $u_1,\ldots,u_n$ de $X'$ à valeurs dans $S$,

\item $n+1$ points $\calV_0,\ldots,\calV_n$ de $\DFib$ à valeurs dans
$S$,

\item pour tous $j=1,\ldots,n$, un isomorphisme
$$t_j: \calV_{j}|_{X\times S-\Gamma(u_j)}\isom
\calV_{j-1}|_{X\times S-\Gamma(u_j)}$$ entre les restrictions de
$\calV_{j-1}$ et $\calV_j$ au complémentaire dans $X\times S$ du
graphe $\Gamma(u_i)$ de $u_j:S\rta X$.
\end{itemize}
telles que pour tout $\bab s\in S(k)$, pour $\bab u_j=u_j(\bab
s)\in \bab X'$,  au-dessus de $\bab s$, la modification
$t_j:\calV_{j}^{\{\bab u_j\}}\isom \calV_{j-1}^{\{\bab u_j\}}$ a
l'invariant $\inv_{\bab u_j}(t_j) \leq \lambda_j$.
\end{definition}

Il est clair que $\DCht_\ul$ peut être construit au moyen de {\em
produit fibré successif}. Pour tout $j=1,\ldots,n$, si on note
$p_j$ et $p'_j$ les deux projections de $\DCht_{\lambda_j}$ sur
$\DFib$ qui envoie $(u,\calV,\calV',t)$ sur $\calV$ et $\calV'$
respectivement, alors on a
$$\DCht_\ul=\DCht_{\lambda_1}\times_{p'_1,\DFib,p_2}\cdots
\times_{p'_{n-1},\DFib,p_n}\DCht_{\lambda_n}.$$

On va aussi considérer le sous-champ $\DHck'_{\und\lambda}$ de
$\DHck_{\und\lambda}$ classifiant les mêmes données mais en
imposant une condition plus forte sur les modifications $t_i$ en
demandant que $\inv_{\bab u_j}(t_j)=\lambda_j$.

\begin{proposition} Soit $\und\lambda=(\lambda_1,\ldots,\lambda_n)$
une suite de $n$ éléments de $\ZZ^d_+$. Le morphisme
$$\DHck_{\und\lambda}\rta{X'}^n\times \DFib$$ qui associe
$((u_j)_{j=1}^n,(\calV_j)_{j=0}^n,(t_j)_{j=1}^n) \mapsto
((u_j),\calV_0)$ est un morphisme représen\-table et projectif. De
plus, $\DHck'_{\und\lambda}$ est un ouvert dense de
$\DHck_{\und\lambda}$, lisse et de purement de dimension relative
$\sum_{j=1}^n 2\rho(\lambda_j)$ sur ${X'}^n\times \DFib$.
\end{proposition}

\dem On se ramène au cas d'une seule modification au moyen
du produit fibré successif, mentionn\'e ci-dessus.
\findem

\bigskip
Sur l'ouvert complémentaire de la réunion des diagonales de
${X'}^n$, on aura un vrai produit fibré.

\begin{proposition}
Au-dessus de l'ouvert $U$  complémentaire de la réunion des
diagonales dans ${X'}^n$, $\DHck_{\und\lambda}$ est canoniquement
isomorphe au produit fibré des morphismes $p_j:\DHck_{\lambda_j}
\rta \DFib$ pour $j=1,\ldots,n$.
\end{proposition}

\dem Il s'agit d'une propriété bien connue des modification de
fibrés vectoriels. Soit $x_1,x_2$ deux points distincts de $\bab
X$. Fixer un fibré vectoriel $\calV$ sur $\bab X$. Alors donner
une "grande" modification
$$t:{\calV'}^{\{x_1,x_2\}}\isom \calV^{\{x_1,x_2\}}$$
telle que $\inv_{x_j}(t)=\lambda_j$ est équivalente à donner deux
modifications indé\-pen\-dantes
$$t'_1:{\calV'}_1^{x_1}\isom \calV^{x_1}\mbox{ et }
t'_2:{\calV'}_2^{x_2}\isom \calV^{x_2}$$
telles que $\inv_{x_j}(t'_j)=\lambda_j$, ou encore équivalente
à deux modifications "en itération"
$$t_2:\calV_2^{x_2}\rta\calV_1^{x_2}\mbox{ et }
t_1:\calV_1^{x_1} \rta \calV^{x_1}$$
avec les mêmes invariants.
L'adaptation est évidente au cas de $\calD$-fibrés si on suppose
en plus que $x_1,x_2$ sont dans $X'$. La proposition s'en déduit.
\findem

\bigskip
Tout comme dans le cas d'une seule modification, on veut
construire un objet "minimal" qui capture toutes les
singularit\'es du morphisme $\DHck_{\und\lambda}\rta {X'}^n\times
\DFib$.

\begin{definition}
Supposons que $\lambda_1,\ldots,\lambda_n\in\NN^d_+$. Le champ
$\calQ_{\und\lambda}$ qui associe \`a tout $\FF_q$-sch\'ema la
cat\'egorie en groupoïdes des donn\'ees suivantes :

- $n$ points $u_1,\ldots,u_n$ de $X'$ \`a valeurs dans $S$,

- un drapeau de $\calD_X\boxtimes\calO_S$-Modules de torsion,
$\calO_S$-plats
$$0=Q_0\subset Q_1\subset Q_2\subset \cdots\subset Q_n$$
tel que $(u_j,Q_i/Q_{j-1})$ est un point de $\calQ_{\lambda_j}$
\`a valeur dans $S$.
\end{definition}

\begin{proposition}
Soit $\lambda\in\NN^d_+$. Soit
$$((u_j)_{j=1}^n;(\calV_j)_{j=0}^n,\calV_{j}^{\{u_j\}}\hfld{t_j}
{}\calV_{j-1}^{\{u_j\}})$$ un objet de $\DHck_{\und\lambda}(S)$.
Alors $t_j$ se prolonge de fa\c con unique en une injection
$\calD_X\boxtimes \calO_S$-lin\'eaire
$t_j:\calV_{j}\hookrightarrow \calV_{j-1}$. De plus, si $Q_j$
d\'esigne le conoyau de $$t_j\circ\cdots\circ t_1:\calV_j\rta
\calV_0$$ alors $((u_j)_{j=1}^n,Q_1\subset\cdots\subset Q_n)$ est
un objet de $\calQ_{\und\lambda}(S)$.

Le morphisme $$\DHck_{\und\lambda}\rta
\calQ_{\und\lambda}\times\DFib$$ qui fait associer \`a
$((u_j);(\calV_j),(t_j))$ la paire $((u_j),(Q_j))\in\calQ_\ul$
comme ci-dessus et le $\calD$-fibr\'e $\calV_0\in\DFib$, est un
morphisme lisse.
\end{proposition}

\dem On se ram\`ene par r\'ecurence au cas $n=1$, donc \`a la
proposition 2 de 1.1, en utilisant la description $\DHck_\ul$
comme un "produit fibré successif", qui se trouve après la
définition 1. \findem

\bigskip
La définition 4 et la proposition 5 se généralisent sans
difficulté à tous $\lambda_1,\ldots,\lambda_n\in\ZZ^d_+$.

\begin{proposition}
Le complexe d'intersection  $\calA_{\und\lambda}$ de
$\calQ_{\und\lambda}$ a les propriétés suivantes.
\begin{enumerate}
\item
Il est localement acyclique par rapport au morphisme
$\calQ_{\und\lambda} \rta X^n$.

\item Ses restrictions aux fibres de ce morphisme sont
des faisceaux pervers irréductibles.
\end{enumerate}
\end{proposition}

\dem Essentiellement la m\^eme démonstration que 1.1.6 marche.
\findem

\bigskip
On a de plus une propriété de factorisation.

\begin{proposition}
Au-dessus de l'ouvert $U$ complémentaire de la réunion des
diagonales dans ${X'}^n$ , $\calQ_{\ul}$ est isomorphe à
$\calQ_{\lambda_1}\times\cdots \times\calQ_{\lambda_n}$. Via cet
isomorphisme, $\calA_{\ul}$ s'identifie à
$\calA_{\lambda_1}\boxtimes\cdots\boxtimes \calA_{\lambda_n}$.
\end{proposition}

\dem Il suffit de combiner la propriété de factorisation avec la
formule de Künneth. \findem

\subsection{$\calD$-chtoucas à modifications multiples}
\setcounter{lemme}{0}

La notion de chtoucas à modifications multiples est due à Drinfeld
\cite{D} et est récemment reprise dans un article de Varshavsky
\cite{Var}.

Soit $\und\lambda=(\lambda_1,\ldots,\lambda_n)$ une
suite de $n$ éléments de $\ZZ^d_+$. On définit le champ
$\DCht_{\und\lambda}$ en formant le produit cartésien
$$
\renewcommand\arraystretch{1.5}
\begin{array}{ccc}
\DCht_{\und\lambda} & \hfld{}{} & \DFib \\
\vfld{}{} & & \vfld{}{} \\
\DHck_{\und\lambda} & \hfld{}{} & \DFib\times \DFib
\end{array}
$$
où la flèche horizontale de bas envoie la donnée de modifications
itérées $((u_j)_{j=1}^n,(\calV_j)_{j=0}^n,(t_j)_{j=1}^n)$ sur
$(\calV_0,\calV_n)$ et la flèche verticale de droite envoie
$\calV$ sur $(\,^\sigma\calV,\calV)$. Puisque fibres par fibres
au-dessus de $S$, $^\sigma\calV$ et $\calV$ sont les fibrés sur
$X$ ayant le même degré, pour que $\DCht_{\und\lambda}$ ne soit
pas vide, il est nécessaire que $\sum_{i=1}^n \sum_{j=1}^d
\lambda_i^j=0$.

On peut expliciter la d\'efinition de $\DCht_{\und\lambda}$
comme suit.

\begin{definition}
Le champ $\DCht_\ul$ associe \`a tout $\FF_q$-sch\'ema la
cat\'egorie en groupoïdes des donn\'ees
$((u_j)_{j=1}^n,(\calV_j)_{j=0}^n, (t_j)_{j=1}^n,\epsilon)$
comprenant :

\begin{itemize}

\item $n$ points $u_1,\ldots,u_n$ de $X'$ \`a valeurs dans $S$,

\item $n+1$ points $\calV_0,\ldots,\calV_n$ de $\DFib$ \`a
valeurs dans $S$,

\item pour tout $j=1,\ldots,n$, une modification
$t_j:\calV_{j}^{\{u_j\}}\isom \calV_{j-1}^{\{u_j\}}$ d'invariant
$\inv_{u_j}(t_j)\leq \lambda_j$,

\item un isomorphisme $\epsilon :\, ^\sigma \calV_0 \isom \calV_n$.
\end{itemize}
\end{definition}

Selon une suggestion de V. Lafforgue, on peut appeler ces objets
{\em les chtoucas \`a mille pattes}. Dans le cas particulier des
chtoucas à deux pattes, et pour
$$\lambda_1=\mu=(1,0,\ldots,0)\mbox{ et }\lambda_2=\mu^\sv=(0,\ldots,0,-1)$$
on retrouve bien la notion de $\calD$-chtoucas usuels de rang $1$
tels qu'ils sont étudiés dans \cite{LL} chapitre 4.

Si les points $u_1,\ldots,u_n \in X'(k)$ sont deux à deux
distincts, la liste des données constituantes d'un
$\calD$-chtouca, ayant ses modifications en $u_1,\ldots,u_n$, peut
se réduire à une seule modification
$$t:\,^\sigma \calV_0|_{X-\{u_1,\ldots,u_n\}}\isom
\calV_0|_{X-\{u_1,\ldots,u_n\}}.$$ En effet la "grande"
modification $t=\epsilon\circ t_n\circ\cdots\circ t_1$ détermine
alors chacune des composantes $\epsilon,t_n,\ldots,t_1$ d'après la
proposition 3 de 1.2.

$\DCht_\ul$ a un nombre infini de composantes connexes paramétrées
par le degré de $\calV_0$.  Soit $a\in \bbA_F^\times$ un idèle de
degré non nul. On pose
$$\calS_\ul:=\DCht_\ul/a^\ZZ.$$
où le quotient par $a^\ZZ$, consiste à ajouter formellement dans
la catégorie $\DCht_\ul(S)$ un isomorphisme entre tout objet avec
l'objet qui s'en déduit par tensorisation avec le fibré en droite
$\calL(a)$ associé à l'idèle $a$. Le morphisme évident
$\DCht_\ul\rta {X'}^n$ induit un morphisme
$$c_\ul:\calS_\ul\rta {X'}^n$$
qu'on appellera {\em le morphisme caractéristique}.

\subsection{Modèles locaux}
\setcounter{lemme}{0}

Le champ des chtoucas usuels $\DCht_{\mu,\mu^\sv}$ est lisse
au-dessus de $X'\times X'$. Ce ne sera plus le cas pour
$\DCht_{\und\lambda}$ pour les suites $\und\lambda=
(\lambda_1,\ldots,\lambda_n)$ dont, au moins, une des composantes
n'est pas minuscule. On peut toutefois montrer que toutes les
singularités de $\DCht_\ul$ apparaissent déjà dans le $\calQ_\ul$
de sorte qu'on peut l'appeler le modèle local de $\DCht_\ul$ par
analogie avec les mod\`eles locaux de vari\'et\'es de Shimura avec
mauvaise r\'eduction. La th\'eorie de d\'eformation de
vari\'et\'es ab\'eliennes de Grothendieck-Messing est ici
remplacée par un lemme de type transversalit\'e du graphe de
l'endomorphisme de Frobenius avec la diagonale, dont voici
l'énoncé précis, voir \cite{LL} section I.2, proposition 1.

\begin{lemme} Considérons un diagramme de champs sur $\FF_q$
$$\renewcommand\arraystretch{1.5}
\begin{array}{ccc}
\calW & \hfld{}{} & \calU \\
\vfld{}{} &\square & \vfld{}{(\Frob_\calU,\id)} \\
\calM & \hfld{}{(\alpha,\beta)} & \calU\times \calU\\
\vfld{}{} && \\
Y & &\\
\end{array}$$
où $Y$ est un $\FF_q$-schéma, où $\calU$ est algébrique et
localement de type fini sur $\FF_q$, $\calM$ est algébrique et
localement de type fini sur $Y$, le morphisme
$(\pi,\alpha):\calM\rta Y\times\calU$ est représentable et où le
carré est $2$-cartésien.

Alors $\calW$ est algébrique et localement de type fini sur $Y$ et
le morphisme diagonal $\calW\rta \calW\times\calW$, qui est
représentable, séparé et de type fini, est partout non-ramifié et
en particulier, quasi-fini.

Si on suppose de plus que $\calU$ est lisse sur $\FF_q$ et que le
morphisme $$(\pi,\alpha):\calM\rta Y\times\calU$$ est lisse et
purement de dimension relative $n$, alors le morphisme $\calW\rta
Y$ est aussi lisse et purement de dimension relative $n$.
\end{lemme}

Soit $\und\lambda=(\lambda_1,\ldots,\lambda_n)\in(\ZZ^d_+)^n$. En
appliquant ce lemme à notre situation où $Y=\calQ_{\und\lambda}$,
$\calU=\DFib$, et où $\calM=\DHck_{\und\lambda}$, on trouve que le
champ $\DCht_{\und\lambda}$ est algébrique, localement de type
fini et que le morphisme
$$\DCht_{\und\lambda}\rta \calQ_{\und\lambda}$$
est un morphisme lisse. On en déduit un morphisme lisse
$$f_{\und\lambda}:\calS_\ul \rta \calQ_{\und\lambda}.$$
L'image inverse $\calF_{\und\lambda}=f_{\und\lambda}^*
\calA_{\und\lambda}$ est alors le complexe d'intersection de
$\calS_\ul$.

\begin{corollaire}
Le complexe d'intersection  $\calS_\ul$ de $\DCht_\ul$
a les propriétés suivantes :

\begin{enumerate}
\item
Il est localement acyclique par rapport au morphisme
caract\'eristique $c_\ul:\calS_\ul \rta {X'}^n$.
\item
Ses restrictions aux fibres de ce morphisme sont à décalage
près, des faisceaux pervers irréductibles.
\end{enumerate}
\end{corollaire}

\subsection{Structure de niveau}
\setcounter{lemme}{0}

Soit $I$ un sous-schéma fermé de $X$. On note $\calD_I$ la
restriction de $\calD$ à $I$.

Soit $\tilde\calV=((u_i),(\calV_i),(t_i),\epsilon)$ un
$\calD$-chtouca \`a valeur dans $S$ avec modi\-fi\-cation
d'invariants born\'es par $\lambda_1,\ldots,\lambda_n$ en
$u_1,\ldots,u_n\in X'(S)$, ou autrement dit, un point de
$\DCht_{\und\lambda}$ \`a valeur dans $S$. Supposons en plus que
$u_1,\ldots,u_n$ \'evitent $I$ c'est-\`a-dire $u_1,\ldots,u_n\in
(X'-I)(S)$. On peut alors d\'efinir ce qu'est une $I$-structure de
niveau de $\tilde\calV$.

Puisque les $u_i$ \'evitent $I$, les modifications $t_i$
d\'efinissent des isomorphismes entre les restrictions des
$\calV_i$ \`a $I\times S$ qu'on notera simplement $\calV_I$; c'est
un $\calD_I\otimes\calO_S$-module localement libre de rang $1$.
L'isomorphisme $\epsilon:\,^\sigma \calV_0\isom \calV_n$ induit un
isomorphisme $\epsilon_I:\,^\sigma\calV_I\isom \calV_I$.

\begin{definition}
Une $I$-structure de niveau de $\tilde\calV$ est la donn\'ee
d'un isomorphisme de $\calD_I\boxtimes\calO_S$-Modules
$$\iota:\calD_I\boxtimes\calO_S\rta \calV_I$$
tel que le diagramme
$$\renewcommand\arraystretch{1.5}
\begin{array}{ccc}
^\sigma(\calD_I\boxtimes \calO_S) & \hfld{^\sigma \iota}{}
&^\sigma \calV_I \\
\vfld{}{} &&\vfld{}{\epsilon_I}\\
\calD_I\boxtimes \calO_S & \hfld{}{\iota} & \calV_I \\
\end{array}$$
est commutatif. Ici, l'isomorphisme vertical gauche est celui qui
se déduit de l'isomorphisme canonique $\Frob_S^*\calO_S\isom
\calO_S$.
\end{definition}

On utilise une notation un peu abusive $\DICht_\ul$ pour désigner
le champ des $\calD$-chtoucas de modification en dehors de $I$ de
position born\'ee par $\ul$, munis d'une $I$-structure de niveau.
Suivant Lafforgue, on a un diagramme 2-cartésien
$$\renewcommand\arraystretch{1.5}
\begin{array}{rcc}
\DICht_\ul & \hfld{}{} & \Spec(\FF_q) \\
\vfld{}{} &\square & \vfld{}{} \\
(X'-I)^n\times_{{X'}^n}\DCht_\ul  & \hfld{}{} & \DITrv\\
\end{array}$$
où $\DITrv$ est le champ en groupoïdes qui associe à chaque
$\FF_q$-schéma $S$ la catégorie des
$\calD_I\boxtimes\calO_S$-Modules $\calF$ localement libre de rang
$1$ muni d'un isomorphism $\phi:\,^\sigma\cal F \isom\calF$. La
flèche horizontale en bas est définie par $$\tilde \calV\mapsto
(\epsilon_I:\,^\sigma\calV_I\isom \calV_I).$$ La flèche verticale
à droite est donnée par l'objet évident $\calF=\calD_I$ sur
$\Spec(\FF_q)$. Rappelons le théorème 2 de \cite{LL} section I.3.

\begin{proposition}
La flèche verticale de droite dans le carré ci-dessus est un
$\calD_I^\times$-torseur par rapport à l'action triviale de
$\calD_I^\times$ sur $\Spec(\FF_q)$. La flèche verticale de gauche
du carré est un morphisme représentable fini étale et galoisien de
groupe de Galois $\calD_I^\times$.
\end{proposition}

Posons $\calS^I_\ul:=\DCht^I_\ul/a^\ZZ$.
On obtient par composition un morphisme lisse
$f^I_\ul:\calS^I_\ul\rta\calQ_\ul$ et un morphisme
$$c_\ul^I:\calS^I_\ul \rta (X'-I)^n$$ qu'on appelle le morphisme
caractéristique de $\calS_\ul$. Nous notons $\calF^I_\ul$ l'image
r\'eciproque $(f^I_\ul)^* \calA_\ul$.

\begin{corollaire}
Le complexe d'intersection  $\calF^I_\ul$ de $\calS_\ul$ a les
propriétés suivantes.
\begin{enumerate}
\item
Il est localement acyclique par rapport au morphisme
caract\'eristique $c^I_\ul:\calS^I_\ul \rta (X'-I)^n$.
\item
Ses restrictions aux fibres de ce morphisme sont à décalage
près, des faisceaux pervers irréductibles.
\end{enumerate}
\end{corollaire}

Soit $v$ un point fermé de $I$ et soit $K^I_v$ le sous-groupe
compact de $G(F_v)$ défini comme le noyau de l'homomorphisme
naturel $[G(\calO_v)\rta G(\calO_I)]$. La donnée de $K^I_v$
définit un schéma en groupe lisse $\calG^I_v$ sur $\calO_v$ de
fibres connexes, de fibre générique $G_{F_v}$ et dont les points
entiers sont
$$\calG^I_v(\calO_v)=K^I_v.$$
Soit $\tilde V$ un $\calD$-chtouca avec une $I$-structure de
niveau défini sur $\bab k$. La structure de niveau induit un
$\calG^I_v\hat\otimes_{\FF_q} k$-torseur au-dessus
$\calO_v\hat\otimes_{\FF_q}k$. Sous l'hypothèse que les pôles et
zéros du chtouca évite $I$, l'isomorphisme $t:\,^\sigma
\calV^T\isom \calV$ permet de descendre le
$\calG^I_v\hat\otimes_{\FF_q} k$-torseur en un $\calG^I_v$-torseur
sur $\calO_v$ puisque $\calG^I_v$ est un schéma en groupes de
fibres connexes.

\subsection{Sur la propreté du morphisme caractéristique}
\setcounter{lemme}{0}

Comme nous a fait remarquer E. Lau, le morphisme caractéristique
$$c_\ul^I:\calS^I_\ul \rta (X'-I)^n$$
n'est pas propre en général. Il est néanmoins propre si on suppose
en plus que l'algèbre à division $\calD$ a un nombre de places
totalement ramifiées plus grand que le nombre de dégénérateurs.
Dans cette section, on va préciser cet énoncé.

Notons toutefois que l'argument principal de l'article ne
nécessite pas le propriété mais seulement de l'acyclicité locale
de la cohomologie. Il est donc probable qu'il puisse être utilisé
aussi dans des situations non propres du moment qu'on peut
construire une bonne compactification analogue à celle de
Lafforgue dans le cas $\GL_n$. Pour le lemme fondamental pour le
changement de base pour $\GL_n$, on peut se restreindre aux cas
des algèbres à division ayant des places totalement ramifiées en
quantité.

Tout élément $\lambda\in\ZZ^d_+$ peut être écrit de façon unique
sous la forme
$$\lambda=\lambda^++\lambda^-$$
où $\lambda_+\in\NN^d_+$ et $\lambda_-\in(-\NN)^d_+$ c'est-à-dire
$$\lambda^+=(\lambda_+^1\geq\cdots\geq\lambda_+^d\geq 0)$$
et
$$\lambda_-=(0\geq\lambda_-^1\geq\dots\geq\lambda_-^d)$$
telle que pour tous $i\in\{1,\ldots,d\}$, on a $\lambda^i$ est ou
bien égal à $\lambda_+^i$ ou bien égal à $\lambda_-^i$. On note
alors
$$||\lambda||:={\rm max}(|\lambda_+|,|\lambda_-|).
$$

\begin{lemme} Soient $\calV$ et $\calV'$ deux fibrés vectoriels de
rang $d$ sur $\bab X_u$ le complété de $\bab X$ en un point
$u\in\bab X$ et $t:\calV^{{u}}\isom \calV'^{{u}}$ une modification
d'invariant $\inv_u(t)\leq \lambda$. Alors $t$ se prolonge en un
morphisme de $\calO_{\bab X_u}$-Modules
$$t:\calV\rta \calV'(||\lambda||u).$$
\end{lemme}

\dem C'est la formule d'inversion des matrices de Cramer. \findem

\begin{proposition} Soit $\ul=(\lambda_1,\ldots,\lambda_n)$ avec
$\lambda_i\in\ZZ^d_+$. Supposons que le nombre de places
totalement ramifiées de $\calD$ est supérieur ou égal à
$$d^2(||\lambda_1||+\cdots+||\lambda_n||).$$
Alors, le morphisme $c_\ul^I:\calS^I_\ul \rta (X'-I)^n$ est un
morphisme propre.
\end{proposition}

\dem Nous allons reprendre l'argument de \cite{LL} IV.1 avec des
modifications appropriées.

Un $\calD$-chtouca de rang $1$ est automatiquement irréductible
est donc $\alpha$-semi-stable pour tout $\alpha\in ]0,1[$ de
sorte que le champ $\DICht_\ul$ est de type fini, voir \cite{LL},
II.2 théorème 8. Il nous suffira donc pour démontrer la propreté
de $c_\ul^I$ de vérifier le critère valuatif de propreté pour ce
morphisme.

Soient $\calO$ un anneau de valuation discrète contenant $\FF_q$,
$K$ son corps des fractions, $\kappa$ son corps résiduel. Notons
$\calX=X\otimes_{\FF_q}\calO$ qu'on voit comme une courbe fibrée
sur $\Spec(\calO)$. Notons $\eta_\kappa$ le point générique de la
fibre spéciale $X\otimes_{\FF_q}\kappa$ de $\calX$. Notons
$\calO'$ l'anneau local de $\calX$ en le point $\eta_\kappa$, $K'$
son corps des fractions, $\kappa'$ est son corps résiduel.

Donnons-nous $n$ morphismes zéros $u_i\in (X'-I)(\calO)$ avec
$i=1,\ldots,n$. Pour alléger les notations, nous allons noter par
les mêmes notations $u_i$ les $K$-points de $X'-I$ qui s'en
d\'eduisent par extension des scalaires de $\calO$ à $K$.

Donnons-nous un point
$$\tilde\calV^K=((u_j)_{j=1}^n,(\calV^K_j)_{j=0}^n,
(t_j)_{j=1}^n,\epsilon)$$ de $\DICht_\ul$ \`a valeurs dans $K$
au-dessus de $(u_j)_{j=1}^n\in (X'-I)^n(K)$. On doit montrer qu'il
peut s'étendre en un point de $\DICht_\ul$ \`a valeur dans $\calO$
au-dessus de $(u_j)_{j=1}^n\in (X'-I)(\calO)$. On va utiliser le
lemme suivant qui est la conjonction de \cite[Proposition 7]{Se}
et \cite[9.4.3]{EGA1}.

\begin{lemme}
Soient $\calX$ un schéma régulier de dimension $2$ et
$j:U\hrta\calX$ une immersion ouverte dont le fermé complémentaire
est de dimension zéro. Alors les foncteurs $j_*$ et $j^*$ entre la
catégorie des $\calO_\calX$-Modules localement libres et la
catégorie des $\calO_U$-Modules  modules localement libres sont
quasi-inverses l'une de l'autre et définissent une équivalence.
\end{lemme}

Prolonger un fibré vectoriel $\calV_j^K$ sur $X\otimes_{\FF_q} K$
consiste donc à prolonger sa fibre générique $V_i$, de $\Spec(K')$
à $\Spec(\calO')$, autrement dit, à choisir un $\calO'$-réseau
dans le $K'$-espace vectoriel $V_j$. Dans le $\calD$-chtouca
$\tilde V^K$, les fibres génériques $V_j$ de $\calV_j^K$ peuvent
être identifiées à l'aide des modifications $t_j$ et nous allons
les noter tous $V$. L'isomorphisme
$\epsilon:\,^\sigma\calV_0^K\isom {\calV'}_0^K$ induit une
application $\sigma$-linéaire injective $\epsilon_V:V\rta V$
c'est-à-dire une application $\epsilon_V$ vérifiant
$\epsilon_V(\alpha v) =\sigma(\alpha) v$ pour tout $\alpha\in
K'=F\otimes_{\FF_q}K$ o\`u $\sigma$ est d\'efini par
$\sigma(a\otimes b)=a\otimes b^q$ pour tout $a$ dans le corps des
fonctions $F$ de $X$ et $b\in K$.

Rappelons le lemme suivant de Drinfeld \cite{D}.

\begin{lemme} Soit $V$ un $K'$-espace vectoriel muni
d'une application $\sigma$-linéaire injective $\epsilon_V:V\rta
V$. Il existe alors des $\calO'$-réseaux stables sous $\epsilon_V$
et parmi ceux-ci il en existe un qui contient tous les autres.
\end{lemme}

Soit $\pi$ un uniformisant de l'anneau local $\calO$ qu'on
l'utilise aussi comme un uniformisant de $\calO'$. Soient $V$ un
$K'$-espace vectoriel muni d'une application $\sigma$-linéaire
injective $\epsilon_V:V\rta V$. Soit $M$ un $\calO'$-réseau
$\epsilon_V$-stable de $V$. Alors $M$ est dit {\em admissible si
l'application induite $$\bab \epsilon_V:M/\pi M\rta M/\pi M$$
n'est pas nilpotente}.

Rappelons un autre lemme de Drinfeld \cite{D}.

\begin{lemme} Gardons les notations du lemme
précédent et notons $M_0$ le $\calO'$-réseau $\epsilon_V$-stable
maximal de $V$. Alors, s'il existe un seul $\calO'$-réseau
$\epsilon_V$-stable admissible dans $V$, alors le réseau $M_0$ est
aussi admissible. De plus, quitte à remplacer $K$ par une
extension finie $L$, $(V,\epsilon_V)$ par $(V\otimes_K
L,\epsilon_V\otimes \sigma_L)$, $V$ contiendra un réseau
$t_V$-stable admissible.
\end{lemme}

Revenons à notre situation où $V$ est la fibre générique des
composantes $\calV_i^K$ du $\calD$-chtouca $\tilde\calV^K$. Quitte
à remplacer $K$ par une extension finie, on peut supposer que le
$\calO'$-réseau $M_0$ $t_V$-stable maximal est admissible. En
utilisant $M_0$, on étend alors canoniquement chaque $\calV^j_K$
en un fibré $\calV_j$ sur $\calX$. Par fonctorialité, chaque
modification $t_j$ s'étend en un isomorphisme entre les
restrictions de $\calV_{j-1}$ et $\calV_j$ en dehors du graphe de
$u_j$. De m\^eme, l'action de $\calO_X$-Algèbre $\calD$ sur
$\calV_j^K$ et ${\calV'}_I^K$ s'étendent en une action $\calD$ sur
$\calV_j$ et $\calV'_j$. En fin, l'isomorphisme
$\epsilon:\,^\sigma \calV_0^K \rta {\calV}_n^K$ s'étend en un
morphisme $\epsilon:\,^\sigma \calV_0 \rta {\calV}_n$. Il reste à
vérifier les trois assertions suivantes :

\begin{proposition}
\begin{enumerate}

\item les $t_j:\calV_{j-1}^{\{u_j\}}\rta \calV_j^{\{u_j\}}$
sont des modifications d'invariant $\inv_{u_j}(t_j)\leq \lambda^j$.

\item $\epsilon:\,^\sigma \calV_0\rta {\calV}_n$ est un
isomorphisme.

\item les $\calD\boxtimes \calO_\calX$-Modules $\calV_j$ sont
localement libres.

\end{enumerate}
\end{proposition}

\dem La première assertion est évidente. En effet $t_i$ est une
modification d'invariant au plus $\lambda_i$ en fibre générique,
il en est de même en fibre spéciale.

Vérifions la seconde assertion. D'après le lemme 3, il suffit de
démontrer que $\epsilon$ induit un isomorphisme entre les fibres
de $^\sigma\calV_0$ et de $\calV_n$ au point générique
$\eta_{\kappa}$ de la fibre spéciale de $\calX$. Bien entendu, il
revient au même de démontrer que l'application $\bab \epsilon:\bab
M_0\otimes_{\kappa',\sigma}\kappa' \rta \bab M_0$ est un
isomorphisme.

Par hypothèse d'admissibilité, l'application linéaire $$\bab
\epsilon:\bab M_0\otimes_{\kappa',\sigma}\kappa' \rta \bab M_0$$
n'est pas nilpotente. Considérons la filtration de $\bab M_0$
définie par les images des puissances de $\bab \epsilon$. Celle-ci
est nécessairement stationnaire puisqu'il s'agit d'une filtration
d'un espace vectoriel de dimension finie. Il existe donc un
sous-espace vectoriel non nul $N$ de $\bab M_0$ qui est égale à
l'image de $\bab \epsilon^n$ pour tous $n$ assez grand. Supposons
que $\bab \epsilon$ n'est pas un isomorphisme, $N$ est alors un
sous-$\kappa'$-espace vectoriel non nul et strict de $\bab M_0$.

Par construction, l'application linéaire restreinte $\bab
\epsilon: N \otimes_{\kappa',\sigma}\kappa' \rta N$ est surjective
donc bijective. Par hypothèse $\bab \epsilon$ commute à l'action
de $D$ de sorte que le sous-espace vectoriel $N$ est $D$-stable.
Pour tout $j=0,\ldots,n$, le sous-espace vectoriel $N$ définit par
saturation un sous-fibré $\bab\calV_{j}^N$ de
$\bab\calV_j=\calV_j\otimes_\calO\kappa$. Puisque $N$ est stable
sous l'action de $D$, les saturés $\bab\calV^N_j$ sont des
$\calD_X\boxtimes \calO_S$-modules.

Les modifications $t_j:\calV_{j-1}^{\{u_j\}}\rta \bab
\calV_j^{\{u_j\}}$ d'invariants au plus $\lambda_j$ induisent des
morphismes injectifs de $\calO_X$-Modules
$$\calV_{j}\rta \calV_{j-1}(||\lambda_j||\bab u_j).$$
Par $S$-platitude, ceci induit un morphisme injectif sur la fibre
spéciale $\calV_{n}\rta \calV_0(\sum_{j=1}^n||\lambda_j||\bab
u_j)$ et donc un morphisme injectif $$\bab\calV_{n}\rta \bab
\calV_0(\sum_{j=1}^n||\lambda_j||\bab u_j).
$$
Par restriction à $N$, on a un morphisme injectif
$$\bab \calV_{n}^N\rta \bab \calV_0^N(\sum_{j=1}^n||\lambda_j||\bab u_j)$$
On en déduit l'inégalité
$$\deg(\bab\calV_n^N)\leq \deg(\bab\calV_0^N)+d^2
\sum_{j=1}^N||\lambda||.$$

Par ailleurs, on a un morphisme $^\sigma\calV_0\rta \calV_N$ qui
par restriction à la fibre spéciale, puis à $N$, induit un
homomorphisme
$$\,^\sigma \bab\calV_0^N\rta \bab\calV_n^N.$$
Cet homomorphisme est aussi injectif puisqu'il l'est génériquement
par la construction de $N$. L'inégalité de degrés obtenue plus
haut implique que son conoyau est de longueur au plus
$d^2\sum_{j=1}^N||\lambda_j||$. Ce conoyau supporté par un
sous-schéma fini de $X$ appelé le dégénérateur de Drinfeld.

Par hypothèse sur le nombre de places totalement ramifiées de $D$,
il existe une place $v$ de $X-X'$ qui évite le dégénérateur où $D$
est totalement ramifiée. Il évite aussi les points
$u_1,\ldots,u_n\in X'$. On a alors un isomorphisme
$$\bab\epsilon:\,^\sigma\bab\calV_{0,v}^N\rta \bab\calV_{0,v}^N$$
de $\calD_v\otimes_{\FF_q}\kappa$-modules. En prenant les vecteurs
fixes par cet isomorphisme on un $D_v$-sous-module ${\rm
Ker}[\bab\epsilon-1,\bab\calV_{0,v}^N]$ de dimension sur $F_v$
égal à la dimension de $N$ sur $F\otimes_{\FF_q}\kappa$, en
particulier strictement plus petite que $d^2$. Ceci n'est pas
possible puisque $D_v$ est une algèbre à division centrale de
dimension $d^2$ sur son centre $F_v$.

On a démontré un isomorphisme $^\sigma \calV_0\isom \calV_n$ sur
$X\times S$. Pour prolonger la structure de $\calD\boxtimes
\calO_S$ sur les $\calV_i$, on peut maintenant invoquer
\cite[proposition 7, I.3]{LL} et conclure la démonstration de la
proposition. \findem

\section{Systèmes locaux fondamentaux}

On va continuer à supposer que l'algèbre à division $D$ a un
nombre de places totalement ramifiées plus grand que
$d^2||\lambda||$. Sous cette hypothèse, le morphisme
caractéristique $c_\lambda^I$ est propre d'après 1.6.

Pour tout niveau $I$, on dispose d'un morphisme propre
$$c^I_\ul:\calS^I_\ul\rta (X'-I)^n$$ et d'un faisceau pervers
$\calF^I_\ul$ sur $\calS^I_\ul$ qui est localement acyclique par
rapport \`a $c^I_\ul$. Les images directes
$\rmR^i(c^I_\ul)_*\calF^I_\ul$ sont donc des {\em syst\`emes
locaux}. On va consid\'erer {\em le système local gradué}
$\calW^I_\ul$ en prenant pour sa $i$-ème composante le système
local $\rmR^i(c^I_\ul)_*\calF^I_\ul$. Cette somme directe est
munie d'une action de l'algèbre de Hecke $\calH^I$.

\subsection{Opérateurs de Hecke et points fixes}
\setcounter{lemme}{0}

On va rappeler la construction des correspondances de Hecke et
leurs actions sur $\calW^I_\ul$. Soit $K^I$ le sous-groupe compact
de $G(\calO_\bbA)=\prod_{v\in|X|} G(\calO_x)$ défini comme le
noyau de l'homomorphisme naturel $[G(\calO_\bbA)\rta G(\calO_I)]$
qui se décompose en produit des sous-groupes compacts $K^I_{v}$ de
$G(\calO_v)$, égaux à $G(\calO_v)$ pour $v\notin|I|$. Pour toutes
places $v\in |X|$, notons $\calH_{v}^I$ l'algèbre des fonctions
$K^I_{v}$-bi-invariantes et à support compacts sur $G(F_v)$ et
$\calH^I$ leur produit restreint.

Rappelons $\calH^I_v$ admet une base $(\phi_{\beta_v})$ formée des
fonctions carac\-té\-ristiques des $K^I_v$-double-classes,
indexées par l'ensemble des double-classes
$$\beta_v\in K^I_v\backslash G(F_v)/K^I_v.$$
Lorsque $v\in |X'-I|$, cet ensemble est en bijection canonique
avec $\ZZ^d_+$.

Considérons la catégorie des triplets
$$(\calV^I_v,{\calV'}^I_v,t)$$
où $\calV^I_v$ et ${\calV'}^I_v$ sont des $\calG^I_v$-torseurs sur
$\calO_v$ et où $t$ est un isomorphisme entre leurs restrictions à
$F_v$. Puisque $\calG^I_v(F_v)=G(F_v)$ et
$\calG^I_v(\calO_v)=K^I_v$, l'ensemble des classes d'isomorphisme
de cete catégorie est l'ensemble des double-classes $\beta_v\in
K^I_v\backslash G(F_v)/K^I_v$.

Soit $v\in |X|$ et $\beta_v\in K^I_v\backslash G(F_v)/K^I_v$ une
double-classe. On va considérer  $\Phi_{\beta_v}$ le champ qui
associe à tout $\FF_q$-schéma $S$ la catégorie en groupoïdes des
triplets $(\tilde\calV,\tilde\calV';t'')$ formées des données
suivantes
\begin{itemize}
\item de deux $\calD$-chtoucas $$\tilde\calV=((u_j)_{j=1}^n,
(\calV_j)_{j=0}^n,(t_j)_{j=0}^n,\epsilon,\iota)$$ et
$$\tilde\calV'=((u_j)_{j=1}^n,(\calV'_j)_{j=0}^n,
(t'_j)_{j=0}^n,\epsilon',\iota)$$ de modifications d'invariant
$\inv_{u_j}(t_j)\leq \lambda_j$ et $\inv_{u_j}(t'_j)\leq
\lambda_j$, tous les deux munis de $I$-structure $\iota$ et
$\iota'$ respectivement, et ayant les mêmes zéros et poles
$$u_1,\ldots,u_n\in (X'-I-\{v\})(S),$$ autrement dit c'est un
$S$-point de
$$\DICht_\ul\times_{(X'-I)^n}
\DICht_\ul\times_{(X'-I)^n} (X'-I-\{v\})^n,$$

\item d'une $\{v\}$-modification $t'':\tilde\calV^{\{v\}}\isom
\tilde \calV'^{\{v\}}$ compatible avec $t$ et $t'$ et qui vérifie
la condition suivante. Le $\calD$-chtouca $\tilde \calV$ muni
d'une $I$-structure de niveau définit en $v$ un
$\calG^I_v$-torseur $\calV^I_v$ sur $\calO_v$ comme dans 1.5. De
même, $\tilde \calV'$ définit un $\calG^I_v$-torseur
${\calV'}^I_v$ sur $\calO_v$. La $\{v\}$-modification $t''$
définit un isomorphisme entre les restrictions de ces deux
$\calG^I_v$-torseurs à $F_v$. On demande que la classe
d'isomorphisme du triplet $(\calV^I_v,{\calV'_v},t'')$ soit la
double-classe $\beta_v\in K^I_v\backslash G(F_v)/K^I_v$.
\end{itemize}

\begin{proposition}
Les projections
$$\pr_1,\pr_2: \Phi_{\beta_v} \rta \DICht_\ul\times_{(X'-I)^n}
(X'-I-\{v\})^n
$$
sont des morphismes finis et étales. De plus, on a $f_\ul\circ
\pr_1=f_\ul \circ \pr_2$ où $f_\ul:\DICht_\ul\rta \calQ_\ul$ est
le morphisme  défini en 1.4.
\end{proposition}

\dem Voir \cite[I.4 proposition 3]{LL} pour la démonstration du
fait que $\pr_1$ et $\pr_2$ sont finis et étales. L'égalité
$f_\ul\circ \pr_1=f_\ul \circ \pr_2$ est claire sur définitions.
\findem

\bigskip
Par passage au quotient par l'action libre de $a^\ZZ$,
on en déduit des morphismes finis étales
$$\pr_1,\pr_2 :\Phi_{\beta_v}/a^\ZZ \rta
\calS^I_\ul\times_{(X'-I)^n}(X'-I-\{v\})^n
$$
De même, on a $f_\ul\circ \pr_1=f_\ul \circ \pr_2$ où
$f_\ul:\calS_\ul\rta \calQ_\ul$ est le morphisme défini en 1.4. On
en déduit un isomorphisme $\pr_1^*\calF_\ul \isom \pr_2^*
\calF_\ul$ et donc une correspondance cohomologique de
$(\calS_\ul,\calF_\ul)$ au-dessus de $(X'-I-\{v\})^n$. La
correspondance $\Phi_{\beta_v}$ détermine ainsi un endomorphisme
de la restriction de $\calW^I_\ul$ à $(X'-I-\{v\})^n$. Puisque
$\calW^I_\ul$ est un système local gradué sur $(X'-I)^n$, on
obtient par prolongement un endomorphisme de $\calW^I_\ul$ qu'on
notera $\Phi_{\beta_v}$ également.

\begin{proposition}
Lorsque $v\in|X|$ et $\beta_v\in K^I_v\backslash G(F_v)/K^I_v$
varient, les $\Phi_{\beta_v}$ définissent une action de $\calH^I$
sur $\calW^I_\ul$.
\end{proposition}

\dem Voir \cite[1.4 théorème 5]{LL} \findem

\subsection{La forme conjecturale du module virtuel}
\setcounter{lemme}{0}

Pour tout système local $\calW$ sur $(X'-I)^n$ muni d'une action
de $\calH^I$, on note $[\calW]$ sa classe dans le groupe de
Grothendieck des $\calH^I\otimes \pi_1(X'-I)^n$-modules. Notons
$[\calW^I_\ul]$ la somme alternée
$$[\calW^I_\ul]=\sum_{i} (-1)^i [\calW^I_{\ul,i}].$$
L'expression suivante est une variante des conjectures de
Langlands sur la cohomologie des variétés de Shimura. Notons
l'absence des termes endoscopiques dans notre cas où le groupe $G$
est une forme intérieure de $\GL_d$.

Soient $\lambda_1,\ldots,\lambda_n\in \ZZ^d_+$ tels que
$\sum_{j=1}^n \sum_{i=1}^d \lambda^i_j=0$.  Alors on devrait avoir
$$[\calW^I_\ul]=\bigoplus_\pi m(\pi)
[\pi^{K_I}\otimes \pr_1^* \calL_{\lambda_1}(\pi)
\otimes\cdots\otimes \pr_n^*\calL_{\lambda_n}(\pi)]$$
où
\begin{itemize}
\item
$\pi$ parcourt l'ensemble des représentations automorphes de
$D^\times$ sur lesquelles $a^\ZZ$ agit trivialement,

\item le sous-groupe compact $K_I$ de $G(\bbA_F)$ est
associé au niveau $I$ de façon usuelle; $\pi^{K_I}$ est
l'espace des vecteurs fixes de $\pi$ sous ce groupe.

\item $\pr_i$ est la projection de $(X'-I)^n$ sur son $i$-ème
facteur.

\item $\calL_{\lambda_i}(\pi)$ est le système local sur $(X'-I)$
obtenu en composant le paramètre de Langlands de $\pi$
$$\sigma_\pi:\pi_1(X'-I)\rta \GL_d(\QQ_\ell)$$
avec la représentation algébrique irréductible de $\GL_d(\QQ_\ell)$
de plus haut poids $\lambda_i$.
\end{itemize}

En particulier, lorsque $n=1$, et $\lambda=\lambda_1$ tel que
$|\lambda|=\sum_{i=1}^d \lambda^i =0$ on devrait avoir la formule
simple
$$[\calW^I_\lambda]=\bigoplus_\pi m(\pi)
[\pi^{K_I}\otimes \calL_{\lambda}(\pi)].$$

\section{L'action du groupe symétrique $\frakS_r$}

Dans ce chapitre nous continuons à supposer que $D$ a un grand
nombre de places totalement ramifiées comme dans l'énoncé de la
proposition 2 de 1.6.

\subsection{Situation A : descente à la Weil}
\setcounter{lemme}{0}

Soit $\lambda\in\ZZ^d_+$ tel que $|\lambda|=0$. Pour tout entier
naturel $r$, considérons les produit cartésien $r$ fois
$$(c^I_\lambda)^r : (\calS^I_\lambda)^r \rta (X'-I)^{r}.$$
sur lesquels le groupe symétrique $\frakS_r$ agit de manière
évidente.

Notons $A$ la somme directe des images directes dérivées $\rmR^i
(c^I_\lambda)^r_* \calF_\lambda^{\boxtimes r}$ qui est
$(\calW^I_\lambda)^{\boxtimes r}$ par Kunneth. Le système local
gradué $A$ est muni en plus de l'action de $(\calH^I)^{\otimes r}$
d'une action de $\frakS_r$ qui relève de l'action de $\frakS_r$
sur $(X'-I)^{r}$.

Soient $x\in (X'-I)(k)$ et $x^r\in (X'-I)^{r}(k)$ le point
diagonal dont toutes les coordonnées sont $x$. L'actions de
$\frakS_r$ sur $(X'-I)^{r}$ définit alors une action de $\frakS_r$
sur le groupe fondamental $\pi_1((X'-I)^{r},x^r)$. Dire que
l'action $\frakS_r$ se relève sur $A$ revient à dire que la
représentation du groupe $\pi_1((X'-I)^{r},x^r)$ sur $A_x$ se
prolonge en une représentation du produit semi-direct
$$\pi_1((X'-I)^{r},x^r)\rtimes \frakS_r.$$
Il y a une équivalence entre la catégorie des systèmes locaux sur
$(X'-I)^{r}$ munis d'une action compatible de $\frakS_r$ et la
catégorie des représentations du produit semi-direct
$\pi_1((X'-I)^{r},x^r)\rtimes \frakS_r$.

La fibre $A_{x^r}$ de $A$ en $x^r$ est alors un espace vectoriel
gradué muni des actions compatibles de
$\pi_1((X'-I)^{r},x^r)\rtimes \frakS_r$ et de $(\calH^I)^{\otimes
r}$. En admettant la formule conjecturale du paragraphe 2.3, on
devrait avoir une égalité dans le groupe de Grothendieck
$$[A_{x^r}]=\Bigl[\bigoplus_{\pi_1,\ldots,\pi_r}
\bigotimes_{i=1}^r m(\pi_i) \pi_i^{K_I}\otimes \bigotimes_{i=1}^r
\calL_{\lambda}(\pi_i)_{x}\Bigr].
$$
où
\begin{itemize}
\item  les $\pi_1,\ldots,\pi_r$ parcourent l'ensemble des représentations
automorphes de $G$ sur lesquelles $a^\ZZ$ agit trivialement,

\item $m(\pi_i)$ est la multiplicité de $\pi_i$,

\item $\calL_{\lambda}(\pi_i)$ est le système local sur $(X'-I)$
associé à la représentation automorphe $\pi_i$ et au copoids
dominant $\lambda$,

\item l'action d'un élément $\tau\in\frakS_r$ envoie le terme
indexé par $(\pi_1,\ldots,\pi_r)$ sur le celui indexé par
$(\pi_{\tau^{-1}(1)},\ldots,\pi_{\tau^{-1}(r)})$.
\end{itemize}

À la place de l'action de tout le groupe symétrique $\frakS_r$, on
aura besoin de regarder que l'action du sous-groupe des
permutations cycliques $\ZZ/r\ZZ$ de $\frakS_r$. {\em On notera
désormais $\tau$ le générateur de ce groupe cyclique}. En tordant
la $\FF_q$-structure évidente sur $(X'-I)^r$ par l'automorphisme
$\tau$, on obtient la $\FF_q$-structure obtenue par la descente à
la Weil de $\FF_{q^r}$ à $\FF_q$.

\subsection{Situation B : $\calD$-chtoucas à modification symétrique}
\setcounter{lemme}{0}

Soit $\lambda\in \ZZ^d_+$ tel que $|\lambda|=0$. Pour tout entier
naturel $r$, considérons la suite $r.\lambda$
$$(r.\lambda):=
(\underbrace{\lambda,\ldots,\lambda}_r).
$$
Supposons que $\calD$ a au moins $rd^2||\lambda||$ places
totalement ramifiées de sorte que le morphisme
$$c^I_{r.\lambda}:\calS^I_{r\lambda}\rta (X'-I)^{r}$$
est propre d'après 1.6. La propreté de $c^I_{r.\lambda}$ et
l'acyclicité locale du complexe d'intersection $\calF_{r.\lambda}$
impliquent que les images directes supérieures $\rmR^i
(c^I_{r.\lambda})_*\calF_{r.\lambda}$ sont des systèmes locaux.
Considérons le système local gradué $\calW^I_{r.\lambda}$ défini
en prenant la somme alternée et notons-le $B$. En plus de l'action
de l'algèbre de Hecke $\calH^I$, $B$ est muni d'une action du
groupe symétrique $\frakS_r$ défini comme suit.

Soit $U^I_{r}$ l'ouvert dans $(X'-I)^{r}$ complémentaire de la
réunion de toutes les diagonales.  L'action restreinte de
$\frakS_r$ à $U^I_{r}$ est libre.

\begin{proposition}
L'action de $\frakS_r$ sur $U^I_{r}$ se relève sur
$\DCht^I_{r.\lambda}\times_{(X'-I)^{r}}U^I_{r}$.
\end{proposition}

\dem D'après la propriété de factorisation, proposition 3 de 1.2,
sur l'ouvert $U^I_{r}$ la donnée d'un point de
$\DCht^I_{r.\lambda}$ au-dessus de $u=(u_i)_{i=1,\ldots,r}\in
U^I_{r}$, est la donnée d'une seule modification $^\sigma
\calV_0^T \rta \calV_0^T$ où $T$ est la réunion des points $u_i$,
ayant un invariant inférieur où égal à $\lambda$ en chaque point
$u_i$. Cette description ne dépend de l'ordre entre les $u^i$. On
en déduit une action du groupe symétrique $\frakS_r$. \hfill
$\square$

\bigskip
L'action de $\frakS_r$ sur $\DCht^I_{r\lambda}\times_{(X'-I)^{r}}
U^I_{r}$ définie ci-dessus, commute de manière évidente à l'action
libre de $a^\ZZ$. On en déduit une action de $\frakS_r$ sur
$\calS^I_{r.\ul}$. Par conséquent, $\frakS_r$ agit sur la
restriction du système local gradué $B$ à $U^I_r$. Mais puisque
c'est un système local, cette action se prolonge en une action de
$\frakS_r$ sur $B$ relevant l'action de $\frakS_r$ sur
$(X'-I)^{r}$.

Soit $x^r\in (X'-I)^r(k)$ le point diagonal choisi en 3.1. La
fibre $B_{x^r}$ de $B$ en $x^r$ est alors munie des actions
commutantes de $\pi_1((X'-I)^{r},x^r)\rtimes \frakS_r$ et de
$\calH^I$. En admettant la formule de la section 2.3, on devrait
avoir une égalité dans le groupe de Grothendieck
$$[B_x]=\bigoplus_{\pi} m(\pi) [\pi^{K_I} \otimes \bigotimes_{i=1}^r
 \calL_{\lambda}(\pi)_{x}]$$
où $\pi$ parcourt l'ensemble des représentations automorphes de
$G$ sur lesquelles $a^\ZZ$ agit trivialement et où l'action de
$\frakS_r$ sur la somme directe devrait se déduire de l'action de
$\frakS_r$ par permutation sur le terme $\bigotimes_{i=1}^r
\calL_{\lambda}(\pi)$.

\subsection{L'argument heuristique de comparaison}
\setcounter{lemme}{0}

Rappelons qu'on s'est donné un élément $\lambda\in\ZZ^d_+$ tel que
$|\lambda|=0$ et d'une algèbre à division $D$ sur $F$ ayant un
nombre de places totalement ramifiées plus grand que
$d^2r||\lambda||$. Cette dernière hypothèse garantit que les
morphismes caractéristiques $c^I_\lambda$ et $c^I_{r.\lambda}$
sont propres.

Soient $A_{x^r}$ et $B_{x^r}$ les espaces vectoriels gradués,
munis d'actions du groupe fondamental $\pi_1((X'-I)^{r},x^r)$, de
l'algèbre de Hecke $\calH^I$ et du groupe symé\-trique, définis
dans 3.1 et 3.2.

\begin{theoreme}
Soit $x\in (X'-I)(k)$ et $x^r$ le point diagonal de $(X'-I)^r$ de
coordonnées $x$. Pour tous $g\in\pi_1((X'-I)^{r},x^r)$ et
$f\in\calH^I$on a l'égalité
$$ {\rm Tr}((f\otimes1\otimes\cdots\otimes 1)
g\tau, A_{x^r})={\rm Tr}(fg\tau,B_{x^r})$$ où $\tau\in\frakS_r$
est la permutation cyclique.
\end{theoreme}

Il nous semble utile de présenter un argument heuristique en se
basant sur les expressions plausibles, voir 2.3, de $[A_x]$ et de
$[B_x]$.

Dans la somme $[A_x]$, $\tau$ envoie le terme indexée par
$(\pi_1,\ldots,\pi_r)$ sur le terme indexé par
$(\pi_r,\pi_1,\ldots,\pi_{r-1})$. Considérons l'ensemble des
termes qui forme une orbite cyclique sous le groupe engendré par
$\tau$. La somme des termes dans cette orbite est stable sous
$(f\otimes 1\otimes\cdots \otimes 1)g\tau$, mais la trace de cet
opérateur sur cette somme est nulle, sauf si l'orbite est
constituée d'un seul élément, c'est-à-dire $\pi_1=\cdots=\pi_r$.
Pour un tel terme, diagonal, la comparaison résulte d'un lemme
général d'algèbre, qui a été sans doute connu de Saito et Shintani
\cite{Sh}.

\begin{lemme}[Saito-Shintani]
Soit $V$ un espace vectoriel de dimension finie sur un corps $K$.
Notons $\tau$ l'endomorphisme de $V^{\otimes r}$ défini par
$$v_1\otimes\cdots\otimes v_r\mapsto v_r\otimes v_1\otimes\cdots\otimes
v_{r-1}.$$
Pour tous endomorphismes $f_1,\ldots,f_r$, on a
$${\rm Tr}(f_1f_2\ldots f_r,V)={\rm Tr}((f_1\otimes\cdots\otimes f_r)\tau,
V^{\otimes r}).$$
\end{lemme}

Nous allons en fait démontrer directement théorème 1 et en déduire
le lemme fondamental pour le changement de base pour des fonctions
sphériques de $\GL_n$ de degré $0$. La démonstration est basée sur
le comptage et la théorie des modèles locaux.

\section{Comptage}

\subsection{Problème de comptage}
\setcounter{lemme}{0}

On va fixer une clôture algébrique $k$ de $\FF_q$. Donnons-nous un
entier naturel $s\geq 1$, et deux fermés finis disjoints $T$ et
$T'$ de $X$. Supposons de plus que $T$ évite le lieu $X-X'$ où
l'algèbre à division $D$ se ramifie.

Considérons la catégorie $\frakC(T,T';s)$ dont les objets sont des
$\calD$-fibrés $\calV$ sur $X\otimes_{\FF_q}k$ muni de deux
modifications :
\begin{itemize}
\item une $T$-modification $t:\,^\sigma \calV^{T}\rta \calV^{T}$,

\item une $T'$-modification $t':\,^{\sigma^s} \calV^{T'} \rta
\calV^{T'}$,
\end{itemize}

\noindent telles que le diagramme
$$\renewcommand\arraystretch{1.5}
\begin{array}{ccc}
^{\sigma^{s+1}} \calV & \hfld{{\sigma^s}(t)}{} & ^{\sigma^s} \calV\\
\vfld{\sigma(t')}{} & & \vfld{}{t'} \\
^{\sigma} \calV & \hfld{}{t} & \calV
\end{array}$$
commute sur le plus grand ouvert de $\bab X$ où toutes les flèches
sont définies.

Les morphismes de la catégorie $\frakC(T,T';s)$ sont des
isomorphismes entre les triplets $(\calV;t,t')$ augmentés
formellement d'un isomorphisme entre $(\calV;t,t')$ et
$(\calV\otimes\calL(a);t\otimes\id_{\calL(a)},
t'\otimes\id_{\calL(a)})$ pour tous $(\calV;t,t')$. Ici,
$\calL(a)$ est le fibré en droites sur $X$ associé à l'idèle $a$
qu'on a fixé d'emblée. Un isomorphisme entre deux objets de cette
catégorie, sera appelé un $a$-isomorphisme.

Soit $I$ un fermé fini de $X$, étranger de $T$. Soit
$K^I=\prod_{v\in|X|}K^I_v$ sous-groupe ouvert compact de
$G(\calO_\bbA)$ défini dans 2.1. Donnons-nous deux fonctions
$$\alpha_T:|T\otimes_{\FF_q}\FF_{q^s}|\rta \ZZ^d_+$$ et
$\beta_{T'}$ qui associe à toute place $v\in|T'|$ une
double-classe
$$\beta_v\in K^I_v\backslash G(F_v)/K^I_v.$$
Considérons la catégorie
$$\frakC^I_{\alpha_T,\beta_{T'}}(T,T';s)$$
dont les objets $(\calV,t,t',\iota)$ consistent en les données
suivantes
\begin{itemize}
    \item un $T$-chtouca  $(\calV,t)$ d'invariant de Hodge
    $\alpha_T$, à valeur dans $k$ et munis de $I$-structure
    de niveau $\iota$. Ici $\alpha_T$ est vu comme la fonction
    $T\otimes_{\FF_q}k$ qui se factorise par
    $\alpha_T:|T\otimes_{\FF_q}\FF_{q^s}|\rta \ZZ^d_+$.
    \item une $T'$-modification $t':\,^{\sigma^s}\calV^{T'}
    \isom \calV^{T'}$, commutant à $t$ et qui en toute place
    $v\in T'$, est de type $\beta_v$ dans le sens de 2.1.
\end{itemize}

Notre problème de base consiste à exprimer le nombre
$$\#\frakC^I_{\alpha_T,\beta_{T'}}(T,T';s)
:=\sum_{(\calV;t,t';\iota) } {1\over \#{\rm
Isom}(\calV;t,t';\iota)}$$ où $(\calV;t,t';\iota)$ parcourt un
ensemble de représentants des classes d'isomorphisme d'objets de
$\frakC^I_{\alpha_{T},\beta_{T'}} (T,T';s)$, en termes
d'inté\-grales orbitales et d'inté\-grales orbitales tordues en
suivant le comptage à la Kottwitz des points des variétés de
Shimura.

\subsection{La fibre générique de $(\calV;t,t')$}
\setcounter{lemme}{0}

Soient $T,T'$ deux sous-schémas fermés finis disjoints de $\bab
X$, $T\subset \bab X'$ et $s\geq 1$ un entier naturel comme
précédemment. Soit $(\calV;t,t')$ un objet de $\frakC(T,T';s)$.
Notons $V$ la fibre générique de $\calV$; c'est un
$D^{op}\otimes_{\FF_q} k$-module libre de rang $1$. Les
modifications $t$ et $t'$ induisent des applications bijectives
$\tau$ et $\tau'$ de $V$ dans $V$ qui sont respectivement
$\id_{D^{op}}\otimes\sigma$-linéaire et
$\id_{D^{op}}\otimes\sigma^s$-linéaire. L'hypothèse de commutation
entre $t$ et $t'$ implique que $\tau\tau'=\tau'\tau$.

\begin{proposition}
Soient $(\calV;t,t')$ un objet de $\frakC(T,T';s)$ et
$(V;\tau,\tau')$ sa fibre générique. Soit $v$ un point fermé de
$X-T$, $V_v$ le complété de $V$ en $v$, $\tau_v:V_v\rta V_v$ la
bijection $\id_{D_v}\hat\otimes\sigma$-linéaire induite de $\tau$.
Alors, il existe un isomorphisme
$$(V_v,\tau_v)\simeq (D_v\hat\otimes_{\FF_q}k,\id_{D_v}\hat\otimes\sigma).$$
De même, soit $v$ un point fermé de $X-T'$. Alors il existe un
isomorphisme $$(V_v,\tau'_v)\simeq
((D_v\hat\otimes_{\FF_q}\FF_{q^s})\hat\otimes_{\FF_{q^s}}k,
\id_{D_v\hat\otimes_{\FF_q}\FF_{q^s}}\hat\otimes\sigma^s).$$
\end{proposition}

\dem Pour une place $v\notin T$, le
$D_v\hat\otimes_{\FF_q}k$-module $V_v$ admet un
$\calD_v\hat\otimes_{\FF_q}k$-réseau $\calV_v$ tel que
$\tau_v(\calV_v)=\calV_v$. Il s'ensuit que $V_v^{\tau_x=1}$ est un
$\calD_v$-module libre de rang $1$ dont le produit tensoriel
complété avec $V_v^{\tau_x=1}\hat\otimes_{\FF_q}k$ est $V_v$. Il
s'ensuit que la paire $(V_v,\tau_v)$ est isomorphe à
$(D_v\hat\otimes_{\FF_q}k,\id_{D_v}\hat\otimes\sigma)$. Le même
argument vaut pour $T'$. \hfill $\square$

\bigskip
Soit $C(T,T';s)$ la catégorie dont les objets sont des triplets
$(V;\tau,\tau')$ où $V$ est un $D\otimes_{\FF_q}k$-module libre de
rang $1$, $\tau$ et $\tau'$ sont des bijections $V\rta V$ qui sont
$\id_D\otimes\sigma$-linéaire et $\id_D\otimes\sigma^s$-linéaire
respectivement qui vérifient les propriétés suivantes

\begin{itemize}

\item $\tau$ et $\tau'$ commutent : $\tau\tau'=\tau'\tau$,

\item pour toute place  $v\notin |T|$, on a $(V_v,\tau_v)$ est
isomorphe à
$$(D_v\hat\otimes_{\FF_q}k,\id_{D_x}\hat\otimes\sigma)$$

\item pour toute place $v\notin |T'|$, on a $(V_v,\tau'_v)$ est
isomorphe à
$$((D_v\hat\otimes_{\FF_q}\FF_{q^s})\hat\otimes_{\FF_{q^s}}k,
\id_{D_v\hat\otimes_{\FF_q}\FF_{q^s}}\hat\otimes\sigma^s).$$
\end{itemize}
\noindent Les flèches de $C(T,T';s)$ sont des isomorphismes entre
les triplets $(V;\tau,\tau')$.

\begin{proposition}
Supposons que $D$ a au moins une place totalement ramifiée. Alors
tout objet $(V;\tau,\tau')$ de $C(T,T';s)$ est irréductible.
\end{proposition}

\dem Soit $x\in|X|$ une place où $D$ est totalement ramifiée.
Puisque $x\notin T$, la paire $(V_x,\tau_x)$ est isomorphe à
$(D_x\hat\otimes k,\id_{D_x}\hat\otimes\sigma)$. Cette paire est
irréductible puisque $D_x$ est une algèbre à division centrale sur
$F_x$. \findem

\bigskip
D'après la proposition 1, on a un foncteur de la catégorie
$\frakC(\bab T,\bab T';s)$ dans la catégorie $C(T,T';s)$ qui
associe à un triplet $(\calV;t,t')$ sa fibre générique
$(V;\tau,\tau')$. {Sans savoir a priori que ce foncteur est
essentiellement surjectif}, on va appeler $C(T,T';s)$ la catégorie
des fibres génériques de $\frakC(\bab T,\bab T';s)$, et chercher à
classifier les classes d'isomorphisme des objets de $C(T,T';s)$.
Pour cela il nous faut rappeler la théorie des $\phi$-espaces de
Drinfeld.

\subsection{Rappels sur les $\phi$-espaces}
\setcounter{lemme}{0}

La théorie des $\phi$-espaces est due à Drinfeld. Une référence
utile est l'appendice de \cite{LRS}.

\begin{definition}
Un {\em $\phi$-espace} est un $F\otimes_{\FF_q} k$-espace
vectoriel $V$ muni d'une application $\id_F\otimes\sigma$-linéaire
bijective $\phi:V\rta V$. Un morphisme de $\phi$-espaces
$\alpha:(V_1,\phi_1)\rta (V_2,\phi_2)$ est une application
$F\otimes_{\FF_q} k$-linéaire $\alpha:V_1\rta V_2$ telle que
$\phi_2\circ\alpha =\alpha\circ \phi_1$. On note $\calF$ la
catégorie des $\phi$-espaces.
\end{definition}

Drinfeld a attaché à tous $\phi$-espaces un invariant qu'il
appelle {\em $\phi$-paire} dont on va rappeler la définition.

Soit $E$ une extension de $F$. Notons $X_{E}$
le normalisé de la courbe $X$ dans $E$.
Notons ${\rm Div}^0(X_{E})$ le groupe des diviseurs de degré
$0$ sur la courbe $X_{E}$.
On a l'homomorphisme de groupes abéliens
$$E^\times \rta {\rm Div}^0(X_{E})$$
qui associe à une fonction rationnelle non nulle $f\in E$ sur
$X_E$, son diviseur ${\rm div}(f)$. Le conoyau de cette
application, le groupe des diviseurs de $X_{E}$ de degré $0$
modulo les diviseurs principaux, est le groupe des $\FF_q$-points
de la jacobienne de $X_E$. En particulier, c'est un groupe abélien
fini. Son noyau, étant le groupe des fonctions constantes non
nulles en $X_{E}$, est aussi un groupe abélien fini. Il s'ensuit
que l'homomorphisme de $\QQ$-espaces vectoriels
$$E^\times\otimes\QQ\rta {\rm Div}^0(X_E)\otimes \QQ$$
qui s'en déduit, est un isomorphisme.

\begin{definition} Une $\phi$-paire est une paire $(E,Q)$ où $E$
est une $F$-algèbre finie et où $Q$ est un élément de
$E^\times\otimes\QQ$.

Soit $\calE'$ la catégorie dont les objets sont des $\phi$-paires
$(E,Q)$, une flèche $(E,Q)\rta (E',Q')$ dans $\calE'$  est un
homomorphisme injectif de $F$-algèbres $\phi:E\rta E'$ tel que
$\phi(Q)=Q'$. On appelle la catégorie des $\phi$-paires, la
catégorie $\calE$ obtenue de $\calE'$ en inversant toutes les
flèches de $\calE'$.
\end{definition}

Soit $(E,Q)$ un objet de $\calE$. Parmi les $F$-sous-algèbres $E'$
de $E$ telles que $Q$ soit dans l'image de l'inclusion évidente
${E'}^\times\otimes\QQ\hookrightarrow E^\times\otimes\QQ$, il en
existe une qui est contenue dans toutes les autres; notons-la
$\tilde E$. Si $E=\tilde E$, la $\phi$-paire $(E,\phi)$ sera dite
{\em minimale}. Si $(\tilde E,Q)$ est une $\phi$-paire minimale,
on vérifie facilement que $\tilde E$ est une $F$-algèbre étale.
Dans $\calE$, toute $\phi$-paire $(E,Q)$ est isomorphe à une
$\phi$-paire $(\tilde E,Q)$ minimale, unique à isomorphismes non
uniques près.

Soit $(E,Q)$ une $\phi$-paire minimale où $E$ est un corps. Choisissons
un plongement de $E$ dans une clôture séparable $F^{sep}$ de $F$. L'image
de $Q$ dans $(F^{sep})^\times\otimes\QQ$ est bien déterminée à
${\rm Gal}(F^{sep}/F)$-conjugaison près. Il est clair qu'on peut
retrouver la classe d'isomorphisme de $(E,Q)$ à partir de cet élément
de $[(F^{sep})^\times\otimes\QQ]/{\rm Gal}(F^{sep}/F)$.

\begin{proposition}
L'ensemble des classes d'isomorphisme
des $\phi$-paires $(E,Q)$ où $E$ est un corps est en bijection
naturelle avec $[(F^{sep})^\times\otimes\QQ]/{\rm Gal}(F^{sep}/F)$.
\end{proposition}

On peut construire un foncteur de la catégorie $\calF$ des
$\phi$-espaces dans la catégorie $\calE$ des $\phi$-paires comme
suit. Soit $(V,\phi)$ un $\phi$-espace. Considérons l'ensemble des
paires $(n,V_n)$ formée d'un entier naturel $n$ et d'un
$F\otimes_{\FF_q} \FF_{q^n}$-sous-espace vectoriel de $V$, stable
par $\phi$, tel que $V=V_n\otimes_{\FF_{q^n}}k$. Cet ensemble
possède un ordre partiel : $(n,V_n) \leq (n',V_{n'})$ si et
seulement si $n$ divise $n'$ et $V_{n'}=V_n\otimes_{\FF_{q^n}}
\FF_{q^{n'}}$. Cet ordre est filtrant.

Soit $(n,V_n)$ une $\FF_{q^n}$-structure de $(V,\phi)$ comme ci-dessus.
La restriction de $\phi^n$ à $V_n$ définit alors un automorphisme linéaire
de $V_n$. Notons $E$ la $F$-sous-algèbre de ${\rm End}(V_n)$ engendrée
par la restriction de $\phi^n$ à $V_n$. Notons $Q$ l'élément de
${\rm Div}^0(X_E)\otimes \QQ$ telle que
$$nQ={\rm div}(\phi^n|_{V_n}).$$
On obtient une $\phi$-paire $(E,Q)$.

Si $(n,V_n)\leq (n',V'_n)$ et $(E',Q')$ est la $\phi$-paire
associée à $(n',V'_n)$, on a visiblement un $F$-homomorphisme
canonique $E'\rta E$ qui envoie $Q'$ sur $Q$. Ceci définit un
isomorphisme canonique dans la catégorie $\calE$. Puisque l'ordre
$(n,V_n)\leq (n',V'_n)$ est filtrant, les différentes $\phi$-paires
$(E,Q)$ qu'on a construit se diffèrent par un unique isomorphisme
de $\calE$. On obtient ainsi un foncteur $\calF\rta \calE$.

Rappelons un théorème de Drinfeld.

\begin{theoreme}

\begin{enumerate}
\item La catégorie des $\phi$-espaces sur $k$ est abélienne
$F$-linéaire et semi-simple.

\item Le foncteur $\calF\rta\calE$ ci-dessus induit une bijection
entre l'ensemble des classes d'isomorphisme de $\phi$-espaces
irréductibles et l'ensemble des classes d'isomorphisme dans
$\calE$ des $\phi$-paires $(E,Q)$ où $E$ est un corps.

\item Pour une extension finie $E$ de $F$ et un élément $Q\in {\rm
Div}^0(X_E)\otimes\QQ$. Supposons que la $\phi$-paire $(E,Q)$ est
minimale. Notons $b$ le plus petit entier naturel tel que $$Q\in
{1\over b}{\rm Div}^0(X_E).$$ Soit $(V,\phi)$ le $\phi$-espace
irréductible qui correspond à la $\phi$-paire $(E,Q)$, on a
$$\dim_{F\otimes_{\FF_q}k}(V)=[E:F]b.$$ L'algèbre ${\rm
End}(V,\phi)$ est une algèbre à division centrale sur $E$ de
dimension $b^2$ et d'invariant $${\rm inv}({\rm End}(V,\phi))\in
{\rm Div}^0(X_E)\otimes\QQ/\ZZ$$ égal à l'image de $-Q\in {\rm
Div}^0(X_E)\otimes\QQ$ dans ${\rm Div}^0(X_E)\otimes\QQ/\ZZ$.
\end{enumerate}
\end{theoreme}

La classe d'isomorphisme d'un $\phi$-espace irréductible
$(V,\phi)$ est donc déter\-minée par une $\phi$-paire minimale
$(E,Q)$, et donc en définitive, est déterminée par l'image de $Q$
dans $[(F^{sep})^\times\otimes\QQ]/\Gal(F^{sep}/F)$.

\subsection{La classe de conjugaison $\gamma_0$}
\setcounter{lemme}{0}

Soit $(V;\tau,\tau')$ un objet de la catégorie $C(T,T';s)$. On
peut produire une bijection $D\otimes_{\FF_q}k$-linéaire de $V$ à
partir des applications semi-linéaires $\tau$ et $\tau'$ en
prenant le composé
$$\gamma=\tau^s {\tau'}^{-1}: V\rta V.$$ En choisissant une
rigidification $V\isom D^{op}\otimes_{\FF_q} k$, $\gamma$ définit
un élément de
$$G(F\otimes_{\FF_q}k)=(D\otimes_{\FF_q} k)^\times.$$ Nous notons
$\gamma_0$ la classe de conjugaison de $\gamma$ dans
$G(F\otimes_{\FF_q}k)$ qui, bien entendu, ne dépend pas de la
rigidification choisie.

\begin{lemme}
La classe de conjugaison $\gamma_0$ de $G(F\otimes_{\FF_q}k)$ est
$\sigma$-invariante.
\end{lemme}

\dem On déduit de l'hypothèse de commutation $\tau\tau'=\tau'\tau$
un carré commutatif
$$\renewcommand\arraystretch{1.5}
\begin{array}{ccc}
V & \hfld{\tau}{} & V\\
\vfld{\gamma}{} & & \vfld{}{\gamma} \\
V & \hfld{}{\tau} & V
\end{array}$$
Que $\gamma$ commute avec un isomorphisme $\sigma$-linéaire
implique que $\gamma$ et $\sigma(\gamma)$ sont conjugués. \findem

\begin{proposition}
La classe de conjugaison $\gamma_0$ dans $G(F\otimes_{\FF_q}k)$
rencontre $G(F)$. De plus, cette intersection est formée d'une
seule classe de $G(F)$-conjugaison.
\end{proposition}

\dem $G$ est une forme intérieure de $\GL(d)$ de sorte qu'il
existe une injection canonique Galois équivariante de l'ensemble
des classes de conjugaison de $G(F\otimes_{\FF_q}k)$ dans
l'ensemble des classes de conjugaison $\GL_d(F\otimes_{\FF_q}k)$.
On peut donc transférer $\gamma_0$ en une classe de conjugaison de
$\GL_d(F\otimes_{\FF_q}k)$ qui est $\sigma$-invariante. Celle-ci
provient donc d'une unique classe de conjugaison de $\GL_d(F)$.

Le problème consiste donc à vérifier que cette classe  de
conjugaison $\tilde \gamma$ de $\GL_d(F)$ provient de $G(F)$. Il
revient au même de trouver un plongement de l'algèbre
$F[\tilde\gamma]$ engendrée par $\tilde\gamma$ dans $D^{op}$.
D'après \cite[lemme 4, III.3]{LL}, ceci est un problème local en
des places $v$ de $X$ où $D$ se ramifie, c'est-à-dire $\gamma_0$
vue comme classe de conjugaison de $\GL_d(F_v)$ provient de
$G(F_v)$. Or, $T$ évite ces places par hypothèses si bien que la
classe de conjugaison $\gamma_0$ vue dans $G(F_v \otimes_{\FF_q}
k)$ vient d'une classe de conjugaison dans $G(F_v)$. \findem

\bigskip
Choisissons un représentant de la classe de conjugaison $\gamma_0$
dans $G(F)=D^\times$ qu'on va encore noter $\gamma_0$. Notons $E=F[\gamma_0]$
la sous-$F$-algèbre de $D$ engendré par $\gamma_0$. Puisque $D$ est une
algèbre à division centrale sur $F$, $E$ est un corps. Notons $X_E$
le normalisé de $X$ dans $E$. Notons $\pi:X_E\rta X$.

\begin{definition}
La classe de conjugaison $\gamma_0$ de $D^\times$ est dite {\em
$(T,T')$-admissible} si le diviseur ${\rm div}(\gamma_0)$ de $X_E$
est supporté par $\pi^{-1}(T\cup T')$ et si en décomposant ${\rm
div}(\gamma_0)$ en somme
$${\rm div}(\gamma_0)={\rm
div}_T(\gamma_0)+{\rm div}_{T'}(\gamma_0)$$
d'un diviseur supporté
par $\pi^{-1}(T)$ et d'un diviseur supporté par $\pi^{-1}(T')$,
les diviseurs ${\rm div}_T(\gamma_0)$ et ${\rm
div}_{T'}(\gamma_0)$ sont tous les deux de degré $0$.
\end{definition}

\begin{proposition}
Soit $(V;\tau,\tau')$ un objet de $C(T,T';s)$. Soit
$\gamma=\tau^s{\tau'}^{-1}$. Soit $\gamma_0$ la classe de
conjugaison de $G(F)$ qui est
$G(F\otimes_{\FF_q}\bab\FF_q)$-conjuguée à $\gamma$. Alors
$\gamma_0$ est une classe de conjugaison de $G(F)$ qui est
$(T,T')$-admissible.

On note encore $\gamma_0\in D^\times$ un représentant de cette
classe de conjugaison $\gamma_0$, soit $E_0$ la $F$-sous-algèbre de
$D$ engendrée par $\gamma_0$. Soit
$${\rm div}(\gamma_0)={\rm div}_T(\gamma_0)+{\rm div}_{T'}(\gamma_0)$$
la décomposition du diviseur sur $X_{E_0}$
associé la fonction rationnelle $\gamma_0$, en la somme d'un diviseur
supporté par $T$ avec un diviseur supporté par $T'$. Alors, dans
la catégorie $\calE$, la $\phi$-paire $$\left(E,{{\rm
div}_T(\gamma_0)\over s}\right)$$ est isomorphe à la $\phi$-paire
associée au $\phi$-espace $(V,\tau)$.
\end{proposition}

\dem Il existe un entier naturel $n$ assez divisible, en
particulier $s|n$, tel qu'il existe un
$D\otimes_{\FF_q}\FF_{q^n}$-module libre $V'$ contenu dans $V$ tel
que $V=V'\otimes_{\FF_{q^n}}k$ et tel que $\tau$ et $\tau'$ laisse
stable $V'$. Il est clair que $(\tau|_{V'})^n$ et
$(\tau'|_{V'})^{n/s}$ dont des automorphismes
$D\otimes_{\FF_q}\FF_{q^n}$-linéaires de $V'$. On a
$$(\gamma|_{V'})^{n/s}=(\tau|_{V'})^{n}(\tau'|_{V'})^{-n/s}.
$$
Soit $E$ la $F$-sous-algèbre de ${\rm End}(V')$ engendrée par
$(\tau|_{V'})^n$ et $(\tau'|_{V'})^n$; c'est un corps puisque que
$(V;\tau,\tau')$ est irréductible d'après la proposition 2 de 4.2.

Soit $X_E$ le normalisé de $X$ dans $E$. On a donc une égalité de
diviseurs sur $X_E$
$${\rm div}((\gamma|_{V'})^{n/s})={\rm div}((\tau|_{V'})^{n})
+{\rm div}((\tau'|_{V'})^{-n/s}).
$$
On va utiliser le lemme suivant.

\begin{lemme} Dans l'expression ci-dessus, ${\rm div}
((\tau|_{V'})^{n})$ est supporté par $T$ et ${\rm
div}((\tau'|_{V'})^{-n/s})$ est supporté par $T'$.
\end{lemme}

\noindent{\em Démonstration du lemme 5}. Soit $x\in |X-T|$. Par
hypothèse $(V_x,\tau_x)$ est isomorphe à
$(D_x\otimes_{\FF_q}k,\id_{D_x}\otimes\sigma)$. Il existe donc un
entier $n''$ divisible par $n$ tel que
$(V'_x,\tau_x)\otimes_{\FF_{q^n}} \FF_{q^{n''}}$ est isomorphe à
$(D_x\otimes_{\FF_q}\FF_{q^{n''}},\id_{D_x} \otimes\sigma)$.
L'automorphisme linéaire de $(V'_x,\tau_x)\otimes_{\FF_{q^n}}
\FF_{q^{n''}}$
$$(\tau_x|_{V'_x\otimes_{\FF_{q^n}}\FF_{q^{n''}}})^{n''}
=((\tau_x|_{V'_x\otimes_{\FF_{q^n}}\FF_{q^{n''}}})^{n'})^{n''/n'}
$$
est donc égal à l'identité. Il s'ensuit que $x$ n'appartient pas
au support de ${\rm div} ((\tau|_{V'})^{n})$. Le même argument
vaut pour $\tau'$ et pour $x\in |X-T'|$. \findem

\bigskip
\noindent{\em Suite de la démonstration de la proposition 4} . Il
est clair, grâce au lemme ci-dessus, que ${\rm
div}((\gamma|_{V'})^n)$ est un diviseur supporté par $T\cup T'$,
dont la partie supportée par $T$ est ${\rm
div}((\tau|_{V'})^{ns})$ et dont la partie supportée par $T'$ est
${\rm div}((\tau'|_{V'})^{-n})$. Puisque ce sont, tous les deux,
des diviseurs principaux, leurs degrés sont nuls. On en déduit que
la class de conjugaison $\gamma_0$ est $(T,T')$-admissible car
dans la définition de la propriété $(T,T')$-admissible, il est
loisible de remplacer l'extension engendré par $\gamma_0$ par
n'importe quelle extension qui la contient, et de remplacer
$\gamma_0$ par une puissance de $\gamma_0$.

Il reste à démontrer que la $\phi$-paire $(E_0,{\rm
div}_{T}(\gamma_0)/s)$ est isomorphe à la $\phi$-paire associée au
$\phi$-espace $(V,\tau)$, dans la catégorie $\calE$. Il existe un
homomorphisme d'algèbres $E_0\rta E$ qui envoie $\gamma_0$ sur
$\gamma$ de sorte qu'on se ramène à démontrer l'assertion pour la
$\phi$-paire $(E,{\rm div}_T(\gamma)/s)$.

Soit $E'$ la $F$-sous-algèbre de ${\rm End}(V')$ engendrée par
$(\tau|_{V'})^{n}$. Par construction la $\phi$-paire associée au
$\phi$-espace $(V;\tau)$ est $(E',{\rm div}((\tau|_{V'})^{n})/n)$.
Il est clair que $E'$ est une sous-extension de $E$ qui est la
$F$-sous-algèbre ${\rm End}(V')$ de engendrée par
$(\tau|_{V'})^{n}$ et $(\tau'|_{V'})^{n}$ de sorte que dans la
catégorie $\calE$, les $\phi$-paires $(E',{\rm
div}((\tau|_{V'})^{n})/n)$ et $(E,{\rm div}((\tau|_{V'})^{n})/n)$
sont isomorphes. Sur $X_E$, l'égalité $${{\rm
div}((\tau|_{V'})^{n})\over n}={{\rm div}_T(\gamma)\over s}$$ se
déduit du lemme 5. \findem

\bigskip
Soit $x\in |T|$. On a une application $s$-norme de l'ensemble des
classes de $\sigma$-conjugaison dans
$G(F_x\otimes_{\FF_q}\FF_{q^s})$ dans l'ensemble des classes de
conjugaison stable de $G(F_x)$.  En une place $x\in |T|$ donc
$x\in |X'|$ , $G$ est isomorphe à ${\rm GL}_d$ de sorte que les
classes de conjugaison stable sont simplement des classes de
conjugaison. On a donc une application $s$-norme
$${\rm N}_{x,s}:G(F_x\otimes_{\FF_q}\FF_{q^s})/\mbox{\small $\sigma$-conj.}
\rta G(F_x)/\mbox{\small conj.}$$

\begin{definition}
Soit $s\in\mathbb N$. La classe de conjugaison $\gamma_0$ de
$G(F)$ est dite une $s$-norme en $T$ si pour tout $x\in |T|$,
$\gamma_0$ vue comme classe de conjugaison de $G(F_x)$ est dans
l'image de ${\rm N}_{x,s}$.

\end{definition}

\begin{proposition}
Soit $(V;\tau,\tau')$ un objet de $C(T,T';s)$. Soit
$\gamma=\tau^s{\tau'}^{-1}$. Soit $\gamma_0$ la classe de
conjugaison de $G(F)$ qui est $G(F\otimes_{\FF_q}k)$-conjuguée à
$\gamma$. Alors $\gamma_0$ est une $s$-norme en $T$.
\end{proposition}

\dem Soit $x\in |T|$. Par hypothèse $x\notin |T'|$. Le complété
$V_x$ de $V$ en $x$ est un $D_x\otimes_{\FF_q}k$-module libre de
rang $1$ muni d'une bijection $\sigma$-linéaire $\tau_x$ et d'une
bijection $\sigma^s$-linéaire $\tau'_x$. Puisque $x\notin T'$, la
paire $(V_x,\tau'_x)$ est isomorphe à
$((D_x\hat\otimes_{\FF_q}\FF_{q^s}) \hat\otimes_{\FF_{q^s}} k,
\id_{D_x\hat\otimes_{\FF_q}\FF_{q^s}}\hat\otimes \sigma^s)$. En
prenant les vecteurs fixes sous $\tau'_x$, on obtient donc un
$D_x\hat\otimes_{\FF_q}\FF_{q^s}$-module libre $V'_x$. Puisque
$\tau_x$ commute à $\tau'_x$, $\tau_x$ induit sur $V'_x$ une
application $\id_{D_x}\hat\otimes\sigma$-linéaire qui définit une
classe de $\sigma$-conjugaison dans
$G(F_x\hat\otimes_{\FF_q}\FF_{q^s})$ qu'on notera $\delta_x$. Par
construction $\tau'_x|_{V'_x}=1$ de sorte que
$\gamma_x|_{V'_x}=\tau_x^s|_{V'_x}$. On en déduit que $\gamma_0$
vue comme classe de conjugaison de $G(F_x)$ est la $s$-norme de la
classe de $\sigma$-conjugaison $\delta_x$. \findem

\subsection{Les classes d'isogénie}
\setcounter{lemme}{0}

On a vu qu'un objet $(V;\tau,\tau')$ de $C(T,T';s)$ détermine une
classe de conjugaison $\gamma_0$ dans $G(F)$ qui est
$(T,T')$-admissible et $s$-norme en $T$. On va voir que,
réciproquement, une classe de conjugaison de $G(F)$ qui est
$(T,T')$-admissible et $s$-norme en $T$, détermine uniquement une
classe d'isomorphisme de $C(T,T';s)$. Ceci constituera l'analogue
de la théorie de Honda-Tate dans notre problème de comptage. Il
s'agit encore d'une variante d'un théorème de Drinfeld, déjà
généralisé par Laumon, Rapoport, Stuhler et par Lafforgue à
d'autres contextes.

\begin{theoreme}
L'application qui associe à un objet $(V;\tau,\tau')$ de
$C(T,T';s)$ la classe de conjugaison $\gamma_0$ de $G(F)$
conjuguée dans $G(F\otimes_{\FF_q}k)$ à
$\gamma=\tau^s{\tau'}^{-1}$, définit une bijection de l'ensemble
des classes d'isomorphisme de $C(T,T';s)$ sur l'ensemble des
classes de conjugaison de $G(F)$ qui sont $(T,T')$-admissibles et
$s$-normes en $T$.
\end{theoreme}

\dem  Il s'agit de vérifier que l'application $(V;\tau,\tau')\mapsto\gamma_0$
est injective et surjective. Vérifions d'abord la {\em surjectivité}.

Soit $\gamma_0$ une classe de conjugaison de $G(F)$ qui est
$(T,T')$-admissible et $s$-norme en $T$. On va construire un objet
$(V;\tau,\tau')$ de $C(T,T';s)$ dont la classe de conjugaison dans
$G(F)$ associée est $\gamma_0$ en plusieurs étapes.

\bigskip
\noindent{\em Étape 1 : le $\phi$-espace irréductible $(W,\psi)$
facteur de $(V,\tau)$}. Prenons un repré\-sen\-tant de la classe
de conjugaison $\gamma_0$ qu'on note $\gamma_0\in D^\times$. Soit
$E=F[\gamma_0]$ la sous-$F$-algèbre de $D$ engendrée par
$\gamma_0$. C'est une extension de corps de $F$ puisque $D$ est
une algèbre à division. Le degré de cette extension $e=[E:F]$
divise l'entier $d$ le rang de $D$.

Notons $X_E$ le normalisé de $X$ dans $E$. D'après l'hypothèse de
$(T,T')$-admissibilité, le diviseur ${\rm div}(\gamma_0)$ se casse
en somme de deux diviseurs de degré $0$ supportés respectivement
par $T$ et $T'$ $${\rm div}(\gamma_0)={\rm div}_{T}(\gamma_0)+{\rm
div}_{T'}(\gamma_0).$$ Soit $E_T$ le plus petit sous-corps de $E$
tel que ${\rm Div}^0(X_{E_T}) \otimes\QQ$ contient encore la
$\QQ$-droite engendrée par ${\rm div}_T (\gamma_0)$. Notons
$e_T=[E_T:F]$.

La paire $(E_T,{\rm div}_{T}(\gamma_0)/s)$ est alors une
$\phi$-paire minimale. Soit $(W,\psi)$ un $\phi$-espace
irréductible correspondant à $(E_T,{\rm div}_{T}(\gamma_0)/s)$ ;
il est bien défini à isomorphisme près. Soit $b$ le plus petit
entier naturel tel que $${{\rm div}_{T}(\gamma_0)\over s}\in
{1\over b}{\rm Div}^0(X_{E_T}).$$ D'après le théorème 4 de 4.3,
${\rm End}(W,\psi)$ est une algèbre à division centrale sur $E_T$
de dimension $b^2$ sur $E_T$ et $${\rm dim}_{F\otimes
k}(W)=be_T.$$

\bigskip
\noindent{\em Étape 2 : structure de $D^{op}$-module sur le bon
multiple de $(W,\psi)$}. Considérons le multiple
$(V,\tau)=(W,\psi)^{r}$ avec $r=d^2/be_T$. L'espace vectoriel $V$
a donc la bonne dimension $d^2$ sur $F\otimes k$. Notons que
$$\dim_{E_T}M_{r}({\rm End}(W,\psi))=r^2 b^2={d^4\over e_T^2}.$$

Pour donner à $V$ une structure de $D^{op}$-module par rapport à
laquelle $\tau$ est $\id_D\otimes\sigma$-linéaire, il suffit de
construire un plongement $$D^{op}\hookrightarrow M_{r}({\rm
End}(W,\psi))$$ ou, ce qui est équivalent, un plongement
d'algèbres simples centrales sur $E_T$ $$D^{op}\otimes_F
E_T\hookrightarrow M_{r}({\rm End}(W,\psi)).$$ Ceci est encore
équivalent à l'existence d'une autre algèbre simple centrale $C$
sur $E_T$ telle que $$(D^{op}\otimes_F E_T)\otimes C \simeq
M_{r}({\rm End}(W,\psi)).$$ Une telle algèbre simple centrale $E$
doit avoir la dimension $$\dim_{E_T}(C)={\dim_{E_T}M_{r}({\rm
End}(W,\psi))\over \dim_{E_T}(D^{op}\otimes_F E_T)}={d^2\over
e_T^2}$$ et l'invariant $${\rm inv}(C)={\rm inv}({\rm
End}(W,\psi))-{\rm inv}(D^{op}\otimes_F E_T)$$ l'égalité étant
prise dans ${\rm Div}^0(X_{E_T})\otimes\QQ/\ZZ$. Pour que $C$
existe il faut et il suffit que $${\rm inv}({\rm
End}(W,\psi))-{\rm inv}(D^{op}\otimes_F E_T) \in {e_T\over d}{\rm
Div}^0(X_{E_T}) / {\rm Div}^0(X_{E_T}).$$

Rappelons un lemme \cite{LL} III.4 lemme 4.

\begin{lemme} Soit $E$ une extension de $F$ et $A$ une algèbre
simple centrale sur $E$ de dimension $d^2$. Soit $E'$ une
extension de $E$ de degré $e'$. Pour qu'il existe un plongement de
$E'\hookrightarrow A$, il faut et il suffit que $e'$ divise $d$ et
que l'image de ${\rm inv}(A)$ dans ${\rm Div}^0(X_{E'})\otimes
\QQ/\ZZ$ soit dans $${\rm inv}(A\otimes_E E')\in {e'\over d}{\rm
Div}^0(X_{E'})/{\rm Div}^0(X_{E'}).$$
\end{lemme}

On a les inclusions $E_T\subset E\subset D^{op}$ de sorte qu'en
vertu du lemme ci-dessus, on a $${\rm inv}(D^{op}\otimes_F E_T)\in
{e_T\over d}{\rm Div}^0(X_{E_T}) / {\rm Div}^0(X_{E_T}).$$ Par
ailleurs, grâce au théorème 4 de 4.3, on sait que l'invariant
$${\rm inv}({\rm End}(W,\psi))\in {\rm
Div}^0(X_{E_T})\otimes\QQ/\ZZ$$ est égal à l'image de $-{\rm
div}_T(\gamma_0)/s\in {\rm Div}^0(X_{E_T})\otimes\QQ$. Il reste
donc à démontrer que $${{\rm div}_T(\gamma_0)\over s} \in
{e_T\over d}{\rm Div}^0(X_{E_T}).$$ On va utiliser le lemme
suivant.

\begin{lemme}
Soit $\gamma_0$ un élément de $D^\times$ dont la classe de
conjugaison est une $s$-norme en $T$. Soit $E$ la $F$-sous-algèbre
de $D$ engendrée par $\gamma_0$. Notons $e=[E:F]$. Alors on a
$${\rm div}_T(\gamma_0)\in {e\over d}{\rm Div}^0(X_E)$$ où ${\rm
div}_T(\gamma_0)$ est la partie ${\rm div}(\gamma_0)$ supporté par
$T$.
\end{lemme}

\dem Soit $D'$ le centralisateur de $E$ dans $D$. C'est une
algèbre à division centrale de rang $d/e$ sur $E$. En remplaçant
$D$ par $D'$ et $F$ par $E$, on peut supposer que l'élément
$\gamma_0$ appartient à $F$. En une place $x\in |T|$, $D_x\simeq
{\rm GL}_d(F_x)$. Il existe donc $\delta\in {\rm
GL}_d(F_x\otimes_{\FF_q}\FF_{q^s})$ tel que $$\gamma_0=\delta
\sigma(\delta)\ldots \sigma^{s-1}(\delta).$$ On en déduit que
$$\det(\gamma_0)={\rm N}_s(\det(\delta))$$ de sorte que la
valuation de $\det(\gamma_0)$ est divisible par $s$. Or vu comme
élément du centre $F_x$, on a ${\rm det}(\gamma_0)=\gamma_0^d$
d'où le lemme. \findem

\bigskip
\noindent{\em Suite de l'étape 2}. L'inclusion $E_T\subset E$
induit un revêtement $\pi:X_E\rta X_{E_T}$ de degré $e/e_T$.
L'application injective $\pi^*:{\rm Div}^0(X_{E_T}) \rta {\rm
Div}^0(X_E)$ induit par tensorisation avec $\QQ$ l'application
usuelle $E_T^\times\otimes\QQ\rta E^\times\otimes \QQ$. Mais on a
aussi l'application $\pi_*:{\rm Div}^0(X_E)\rta {\rm
Div}^0(X_{E_T})$. Pour tout $Q\in {\rm Div}^0(X_{E_T})$, on a
$\pi_*(\pi^*(Q))={e\over e_T}Q$. On en déduit que si $$Q\in {\rm
Div}^0(X_E)\cap {\rm Im}[{\rm Div}^0(X_{E_T})\otimes\QQ]$$ alors
$\pi^*(\pi_*(Q))$ est aussi égale à ${e\over e_T}Q$.

D'après le lemme, il existe un diviseur $Q\in {\rm Div}^0(X_E)$
tel que ${{\rm div}_T(\gamma_0)/s}= {e\over d}Q$. En particulier,
$Q$ est dans ${\rm Div}^0(X_E)\cap {\rm Im}[{\rm Div}^0
(X_{E_T})\otimes\QQ]$ si bien que $Q={e_T\over e}\pi^*(\pi_*(Q))$.
On obtient donc une égalité de diviseurs de $X_{E_T}$ $${{\rm
div}_T(\gamma_0)\over s}={e_T\over d}\pi_*(Q)$$ avec $\pi_*(Q)\in
{\rm Div}^0(X_{E_T})$. C'est ce qu'il fallait pour qu'il existe
une algèbre simple centrale $C$ sur $E_T$ telle que $$
D^{op}\otimes_F C\simeq {\rm End}((W,\psi)^r).$$

\bigskip
\noindent{\em Étape 3 : plonger $E$ dans les endomorphismes de
$(V,\tau)$}. On a obtenu un $D^{op}\otimes_{\FF_q}k$-module $V$
muni d'une bijection $\id_D\otimes \sigma$-linéaire $\tau$. De
plus, l'anneau des endomorphismes de $(V,\tau)$ est égal à $C$. Il
reste à construire une bijection $\id_D\otimes \sigma^s$-linéaire
$\tau'$ qui commute à $\tau$. Il revient au même de construire une
action de $\gamma_0$ ou ce qui est équivalent de montrer
l'existence d'un plongement $E\hookrightarrow C$. D'après le lemme
2, il suffit de démontrer que $[E:E_T]=e/e_T$ divise $d/e_T$ et
que
$${\rm inv}(C)\in {e\over d}{\rm Div}^0(X_E)/{\rm Div}^0(X_E).$$
La première divisibilité est automatique. Par ailleurs, dans ${\rm
Div}^0(X_E)\otimes \QQ/\ZZ$, on a l'égalité
$${\rm inv}(C)=-{{\rm div}_T(\gamma_0)\over s}-{\rm inv}(D^{op}).$$
Puisque $E$ peut être plongé dans $D^{op}$, on a d'après le lemme
3
$${\rm inv}(D^{op})\in {e\over d}{\rm Div}^0(X_E)/{\rm Div}^0(X_E).$$
Par ailleurs, on a $${{\rm div}_T(\gamma_0)\over s}\in {e\over
d}{\rm Div}^0(X_E)$$ d'après le lemme 3. Donc, on peut plonger $E$
dans $C={\rm End}(V,\tau)$ et obtenir une bijection $\id_D\otimes
k$-linéaire de $V$ commutant à $\tau$. Posons
$\tau'=\gamma_0\tau^s$.

\bigskip
\noindent{\em Étape 4 : le triplet $(V,\tau,\tau')$ est bien un
objet de $C(T,T';s)$.} Soit $n$ un entier assez divisible tel
qu'il existe un $D\otimes_{\FF_q} \FF_{q^n}$-sous-module de $V'$,
stable par $\tau$ et $\tau'$, tel que $V=V'\otimes_{\FF_{q^n}}k$.
Notons $E'$ la $F$-sous-algèbre de ${\rm End}(V')$ engendrée par
$\tau^n|_{V'}$ et ${\tau'}^n|_{V'}$. On a dans $E'$ l'égalité
$$\gamma_0^n=\tau^{ns}|_{V'} {\tau'}^{-n}|_{V'}.$$
Les paires
$(E,\gamma_0^n)$ et $(E',\gamma_0^n)$ sont isomorphes dans la
catégorie $\calE$ de sorte que le diviseur ${\rm div}(\gamma_0^n)$
sur $X_{E'}$, est supporté par $T\cup T'$. Par construction de
$\tau$, la partie de ${\rm div}(\gamma_0^n)$ supportée par $T$ est
${\rm div}(\tau^{ns}|_{V'})$; il s'ensuit que la partie supportée
par $T'$ est ${\rm div}({\tau'}^{-n}|_{V'})$.

Soit $x\in |X-T|$. Par ce qui précède, l'élément
$\tau_x^{ns}|_{V'_x}$ est une unité de $E'_x$. L'application
$\id_D\otimes \sigma$-linéaire $\tau_x|_{V'_x}$ du
$D_x\otimes_{\FF_q}\FF_{q^n}$-module $V'_x$ a une $ns$-norme
$\tau_x^{ns}|_{V'_x}$ compacte. D'après Kottwitz, $\tau_x|_{V'_x}$
doit fixer un $\calD_x\hat\otimes_{\FF_q}\FF_{q^n}$-réseau de
sorte que la paire $(V_x,\tau_x)$ est isomorphe à la paire
$(D_x\hat\otimes_{\FF_q}k, \id_{D_x}\hat\otimes\sigma)$.

L'assertion sur $\tau'$ se démontre de même.

\bigskip\noindent{\em Étape 5 : la classe de conjugaison
associée à $(V;\tau,\tau')$ est bien $\gamma_0$}.
C'est évident sur la construction de $\tau'= \gamma_0 \tau^s$.

\bigskip
On a donc démontré que l'application $(V;\tau,\tau')\mapsto
\gamma_0$ de l'ensemble des classes d'isomorphisme des objets de
$C(T,T';s)$ sur l'ensemble des classes de conjugaison de $G(F)$
qui sont $(T,T')$-admissibles et $s$-norme en $T$, est {\em
surjective}. Nous allons maintenant vérifier qu'elle est {\em
injective}.

Soit $\gamma_0$ une classe de conjugaison de $G(F)$ qui est
$(T,T')$-admissible et $s$-norme en $T$. Soit $E$ la
$F$-sous-algèbre de $D$ engendrée par un représentant $\gamma_0$
de la classe de conjugaison $\gamma_0$.

Soit $(V;\tau,\tau')$ un objet de $C(T,T';s)$ dont la classe de
conjugaison associée dans $G(F)$, est $\gamma_0$. Puisque
$(V,\tau)$ est irréductible comme $D\otimes_{\FF_q}k$-module muni
d'une bijection $\id_D\otimes_{\FF_q}k$-linéaire, en oubliant
l'action de $D$, comme $\phi$-espace, $(V,\tau)$ est isotypique.
Le $\phi$-espace irréductible $(W,\psi)$, facteur de $(V,\tau)$,
est complètement déterminé par $\gamma_0$, car la $\phi$-paire qui
lui est associée par le théorème 4 de 4.3, est isomorphe dans
$\calE$ à la $\phi$-paire $(E,{\rm div}_T(\gamma_0))$. Ainsi
$(W,\psi)$ est nécessairement le $\phi$-espace irréductible
construit dans l'étape 1.

Comme $\phi$-espace, $(V,\tau)$ est donc nécessairement isomorphe à
$(W,\psi)^r$ où l'entier $r$ est défini comme dans l'étape 2. La
structure de $D^{op}$-module est alors donné par un homomorphisme
d'algèbres
$$D^{op}\hookrightarrow {\rm End}{(W,\psi)^r}.$$
Dans l'étape 2, on a démontré que celui-ci existe. Il est de plus unique à
automorphisme intérieur près, d'après le théorème de Skolem-Noether.

Il reste à construire une bijection $\id_D\otimes\sigma^s$-linéaire
$\tau'$ de $V$ commutant à $\tau$ telle que l'automorphisme
$\tau^s{\tau'}^{-1}$ est conjugué dans $G(F\otimes_{\FF_q}k)$ à
$\gamma_0$. Il revient au même de plonger l'algèbre $E=F[\gamma_0]$
dans le commutant $C$ de $D^{op}$ dans ${\rm End}{(W,\psi)^r}$.
On a démontré dans l'étape 3 que ce plongement existe. Il est
unique à automorphisme intérieur près, de nouveau d'après le théorème
de Skolem-Noether. \findem

\subsection{Le groupe des automorphismes de $(V;\tau,\tau')$}
\setcounter{lemme}{0}

Soit $\gamma_0$ une classe de conjugaison de $G(F)$ qui est
$(T,T')$-admissible et $s$-norme en $T$. Soit $(V;\tau,\tau')$
l'objet de $C(T,T';s)$ qui lui correspond et qui est bien défini à
isomorphisme près. On va démontrer que le groupe des
automorphismes de $(V,\tau,\tau')$ est le groupe des $F$-points
d'une forme intérieure $J_{\gamma_0}$ du centralisateur
$G_{\gamma_0}$ de $\gamma_0$ dans $G$. D'après le principe de
Hasse pour les groupes adjoints, la classe d'isomorphisme de cette
forme intérieure est complètement déterminée par celle les groupes
locaux $J_{\gamma_0,v}$, forme intérieures de $G_{\gamma_0,v}$
qu'on va décrire explicitement.

L'automorphisme $\gamma=\tau^s {\tau'}^{-1}$ de $V$ détermine une
classe de conjugaison de $G(F\otimes_{\FF_q} k)$ dont on choisit
un représentant noté aussi $\gamma$. Soit $G_\gamma(F\otimes_{\FF_q}k)$
le centralisateur de $\gamma$ dans $G(F\otimes_{\FF_q} k)$. Puisque
$\tau$ commute avec $\gamma$, $\tau$ détermine un automorphisme
de $G_\gamma(F\otimes_{\FF_q}k)$. Formons le produit semi-direct
$$G_\gamma(F\otimes_{\FF_q}k)\rtimes\langle\tau\rangle.$$
Il est clair que ${\rm Aut}(V,\tau,\tau')$ est le groupe des points
fixes de $G_\gamma(F\otimes_{\FF_q}k)$ sous l'action de $\tau$.
On veut construire un $F$-groupe $J_{\gamma_0}$ dont le groupe des
$F$-points est ${\rm Aut}(V,\tau,\tau')$.

On choisit un représentant $\gamma_0\in G(F)$ dans la classe de
conjugaison $\gamma_0$. Le centralisateur $G_{\gamma_0}$ étant
défini sur $F$, l'endomorphisme de Frobenius agit comme un
automorphisme sur $G_{\gamma_0}(F\otimes_{\FF_q}k)$. Formons le
produit semi-direct
$$G_{\gamma_0}(F\otimes_{\FF_q}k)\rtimes\langle \sigma\rangle.$$
Il est clair que le groupe $G_{\gamma_0}(F)$ est le groupe des points
fixes de $G_{\gamma_0}(F\otimes_{\FF_q}k)$ sous l'action de $\sigma$.

\begin{lemme}
Pour $r$ assez divisible, on a un isomorphisme de produits semi-directs
$$G_\gamma(F\otimes_{\FF_q} k)\rtimes\langle\tau^r\rangle
\isom G_{\gamma_0}(F\otimes_{\FF_q}k)\rtimes\langle \sigma^r\rangle.$$
\end{lemme}

\dem Prenons une $\FF_q$-structure de $D^{op}\otimes_{\FF_q}
k$-module $V$. Pour une extension assez grande de $\FF_q$, $\tau$
et $\tau'$ vont être définis sur cette extension si bien qu'il va
en être de même de $\gamma$. On peut donc supposer
$\gamma=\gamma_0$.

Soit $E=F[\gamma]$. Puisque l'homomorphisme $E^\times\rta {\rm
Div}^0(X_E)$ a un noyau et un conoyau finis, on peut décomposer
$$\gamma^r = \delta^{s}{\delta'}^{-1}$$
avec $\delta,\delta'\in E^\times$ tels que ${\rm div}(\delta)$
soit supporté par $T$, et ${\rm div}(\delta')$ par $T'$  en
utilisant l'hypothèse que $\gamma$ est $(T,T')$-admissible et
$s$-norme en $T$, pour un entier $r$ assez divisible. L'assertion
d'injectivité du théorème 1 de 4.5 montre que les triplets
$(V,\tau^r,{\tau'}^r)$ et $(V,\delta\sigma^r, \delta'
\sigma^{rs})$ sont isomorphes car ils correspondent tous les deux
à la classe de conjugaison de $\gamma^r$. Il ne reste qu'à
remarquer que $\delta\in E^\times$ appartient au centre de
$G_\gamma(F\otimes_{\FF_q}k)$ de sorte que les actions de
$\sigma^r$ et de $\delta\sigma^r$ sur
$G_\gamma(F\otimes_{\FF_q}k)$ sont les mêmes. La proposition s'en
déduit. \findem

\begin{proposition}
Soit $\gamma_0$ une classe de conjugaison de $G(F)$ qui est
$(T,T')$-admissible et $s$-norme en $T$. Soit $(V;\tau,\tau')$ la
classe d'isomorphisme de la catégorie des fibres génériques
$C(T,T';s)$ qui lui correspond par le théorème 1 de 4.5. Il existe
une forme intérieure $J_{\gamma_0}$ de $G_{\gamma_0}$ telle que
$$J_{\gamma_0}(F)={\rm Aut}(V;\tau,\tau').$$
\end{proposition}

\dem Le triplet $(V;\tau,\tau')$ induit une paire $(V;\gamma)$ où
on peut voir $V$ comme un $G$-torseur sur $F\otimes_{\FF_q}k$ et
où $\gamma:V\rta V$ est un automorphisme de ce $G$-torseur.
L'hypothèse que $\gamma$ et $\gamma_0$ soient conjugués dans
$G(F\otimes_{\FF_q}k)$, implique que la paire $(V,\gamma)$ est
isomorphe à $(G\otimes_{\FF_q} k, \gamma_0)$. L'ensemble des
isomorphismes
$$\calL:={\rm Isom}((G\otimes_{\FF_q} k,\gamma_0),(V,\gamma))$$
définit donc un $G_{\gamma_0}$-torseur sur $F\otimes_{\FF_q}k$.

Le diagramme commutatif
$$\renewcommand\arraystretch{1.5}
\begin{array}{ccc}
^\sigma V & \hfld{\tau}{} & V\\
\vfld{\sigma(\gamma)}{} & & \vfld{}{\gamma}\\
^\sigma V & \hfld{}{\tau} & V
\end{array}$$
induit un isomorphisme $\tau:\,^\sigma \calL\isom \calL$ et donc
une classe dans
$$\rmH^1(\langle\sigma\rangle, G_{\gamma_0}(F\otimes k)).$$
La proposition se ramène à démontrer que l'image de cette
classe dans
$$\rmH^1(\langle\sigma\rangle, G_{\gamma_0}^{\rm ad}(F\otimes k))$$
provient d'un cocycle continu, c'est-à-dire d'un élément de
$\rmH^1(\hat\ZZ, G_{\gamma_0}^{\rm ad}(F\otimes k))$. ce qui est
le contenu du lemme précédent. \findem

\bigskip
D'après le principe de Hasse pour les groupes adjoints,
la forme intérieure $J_{\gamma_0}$ de $G_{\gamma_0}$ est complètement
déterminée par les formes locales $J_{\gamma_0,v}$ de $G_{\gamma_0,v}$
qu'on peut décrire facilement à partir de $(V;\tau,\tau')$.

\begin{itemize}
\item {\em Soit $v$ une place de $X$ en dehors de $T$}. Par
hypothèse la paire $(V_v,\tau_v)$ est isomorphe à
$(D^{op}_v\hat\otimes_{\FF_q}k, \id_{D_v}\hat\otimes \sigma)$. Les
vecteurs fixes sous $\tau_v$ forment donc un $D^{op}_v$-module
libre $V_x^{\tau_v=1}$. Puisque $\tau'_v$ commute à $\tau$, il
laisse stable $V_v^{\tau_v=1}$ et sa restriction à
$V_v^{\tau_v=1}$ est une bijection $D^{op}_x$-linéaire dont la
classe de conjugaison est $\gamma_0$. Un automorphisme
$D^{op}_v\hat\otimes_{\FF_q} k$-linéaire de $V_v$ qui commutent à
$\tau_v$ laissent stable $V_v^{\tau_v=1}$. Si de plus, il commute
à $\tau'$, sa restriction à $V_v^{\tau_v=1}$ doit commuter à la
restriction de $\tau'_v$. Cela définit une bijection de
$J_{\gamma_0}(F_v)$ avec le centralisateur de $\gamma_0$ dans
$G(F_x)$.

\item Pour une place $x\in |T|$, on a nécessairement
$x\notin|T'|$,de sorte que la paire $(V_x,\tau'_s)$ est isomorphe
à
$$((D_x\hat\otimes_{\FF_q}\FF_{q^s})\hat\otimes_{\FF_{q^s}}k,
\id_{D_x\hat\otimes_{\FF_q}\FF_{q^s}}\hat\otimes \sigma^s).
$$
Les vecteurs fixes sous $\tau'_x$ forme un
$D_x\hat\otimes_{\FF_q}\FF_{q^s}$-module $V_x^{\tau'_x=1}$ libre
de rang $1$. Puisque $\tau_x$ commute à $\tau'_x$, $\tau_x$
définit une bijection $\sigma$-linéaire de $V_x^{\tau'_x=1}$ dans
lui-même et donc une classe de $\sigma$-conjugaison $\delta_x$ de
du groupe $G(F_x\hat\otimes_{\FF_q}\FF_{q^s})$. On déduit de la
relation $\gamma_0=\tau^s{\tau'}^{-1}$ que la $s$-norme de
$\delta_x$ est égale à $\gamma_0$. Le groupe $J_{\gamma_0}(F_x)$
est le alors le centralisateur tordu de $\delta_x$. \findem
\end{itemize}

\begin{corollaire}
Si en toutes les places $v$ de $F$, $J_{\gamma_0,v}$ est isomorphe
à $G_{\gamma_0,v}$ alors $J_{\gamma_0}$ est isomorphe à $G_{\gamma_0}$.
\end{corollaire}

\subsection{Compter les réseaux}
\setcounter{lemme}{0}

On est maintenant en mesure de démontrer une formule de comptage
pour la catégorie $\frakC^I_{\alpha_T,\beta_{T'}}(T,T';s)$ où
$\alpha_T$ est une fonction
$\alpha_T:|T\otimes_{\FF_q}\FF_{q^s}|\rta \ZZ^d_+$ et où
$\beta_{T'}$ est une fonction qui associe à toute place $v\in
|T'|$ une double classe $\beta_v\in K^I_v\backslash G(F_v)/K^I_v$.

On a construit un foncteur $$\frakC(T,T';s)\rta C(T,T';s)$$ qui
associe à un triplet $(\calV;t,t';\iota)$ sa fibre générique
$(V;\tau,\tau')$. On a ensuite démontré que l'application qui
associe à une classe d'isomorphisme d'objets $(V;\tau,\tau')$ de
$C(T,T';s)$, l'unique classe de conjugaison $\gamma_0$ de $G(F)$,
conjuguée dans $G(F\otimes_{\FF_q}k)$ à
$\gamma=\tau^s{\tau'}^{-1}$, induit une bijection de l'ensemble
des classes d'isomorphisme de $C(T,T';s)$ sur l'ensemble des
classes de conjugaison de $G(F)$ qui sont $(T,T')$-admissibles et
$s$-normes en $T$. On sait aussi qu'il existe une forme intérieure
$J_{\gamma_0}$ du centralisateur $G_{\gamma_0}$ telle que
$J_{\gamma_0}(F)$ est le groupe des automorphismes de
$(V;\tau,\tau')$ et pour toute place $x$ de $F$,
$J_{\gamma_0}(F_x)$ est le groupe des automorphismes de
$(V_x;\tau_x,\tau'_x)$.

\begin{theoreme}
 Le nombre pondéré des classes d'isomorphisme
$$
\#\frakC^I_{\alpha_T,\beta_{T'}}(T,T';s)
=\sum_{(\calV;t,t';\iota)}{1\over \#{\rm Isom}(\calV;t,t')}
$$
où $(\calV;t,t';\iota)$ parcourt un ensemble de représentants de
classes d'isomorphisme de $\frakC^I_{\alpha_T,\beta_{T'}} (\bab
T,\bab T';s)$, est fini et est égal à
$$
\displaystyle\sum_{(\gamma_0;\delta_T)} {\rm vol}(a^\ZZ
J_{\gamma_0,\delta_T}(F)\backslash J_{\gamma_0,\delta_T}(\mathbb
A_F)) \prod_{v\notin |T|} {\mathbf O}_{\gamma_0}(\phi_{\beta_v})
\prod_{x\in |T|} \mathbf{TO}_{\delta_x}(\phi_{\alpha_x}).
$$
Ici :
\begin{itemize}
\item la somme est étendue sur l'ensemble des paires
$(\gamma_0,\delta_T)$ formées d'une classe de conjugaison
$\gamma_0$ de $D^\times$ qui est $(T,T')$-admissible et d'une
collection $\delta_T=(\delta_x)_{x\in T}$ de classe de
$\sigma$-conjugaison $\delta_x$ de
$G(F_{x}\otimes_{\FF_q}\FF_{q^s})$ dont la norme est $\gamma_0$ ;

\item le groupe $J_{\gamma_0,\delta_T}$, défini dans 4.6, est la
forme intérieure de $G_{\gamma_0}$ dont les $F$-points forment le
groupe des automorphismes du triplet $(V;\tau,\tau')$
correspondant à la classe de conjugaison $\gamma_0$ et dont les
$F_x$-points forment le groupe des automorphismes ${\rm Aut}
(V_x;\tau_x,\tau'_x)$ ;

\item pour toute place $x\in|T|$, on définit
$$\phi_{\alpha_x}=\bigotimes_{y\in |x\otimes_{\FF_q}\FF_{q^s}|}
\phi_{\alpha_T(y)}
\in\bigotimes_{y\in |x\otimes_{\FF_q}\FF_{q^s}|}
\calH_{y}=\calH(G(F_x \otimes_{\FF_q}\FF_{q^s})) ;
$$
\item pour toute place $v\in |T'|$, on définit
$\phi_{\beta_v}\in\calH^I_v$ comme la fonction caractéristique de
la double classe $\beta_v\in K^I_v\backslash G(F_v)/K^I_v$. En une
place $v\notin |T|\cup |T'|$, la fonction $\phi_{\beta_v}$ désigne
la fonction caractéristique de $G(\calO_v)$ dans $G(F_v)$.

\end{itemize}
\end{theoreme}

\dem Soit $\gamma_0$ une classe de conjugaison de $G(F)$ qui est
$(T,T')$-admissible et $s$-norme en $T$. Soit $(V;\tau,\tau')$ un
objet de $C(T,T';s)$ qui correspond à $\gamma_0$. Donner un objet
$$(\calV;t,t';\iota)\in {\mathfrak C}^I_{\alpha_T,\beta_{T'}}
(T,T';s)
$$
de fibre générique $(V;\tau,\tau')$ revient à donner en toutes les
places $v\in |X|$, un $\calD_v\hat\otimes_{\FF_q}k$-réseau
$\calV_v$, plus une $I$-structure de niveau  tels que
\begin{itemize}
    \item en une place $v\notin |T|$, $\calV_v$ est fixe par $\tau$
    \item en une place $v\notin |T'|$, $\calV_v$ est fixe par
    $\tau'$
    \item en une place $x\in |T|$, la position relative entre
    $\calV_x$ et $\tau(\calV_x)$ est donnée par $\alpha_x$
    \item en une place $v\in |T'|$, donc hors de $T$, les
    descentes à l'aide de $\tau$ de $F_v\hat\otimes_{\FF_q}k$ à $F_v$
    de $\calV_v$ et $\tau'(\calV_v)$, munis de leurs structure de niveaux, sont
    en position relative $\beta_v$.
\end{itemize}
Il s'agit de compter le nombre de ces réseaux modulo la relation
d'équivalence induite par la multiplication scalaire par l'idèle
$a\in \bbA_F^\times$ et par l'action de
$${\rm Aut}(V;\tau,\tau')=J_{\gamma_0,\delta_T}(F).$$
Que ce nombre est égal à l'expression intégrale dans l'énoncé du
théorème, résulte de la définition même des intégrales
orbitales.\findem

\section{L'identité de changement de base}

On va donner la démonstration du théorème 1 de 3.3 dans les
paragraphes 5.1-5.4 de ce chapitre. On va en tirer une nouvelle
démonstration du lemme fondamental pour le changement de base en
5.7.

Rappelons qu'on a supposé que $D$ a plus de $d^2 r||\lambda||$
places totalement ramifiées. Ceci implique la propreté des deux
morphismes caractéristiques $c^I_\lambda$ et $c^I_{r.\lambda}$
d'après 1.6.

\subsection{Théorème de densité de Chebotarev}
\setcounter{lemme}{0}

Soit $U$ l'ouvert de $(X'-I)^{r}$ complémentaire de la réunion des
diagonales. Le groupe symétrique $\frakS_r$ agissant sur
$(X'-I)^{r}$ laisse stable $U$, et son action y est libre. En
particulier, l'action sur $U$ du sous-groupe cyclique $\la\tau\ra$
engendré par la permutation cyclique $\tau$, est libre.
Considérons le quotient $\Utau$ de $U$ par l'action de
$\la\tau\ra$. Ses $k$-points sont les orbites de $\la\tau\ra$ dans
$U(k)$. Pour tout point géométrique $u\in U(k)$, on note
$[u]\in\Utau(k)$ son orbite sous l'action de $\la \tau\ra$.

Choisissons un point géométrique $v$ de $U$. Le revêtement
galoisien $U \rta \Utau$ induit une suite exacte de groupes
fondamentaux
$$1\rta \pi_1(U,v) \rta \pi_1(\Utau,[v])\rta
\la\tau\ra\rta 1.
$$
Cette suite n'est pas nécessairement scindable.

Par la descente, la catégorie des systèmes locaux sur $U$ munis
d'une action compatible de $\la\tau\ra$ est équivalente à la
catégorie des systèmes locaux sur $\Utau$ et donc équivalente à la
catégorie des représentations continue du groupe
$\pi_1(\Utau,[v])$.

Comme dans 3.1, fixons un point $x\in (X'-I)(k)$ et notons $x^r\in
(X'-I)(k)$ le point diagonal de coordonnées $x$. En prenant la
fibre en $x^r$, la catégorie des systèmes locaux sur $(X'-I)^{r}$
munis d'un action compatible de $\la\tau\ra$ est équivalente à la
catégorie des représentations continues du produit semi-direct
$\pi_1((X'-I)^{r},x^r)\rtimes \la\tau\ra$.

Le foncteur de restriction de $(X'-I)^{r}$ à $U$ fournit un
homomorphisme de groupes
$$\pi_1(\Utau,[v]) \rta \pi_1((X'-I)^r,x^r)\rtimes \la\tau\ra$$
bien déterminé à un automorphisme intérieur près. Le foncteur de
restriction étant fidèle, l'homomorphisme qui s'en déduit est
surjectif. Il est de plus compatible avec la projection sur le
groupe cyclique $\la\tau\ra$. Le choix de ce chemin de $v$ à $x^r$
nous permet de reformuler le théorème 1 de 3.3 comme suit. Pour
tout $g\in \pi_1(\Utau,v)$ dont l'image dans le groupe cyclique
$\la\tau\ra$ est $\tau$, pour tout $f\in \calH^I$, on a l'égalité
$$\Tr((1\otimes\cdots\otimes 1\otimes f)g,[A_{[v]}])
=\Tr(fg,[B_{[v]}]).$$

Un point fermé de $\Utau$ de corps résiduel $\FF_{q^s}$ est une
orbite $[u]$ telle que $\sigma^s$ soit la plus petite puissance de
$\sigma$ laissant stable $[u]$. Soit $u\in U(k)$ un point
géométrique de $U$ tel que $[u]$ soit un point fermé de $\Utau$ de
corps résiduel $\FF_{q^s}$. Il existe un unique $j\in \ZZ/r\ZZ$
tel que
$$\sigma^s(u)=\tau^j(u).$$
Puisque les opérateurs $\sigma$ et $\tau$ commute, on a
$\sigma^s(u')=\tau^j(u')$ pour tout $u'$ dans l'orbite $[u]$. Le
point fermé $[u]$ de $\Utau$ définit une classe de conjugaison
$\Frob_{[u]}$ de $\pi_1(\Utau,[v])$ dont l'image dans $\la\tau\ra$
est exactement $\tau^j$.

{\em Grâce au théorème de Chebotarev}, pour démontrer le théorème
1 de 3.3, il suffit de démon\-trer que pour tout point fermé $[u]$
de $\Utau$ tel que la classe de conjugaison $\Frob_{[u]}\in
\pi_1(\Utau)$ associée ait l'image $\tau$ dans $\la\tau\ra$ et
pour tout opérateur de Hecke $f\in \calH^I$, on a l'égalité
$$ \Tr((1\otimes\cdots\otimes 1\otimes f)\times \Frob_{[u]},
[A_{[u]}])=\Tr(f\times \Frob_{[u]},[B_{[u]}]).$$

\begin{definition}
Un point fermé $[u]$ de $\Utau$ est dit {\bf cyclique} si l'image
de la classe de conjugaison $\Frob_{[u]}\subset \pi_1(\Utau)$ dans
$\la\tau\ra$ est le générateur $\tau$.
\end{definition}

\subsection{Description des points cycliques}
\setcounter{lemme}{0}

\begin{proposition}
Soient $u=(u_i)\in U(k)$, $[u]\in\Utau(k)$ son image dans $\Utau$.
Supposons que le point fermé de $\Utau$ supportant $[u]$ soit de
degré $s$ et soit cyclique. Alors, il existe un unique point fermé
$x\in |X'|$ de degré $rs$ supportant $\{u_1,\ldots,u_r\}$ et il
existe un unique point fermé $y\in |X'\otimes_{\FF_q}\FF_{q^s}|$
au-dessus de $x$ tel que
$y\otimes_{\FF_{q^s}}k=\{u_1,\ldots,u_r\}$.
\end{proposition}

\dem Puisque $[u]$ est un point cyclique de corps résiduel
$\FF_{q^s}$, on a
$$\sigma^s(u_i)=u_{i+1}.$$
L'ensemble $\{u_1,\ldots,u_r\}$ définit donc une orbite de $\la
\sigma^s\ra$ dans $X'(k)$ de sorte qu'il existe un unique point
fermé $y\in |X'\otimes_{\FF_q}\FF_{q^s}|$ tel que
$y\otimes_{\FF_{q^s}}k=\{u_1,\ldots,u_r\}$. Puisque $u_1,\ldots
u_r$ sont deux à deux distincts, le corps résiduel de
$X'\otimes_{\FF_q}\FF_{q^s}$ est nécessairement $\FF_{q^{rs}}$.

Soit $x$ le point fermé de $X$ en-dessous de $y$. Il reste à
démontrer que le corps résiduel de $X$ en $x$ est aussi
$\FF_{q^{rs}}$. Supposons qu'il y a $s'$ points fermés de
$X\otimes_{\FF_q}\FF_{q^s}$ au-dessus de $x$. L'extension
$F\otimes_{\FF_q}\FF_{q^s}$ étant de degré $s$, le degré de
décomposition de la place $x$ pour cette extension est $s'$ si
bien que le degré d'inertie est $s''$ avec $s=s's''$. En
particulier, le corps résiduel de $x$ est $\FF_{q^{rs'}}$.

Le groupe $\langle\sigma\rangle$ agit transitivement sur
l'ensemble des places de $F\otimes_{\FF_q}\FF_{q^s}$ au-dessus de
$x$ et le stabilisateur de chacune des ces places est
$\sigma^{s'}$. On a donc
$$\FF_{q^{rs'}}\otimes_{\FF_{q^{s'}}}\FF_{q^s}=\FF_{q^{rs}}$$
ce qui implique que $r$ et $s''$ sont premiers entre eux.

De plus, le point fermé $y$ étant fixe par $\sigma^{s'}$,
l'opérateur $\sigma^{s'}$ permute l'ensemble des points
géométriques $\{u_1,\ldots u_r\}$ au-dessus de $y$ et consiste
donc en une certaine permutation $\tau'$ de l'ensemble
$\{1,\ldots,r\}$ avec ${\tau'}^{s''}=\tau$. Puisque $s''$ est
premier à $r$, il existe $r''\in\NN$ tel que $s'' r''$ soit congru
à $1$ modulo $r$. On en déduit que $\tau'=\tau^{r''}$ donc $\tau'$
appartient à $\la\tau\ra$. Donc l'orbite $[u_i]$ de $(u_i)$ sous
l'action $\la\tau\ra$ est stable par $\sigma^{s'}$.

Puisque $s$ est par hypothèse, le plus petit entier ayant cette
propriété, $s'$ est nécessairement égal à $s$. On en déduit que le
corps résiduel de $X$ en $x$ est bien $\FF_{q^{rs}}$. \findem

\subsection{Calcul des traces sur $[B]$ aux points cycliques}
\setcounter{lemme}{0}

Soit $(u_i)=u\in U(k)$ dont l'image $[u]$ dans $\Utau$ est un
point cyclique de degré $s$. D'après 5.2, $u_1,u_2,\ldots,u_r$
sont supportés par une place $x$ de $X$ de degré $rs$ et sont les
points géométriques au-dessus d'un point fermé $y$ de
$X\otimes_{\FF_q}\FF_{q^s}$ au-dessus de $x$. On va prendre
$T=\{x\}$.

La fibre du morphisme $c^I_{r.\lambda}:\calS^I_{r.\lambda}\rta
(X'-I)^{r}$ au-dessus de $u$ est l'ensemble
$$\calS^I_{r.\lambda}(u)(k)=
\{t:\,^\sigma \calV^{\bab T} \isom \calV^{\bab T}\mid {\rm
inv}(t)\leq\sum_{i=1}^r \lambda{u_i}\}/a^\ZZ,$$ où
$\calV\in\DFib(k)$ et où $\bab T=T\otimes_{\FF_q}k$ contient $rs$
points parmi lesquels $\{u_1,\ldots,u_r\}\subset \bab T$. Le
diviseur $u_1+\cdots+u_r$ est fixe par $\sigma^s$ de sorte qu'on a
une $\FF_{q^s}$-structure sur $\calS^I_{r.\lambda}(u)$. On peut
appliquer la formule de trace de Grothendieck-Lefschetz pour
obtenir la formule suivante.

\begin{proposition}
Pour tout fermé $T'$ de $X$, disjoint de $T$, et pour toute
fonction $\beta_{T'}$ qui à toute place $v\in |T'|$ associe une
double-classe $$\beta_v\in K^I_v\backslash G(F_v)/K^I_v,$$ on a
l'égalité
$$\begin{array}{c}
\Tr(\sigma^s\circ
\Phi_{\beta_{T'}},\rmR\Gamma(\calS^I_{r.\lambda}(u),
\calF_{r.\lambda}))\\ = \sum_{z\in {\rm Fix}(\sigma^s\circ
\Phi_{\beta_{T'}})}\#{\rm Aut}(z)^{-1} \Tr(\sigma^s\circ
\Phi_{\beta_{T'}}, (\calF_{r.\lambda})_z).
\end{array}$$
\end{proposition}

Par définition, un objet de la catégorie des points fixes de la
correspondance $\sigma^s\circ \Phi_{\beta_{T'}}$ dans
$\calS^I_{r.\lambda}(u)$ consiste en un chtouca $(\calV,t,\iota)$
d'invariant de Hodge $\inv(t)\leq \sum_{i=1}^r \lambda u_i$ et
muni d'une $I$-structure de niveau $I$, plus d'une
$T'$-modification $t':\,^{\sigma^s}\calV^{T'}\isom \calV^{T'}$,
commutant à $t$ en qui en une place $v\in T'$ ayant le type
$\beta_v$. Notons que les $u_i$ sont les points géométriques
au-dessus d'un seul point fermé $y\in x\otimes_{\FF_q}\FF_{q^s}$
d'après 5.2.

\begin{proposition}
Soit $(\calV,t,t',\iota)\in {\rm Fix}(\sigma^s\circ
\Phi_{\beta_{T'}}, \calS^I_{r.\lambda}(u))$, alors il existe
$\alpha\in\ZZ^d_+$ avec $\alpha\leq\lambda$ tel qu'on a
$$\inv(t)=\sum_{i=1}^r \alpha u_i.$$ De plus, on a
$$\Tr(\sigma^s\circ \Phi_{\beta_{T'}}, (\calF_{r.\lambda})_{(\calV;t,t')})
=P_{\lambda,\alpha}(q^{rs})
$$
où $P_{\lambda,\alpha}(q^{rs})$ est définie comme la trace de
$\sigma_{q^{rs}}$ sur la fibre de $\calA_{\lambda,u_i}$ au-dessus
du point $\alpha(u_i)$ de $\calQ_{\lambda}(u_i)$, correspondant au
copoids $\alpha$ et pour n'importe quel $i=1,\ldots,r$.
\end{proposition}

\dem Puisque $\sigma^s(u_i)=u_{i+1}$ pour tous $i=1,\ldots,r$, on
a une $\FF_{q^s}$-structure de la fibre
$\calQ_{r.\lambda}(u)=\prod_{i=1}^n\calQ_{\lambda,u_i}$ de
$\calQ_{r.\lambda}$ au-dessus de $(u)$, donnée par
$\sigma^s(Q^i_j)=({Q'}^i_j)$ avec ${Q'}^i_j=Q^{i+1}_j$. Le
morphisme
$$f_{r.\lambda}:\calS^I_{r.\lambda}(u)
\rta \calQ_{r.\lambda}(u)$$ est compatible aux
$\FF_{q^s}$-structures. Un point $(\calV;t,t')$ fixe par
$\sigma^s\circ \Phi_{\beta_{T'}}$ doit nécessairement s'envoyer
sur un point $Q$ de $\calQ_{r.\lambda}(u)$ fixe par $\sigma^s$. On
a aussi
$$\Tr(\sigma^s\circ \Phi_{\beta_{T'}}, (\calF_{r.\lambda})_{(\calV;t,t')})
=\Tr(\sigma^s,(\calA_{r.\lambda})_Q).$$

Il est clair que la $\FF_{q^s}$-champ $\prod_{i=1}^r
\calQ_{\lambda,u_i}$ est la restriction de scalaires à la Weil du
$\FF_{q^{rs}}$-champ $\calQ_{\lambda,u_i}$ pour n'importe quel
$i$. Le $\FF_{q^s}$-point $Q$ de $\prod_{i=1}^r
\calQ_{\lambda,u_i}$ correspond donc à un $\FF_{q^{rs}}$-point
$Q'$ de $\calQ_{\lambda,u_i}$. On a alors
$$\Tr(\sigma^s,(\calA_{r.\ul})_Q)=\Tr(\sigma^{rs},(\calA_{\ul})_{Q'}).$$
Il s'agit là d'une propriété complètement générale de la
restriction à la Weil dont voici l'énoncé précis. Sa démonstration
se ramène essentiellement à la formule de Saito-Shintani dans 3.3.

\begin{proposition}
Soient $Z$ un schéma défini sur $\FF_{q^s}$ et $\mathcal G$ un
complexe de $\QQ_\ell$-faisceaux sur $Z$. Soit
$Z'=\rmR_{\FF_{q^s}/\FF_q}Z$ et $\mathcal
G'=\rmR_{\FF_{q^s}/\FF_q}\mathcal G$ les restriction de scalaires
à la Weil de $Z$ et de $\mathcal G$. Soient $z\in Z(\FF_{q^s})$ un
$\FF_{q^s}$-point de $X$ et $z'$ le $\FF_{q}$-point de $Z'$ qui
lui correspond. Alors, on a
$$\Tr(\sigma^s,\calF_z)=\Tr(\sigma,\calF'_{z'}).$$
\end{proposition}

\begin{corollaire}
On a la formule
$$\Tr(\sigma^s\circ \Phi_{\beta_{T'}},\rmR\Gamma(\calS^I_{r.\lambda}(u),
\calF_{r.\lambda})) =\sum_{\alpha\leq\lambda}
P_{\lambda,\alpha}(q^{rs})
\,\#\frakC^I_{\alpha_{y},\beta_{T'}}(T,T';s).
$$
où $\alpha_y$ est la fonction $|T\otimes_{\FF_q}\FF_{q^s}|\rta
\ZZ^d_+$ qui prend la valeur $\alpha$ en $y$ et $0$ sur les $s-1$
autres points de $|T\otimes_{\FF_q}\FF_{q^s}|$.
\end{corollaire}

En appliquant le théorème 1 de 4.7, et en tenant compte du fait
que l'extension $F\otimes_{\FF_q}\FF_{q^s}$ de $F$ est totalement
décomposée en $x$, on obtient la formule suivante pour le cardinal
$\#\frakC^I_{\alpha_y,\beta_{T'}}(T,T';s)$
$$\#\frakC^I_{\alpha_y,\beta_{T'}}(T,T';s)= \sum_{\gamma_0}{\rm
vol}(G_{\gamma_0}(F)a^\ZZ\backslash G_{\gamma_0}(\bbA_F))
\prod_{v\not=x}{\mathbf O}_{\gamma_0}(\phi_{\beta_v}) {\mathbf
O}_{\gamma_0}(\phi_{\alpha_x}).
$$
où $\phi_{\alpha_x}\in\calH(G(F_x))$ est la fonction
caractéristique de la $G(\calO_x)$-double classe correspondant à
$\alpha$.

\begin{corollaire}
On a la formule
$$\begin{array}{c}
\Tr(\sigma^s\circ
\Phi_{\beta_{T'}},\rmR\Gamma(\calS^I_{r.\lambda}(u),
\calF_{r.\lambda}))=\\
\sum_{\gamma_0} {\rm vol}(a^\ZZ G_{\gamma_0}(F)\backslash
G_{\gamma_0}(\bbA_F)) \prod_{v\not= x} {\mathbf
O}_{\gamma_0}(\phi_{\beta_v}) {\mathbf
O}_{\gamma_0}(\psi_{\lambda_x}).
\end{array}$$
\end{corollaire}

En effet, on a $\psi_{\lambda_x}=\sum_{\alpha\leq\lambda}
P_{\lambda,\alpha}(q^{rs}) \phi_{\alpha_x} $ car le point $x$ est
de degré $rs$.

\subsection{Calcul des traces sur $[A]$ aux points cycliques}
\setcounter{lemme}{0}

Prenons le même point $u=(u_1,\ldots,u_r)$ de $U$ comme dans 5.3
dont l'image $[u]$ dans $\Utau$ définit un point cyclique de degré
$s$. Les points $u_1,\ldots,u_r$ sont supportés par un point fermé
$y$ de $X\otimes_{\FF_q}\FF_{q^s}$ de degré $rs$ sur $\FF_q$.
Quant à $y$, il est au-dessus d'un point fermé $x$ de $X$ de degré
$rs$ aussi. On va noter $T=\{x\}$.

La fibre du morphisme $(c^I_\lambda)^r :(\calS^I_\lambda)^r \rta
(X'-I)^{r}$ au-dessus de $u$ est munie d'une $\FF_{q^s}$-structure
donnée par
$$\sigma^s(\tilde\calV_1,\ldots,\tilde\calV_r)=
(\sigma^s(\tilde\calV_r),\sigma^s(\tilde\calV_1),\ldots,\sigma^s
(\tilde\calV_{r-1}))$$ avec $\tilde\calV_i\in
\calS^I_\lambda(u_i)$. On peut relever cette $\FF_{q^s}$-structure
à la restriction de $\calF_{\lambda}^{\boxtimes r}$ à cette fibre.

Soient $T'$ et $\beta_{T'}$ comme dans le paragraphe précédent. En
appliquant la formule des traces de Grothendieck-Lefschetz, on
obtient la formule.

$$\begin{array}{c}
{\rm Tr}(\sigma^s\circ (1\otimes \cdots\otimes 1\otimes
\Phi_{\beta_{T'}}),
\rmR\Gamma((\calS^I_\lambda)^r(u),\calF_{\lambda}^{\boxtimes r}))\\
=\sum_{z\in {\rm Fix}}\#{\rm Aut}(z)^{-1} \Tr(\sigma^s\circ
(1\otimes\cdots\otimes 1 \otimes \Phi_{\beta_{T'}}),
(\calF_\lambda^{\boxtimes r})_z).
\end{array}$$
où ${\rm Fix}$ désigne l'ensemble des classes d'isomorphisme
d'objets de la catégorie des points fixes ${\rm Fix}(\sigma^s\circ
(1\otimes\cdots\otimes 1\otimes
\Phi_{\beta_{T'}}),(\calS^I_\lambda)^r(u))$.

Les objets de cette catégorie des points fixes sont constitués de
\begin{itemize}
\item $\tilde\calV_i\in \DCht^I_\lambda((u_i))$ pour tous
$i=1,\ldots,r$, \item des isomorphisme
$\tilde\calV_2\isom\sigma^s(\tilde\calV_1),\ldots,
\tilde\calV_r\isom \sigma^s(\tilde\calV_{r-1})$ et une
$T'$-modifi\-ca\-tion $\sigma^s(\calV^r)^{T'}\isom \calV_1^{T'}$
qui en toute place $v\in |T'|$ soit de type $\beta_v$.
\end{itemize}
Il revient au même de se donner un point
$\tilde\calV=(\calV,t,\iota)\in\DCht^I_\lambda(u_1)(k)$ et d'une
$T'$-modification $t':\,^{\sigma^{rs}}\tilde\calV^{T'} \rta
\tilde\calV^{T'}$ qui ait le type $\beta_v$ en toute place
$v\in|T'|$.

Comme dans le paragraphe précédent, on a les assertions suivantes.

\begin{proposition}
Soit  $z$ un point fixe de $\sigma^s\circ (1\otimes\cdots\otimes
1\otimes\Phi_{\beta_{T'}})$ dans $(\calS^I_\lambda)^r(u)$. A $z$
est associé $(\calV,t,\iota)\in\DCht^I_\lambda(u_1)(k)$ muni de
$t':\,^{\sigma^{rs}}\tilde\calV^{T'} \rta \tilde\calV^{T'}$ comme
précédemment. Alors on a la formule
$$\Tr(\sigma^s\circ (1\otimes\cdots\otimes 1\otimes \Phi_{\beta_{T'}}),
(\calF_\ul^{\boxtimes r})_z)= P_{\lambda\alpha}(q^{rs})$$ où
$\inv(t)= \alpha u_1$.
\end{proposition}

\begin{corollaire}
On a la formule
$$\begin{array}{c}
\Tr(\sigma^s\circ (1\times \cdots\times 1\times
\Phi_{\beta_{T'}},\rmR\Gamma((\calS^I_{\lambda})^r(u),
(\calF^I_\lambda)^{\boxtimes r}))\\
=\sum_{\alpha\leq\lambda} P_{\lambda\alpha}(q^{rs})
\,\#\frakC^I_{\alpha_{u_1},\beta_{T'}}(T,T';rs)
\end{array}$$
où $\alpha_{u_1}$ est la fonction
$|T\otimes_{\FF_q}\FF_{q^{rs}}|\rta\ZZ^d_+$ qui prend la valeur
$\alpha$ en $u_1$ et $0$ ailleurs.
\end{corollaire}

Les catégories
$$
\frakC^I_{\alpha {u_1},\beta_{T'}}(T,T';rs)\quad
\mbox{et}\quad\frakC^I_{\sum_{i=1}^r\alpha
{u_i},\beta_{T'}}(T,T';s)
$$
ne semblent pas équivalentes de façon naturelle. Néanmoins, la
formule de comptage leur donne le même cardinal. En effet, en
appliquant le théorème 1 de 4.7, et en tenant compte du fait que
l'extension $F\otimes_{\FF_q}\FF_{q^{rs}}$ de $F$ est totalement
décomposée en $x$, on obtient la formule
$$\#\frakC^I_{\alpha_{u_1},\beta_{T'}}(T,T';rs)=
\sum_{\gamma_0}{\rm vol}(G_{\gamma_0}(F)a^\ZZ\backslash
G_{\gamma_0}(\bbA_F)) \prod_{v\not=x}{\mathbf
O}_{\gamma_0}(\phi_{\beta_v}) {\mathbf
O}_{\gamma_0}(\phi_{\alpha_x})
$$
où $\phi_{\alpha_x}\in\calH(G(F_x))$ est la fonction
caractéristique de la $G(\calO_x)$-double classe correspondant à
$\alpha_x$.

On en déduit l'égalité
$$ \Tr(\Frob_{[u]}\times (1\otimes\cdots\otimes 1\otimes f)\times,
[A_{[u]}])=\Tr(\Frob_{[u]}\times f,[B_{[u]}])$$ pour tous les
points cyclique $u=(u_i)$. Le théorème 1 de  3.3 s'en déduit par
le théorème de densité de Chebotarev comme il a été expliqué dans
5.1. \findem

\subsection{Calcul sur la petite diagonale : situation $A$}
\setcounter{lemme}{0}

Soit $x$ un point fermé de degré $1$ de $X'-I$, $\bab x$ le point
géométrique au-dessus de $x$ et notons
$$r.{\bab x}:=(\underbrace{\bab x,\ldots,\bab x}_{r})\in (X'-I)^r(k)$$
le point sur la petite diagonale correspondant.

La fibre $(\calS^I_\lambda)^r(r.\bab x)$ est la catégorie des
$(\tilde\calV_i)_{i=1}^r$ où $\tilde\calV_i
=(\calV_i,t_i,\iota_i)\in \calS^I_\lambda(\bab x)$. Puisque $\bab
x$ est défini sur $\FF_{q}$, cette fibre a une $\FF_{q}$-structure
évidente. Elle est de plus munie d'un automorphisme $\tau$ d'ordre
$r$ qui permute les facteurs de manière cyclique.

Soient $T'$ un fermé de $X'$ évitant $x$, et $\beta_{T'}$ une
fonction qui associe à toute $v\in |T'|$ une double-classe
$\beta_v\in K^I_v \backslash G(F_v)/K^I_v$. Les points fixes de la
correspondance $\tau\circ\sigma\circ (1\times\cdots\times 1\times
\Phi_{\beta_{T'}})$ dans la fibre de $(\calS^I_\lambda)^r(r.\bab
x)$ sont les $(\tilde\calV_i)_{i=1}^r$ muni
\begin{itemize}
\item
 des isomorphismes
$^{\sigma}\calV_1\isom \calV_2,\ldots,\,^{\sigma}\calV_{r-1}\isom
\calV_r$

\item
 d'une $T'$-modification $t':\,^{\sigma}\calV_r^{\bab T'} \isom
\calV_1^{\bab T'}$, commutant à $t$, qui en toute place $v\in T'$
ait le type $\beta_v$.
\end{itemize}
Il revient donc au même de donner $\calV\in\calS^I_\lambda(\bab
x)$ muni d'une modification $t':\,^{\sigma^{r}} \calV^{\bab
T'}\isom \calV^{\bab T'}$ qui en toute place $v\in T'$ ait le type
$\beta_v$.

\begin{corollaire}
On a
$$\begin{array}{c}
\Tr(\tau\circ \sigma \circ (1\times\cdots\times 1\times
\Phi_{\beta_{T'}}),[A_{r.\bab x}])=\\
\sum_{(\gamma_0,\delta_x)}{\rm vol}(J_{\gamma_0}(F)a^\ZZ
\backslash J_{\gamma_0}(\bbA_F)) \prod_{v\not=x} {\mathbf
O}_{\gamma_0}(\phi_{\beta_v})
{\mathbf{TO}}_{\delta_x}(\psi_{\lambda_{x\otimes_{\FF_q}\FF_{q^r}}})
\end{array}$$
où
\begin{itemize}
\item  $(\gamma_0,\delta_x)$ parcourt l'ensemble des paires
constituées d'une classe de conjugaison $\gamma_0$ de $G(F)$ et
d'une classe de $\sigma$-conjugaison
$$\delta_x\in G((F\otimes_{\FF_q}\FF_{q^{r}})_x)$$
dont la norme est $\gamma_0$,

\item la fonction $\psi_{\lambda_{x\otimes_{\FF_q}\FF_{q^r}}}$ est
la fonction sphérique correspondant à la double classe $\lambda$
dans l'algèbre de Hecke $G((F\otimes_{\FF_q}\FF_{q^{r}})_x)$.
\end{itemize}
\end{corollaire}

\subsection{Calcul sur la petite diagonale : situation $B$}
\setcounter{lemme}{0}

Soit $x,\bab x$ et $r.\bab x$ comme dans 5.5. On va maintenant
calculer la trace sur la fibre de $B$ au-dessus du point $r.\bab
x=(\bab x,\ldots,\bab x)$ ; ce calcul est un peu plus difficile
que dans le cas $A$ car on doit utiliser la théorie des modèles
locaux et l'interprétation géométrique de l'homomorphisme de
changement de base, voir \cite{N}.

Rappelons que l'action de $\la\tau\ra$ sur la fibre $[B]_{r.\bab
x}$ a été obtenu par prolongement d'une action définie
géométriquement sur restriction de $[B]$ sur l'ouvert $U$. A
priori, on ne connaît pas la trace du Frobenius tordu par $\tau$
en dehors de l'ouvert $U$. Les modèles locaux sont utiles dans ces
situations de mauvaise réduction.

Soit ${X'}^{[r]}={X'}^r/\frakS_r$ la $r$-ème puissance symétrique
de $X'$. Un $k$-point $\bab T$ de ${X'}^{[r]}$ est un diviseur
effectif de degré $r$ de $X'$. Pour tout $\lambda\in\ZZ^d_+$, on
pose $\lambda\bab T:= \sum (m_i\lambda) x_i$ où $\bab T=\sum
m_ix_i$. On a un morphisme
$$c^I_{[r\lambda]}:\calS^I_{[r\lambda]} \rta {X'}^{[r]}$$
où la fibre de $c^I_{[r\lambda]}$ au-dessus d'un point $\bab T\in
{X'}^{[r]}$ est la catégorie en groupoïdes des triplets
$(\calV,t,\iota)$ où $\calV\in\DFib(k)$, où
$t:\,^\sigma\calV^{\bab T} \isom\calV^{\bab T}$ est une $\bab
T$-modification d'invariant $\inv(t)\leq \lambda\bab T$ et où
$\iota$ est une $I$-structure de niveau.

On a un diagramme commutatif
$$\renewcommand\arraystretch{1.5}
\begin{array}{ccc}
\calS^I_{r.\lambda} & \hfld{c^I_{r.\lambda}}{} & {X'}^r \\
\vfld{\pi_\calS}{} & & \vfld{}{\pi} \\
\calS^I_{[r\lambda]} &\hfld{}{c^I_{[r\lambda]}} & {X'}^{[r]}
\end{array}$$
où $\pi_\calS$ envoie un point
$$((x_i)_{i=1}^r,(\calV_i)_{i=0}^r\,,\, t_i:\calV_i^{x_i}\isom \calV_{i-1}^{x_i}
\,,\,\epsilon:\,^\sigma\calV_0\isom \calV_r, \iota)$$ de
$\calS^I_{r.\lambda}$ au-dessus de $(x_1,\ldots,x_n)$, sur le
point $(t:\,^\sigma\calV_0^{\bab T}\rta \calV_0^{\bab T},\iota)$
de $\calS^I_{[r\lambda]}$ au-dessus de $\bab T=x_1+\cdots+x_r\in
{X'}^{(r)}$, où $t$ est la modification composée $t:=\epsilon\circ
t_r\circ\cdots\circ t_1$. D'après la propriété de factorisation du
paragraphe 1.2, ce diagramme est cartésien au-dessus de l'ouvert
$\pi(U)$ où $U$ est le complémentaire dans ${X'}^r$ de la réunion
des diagonales.

On a de plus un diagramme cartésien
$$\renewcommand\arraystretch{1.5}
\begin{array}{ccc}
\calS^I_{r.\lambda} & \hfld{f^I_{r.\lambda}}{} & \calQ_{r.\lambda} \\
\vfld{\pi_\calS}{} & & \vfld{}{\pi_\calQ} \\
\calS^I_{[r\lambda]} &\hfld{}{f^I_{[r\lambda]}} & \calQ_{[r\lambda]}
\end{array}$$
où $Q_{r.\lambda}$ est le champ des
$$\{(x_1,\ldots,x_r);(0=Q_0\subset Q_1\subset Q_2\cdots\subset Q_r)\}$$
où $x_1,\ldots,x_r\in X'$, où $Q_i$ sont des $\calD$-modules de
torsion tels que $Q_i/Q_{i-1}$ est supporté par $x_i$ et
$\inv_{x_i}(Q_i/Q_{i-1})\leq \lambda$. Le champ
$\calQ_{[r\lambda]}$ classifie les paires $(T,Q)$ où $T\in
X^{(r)}$ et $Q$ est un $\calD$-module de torsion tel que
$\inv(Q)\leq \lambda T$. Le morphisme $\pi_\calQ:
\calQ_{r.\lambda} \rta \calQ_{[r\lambda]}$ est défini par
$$\{(x_1,\ldots,x_r);(0=Q_0\subset Q_1\subset Q_2\cdots\subset Q_r)\}
\mapsto (T,Q)$$
avec $T=x_1+\cdots+x_r$ et $Q=Q_r$.

\begin{proposition}
Le morphisme $\pi_\calQ:\calQ_{r.\lambda}\rta \calQ_{[r\lambda]}$
est un morphisme petit au sens stratifié. En particulier, l'image
directe $\rmR(\pi_\calQ)_* \calA_{r.\lambda}$ est un faisceau
pervers prolongement intermédiaire de sa restriction à l'ouvert
où il est un système local.
\end{proposition}

Notons que la notion du morphisme petit au sens stratifié de
Mirkovic et Vilonen, pour les schémas, voir \cite{MV}, est
conservé par un changement de base lisse. Il est donc bien défini
pour les champs d'Artin. Il revient donc au même de démontrer la
proposition après un changement de base lisse surjectif comme dans
1.1.5 pour se ramener aux Grassmanniennes affines dans quel cas
elle est démontrée dans \cite{MV} et \cite{NP}.

Au-dessus de l'ouvert $\pi(U)$, le diagramme
$$\renewcommand\arraystretch{1.5}
\begin{array}{ccc}
 \calQ_{r.\lambda} & \hfld{}{} & {X'}^r \\
\vfld{}{} & & \vfld{}{} \\
\calQ_{[r\lambda]} &\hfld{}{} &  {X'}^{[r]}
\end{array}$$
est cartésien de sorte qu'on a une action de $\frakS_r$ sur la
restriction de $\rmR(\pi_\calQ)_* \calA_{r.\lambda}$ à l'image
inverse de $\pi(U)$ dans $\calQ_{[r\lambda]}$. On en déduit une
action de $\frakS_r$ sur tout le faisceau pervers
$\rmR(\pi_\calQ)_* \calA_{r.\lambda}$ par fonctorialité du
prolongement intermédiaire.

D'après \cite{N}, l'action induit de la permutation cyclique
$\tau$ sur ce faisceau pervers est reliée à l'homomorphisme de
changement de base d'algèbres de Hecke. La fibre de
$\calQ_{[r\lambda]}$ au-dessus du point $T=rx$ ont les points
$Q_{\alpha(x)}$ paramétrés par $\alpha\leq r\lambda$, tous définis
sur $\FF_{q}$. On définit une fonction de l'algèbre de Hecke
$\calH_x$ en prenant la combinaison linéaire
$$\psi_{[r\lambda](x)}=\sum_{\alpha\leq r\lambda} \Tr(\sigma\circ \tau,
(\rmR(\pi_\calQ)_*
\calA_{r.\lambda})_{Q_{\alpha(x)}})\phi_{\alpha_x}
$$
où $\phi_{\alpha_x}\in\calH(G(F_x))$ est la fonction
caractéristique de la double classe correspondant à $\alpha$. On a
l'énoncé suivant, voir \cite{N} théorème 3.

\begin{proposition}
Soit $y$ un point fermé de $X\otimes_{\FF_q}\FF_{q^{r}}$ au-dessus
de $x$. Soit $\calH_y$ l'algèbre de Hecke du groupe
$G((F\otimes_{\FF_q}\FF_{q^{r}})_y)$ et ${\mathbf b}:\calH_y\rta
\calH_x$ l'homomorphisme de changement de base. Alors on a
l'égalité
$$\psi_{[r\lambda](x)}={\mathbf b}(\psi_{\lambda(y)}).$$
\end{proposition}

Avant de continuer notre discussion, rappelons qu'on a un
diagramme :
$$\renewcommand\arraystretch{1.5}
\begin{array}{ccccc}
\calQ_{r.\lambda} & \hflg{f^I_{r.\lambda}}{} &
\calS^I_{r.\lambda} & \hfld{c^I_{r.\lambda}}{} & (X'-I)^r\\
\vfld{\pi_\calS}{} & & \vfld{}{\pi_\calQ} & \searrow{\scriptstyle h}& \vfld{}{\pi}\\
\calQ_{[r\lambda]}&\hflg{}{f^I_{[r\lambda]}} &
\calS^I_{[r\lambda]} & \hfld{}{c^I_{[r.\lambda]}} & (X'-I)^{[r]}
\end{array}$$
dont le carré de gauche est cartésien et dont le carré de droite
est commutatif. Notons $h$ la flèche diagonale du carré de droite.

Le complexe d'intersection
$\calF^I_{r.\lambda}=(f^I_{r.\lambda})^* \calA_{r.\lambda}$ étant
localement acyclique par rapport au morphisme $c^I_{r.\lambda}:
\calS^I_{r.\lambda}\rta (X'-I)^r$, les faisceaux de cohomologie
$$B_i=\rmR^i (c^I_{r.\lambda})_* \calF^I_{r.\lambda}$$
sont des systèmes locaux. Rappelons qu'on a une action de
$\frakS_r$ sur la restriction de $B_i$ à l'ouvert $U$ compatible à
l'action de ce groupe sur $U$. Cette action s'étend à tout $B_i$
qui est un système local.

Puisque $\pi:(X'-I)^r\rta (X'-I)^{[r]}$ est un morphisme fini, on
a
$$\rmR^i h_* \calF^I_{r.\lambda} = \pi_* B_i.$$
En particulier, $\rmR^i h_* \calF^I_{r.\lambda}$ est muni d'une
action de $\frakS_r$ déduite de celle sur $B_i$.

Par ailleurs, $\rmR h_*\calF^I_{r.\lambda}= \rmR
c^I_{[r\lambda],*} \rmR(\pi_\calQ )_*\calF^I_{r.\lambda}$ où
$\rmR(\pi_\calQ )_*\calF^I_{r.\lambda}=(f^I_{[r.\lambda]})^*
\calA_{[r.\lambda]}$, de sorte que $\rmR(\pi_\calQ
)_*\calF^I_{r.\lambda}$ hérite l'action canonique de $\frakS_r$
sur le faisceau pervers $\calA_{[r.\lambda]}$. On en déduit donc
une action de $\frakS_r$ sur tout le complexe $\rmR h
_*\calF^I_{r.\lambda}$, et en particulier sur ses faisceaux de
cohomologie.

\begin{proposition}
Les deux actions de ${\mathfrak S}_r$ sur $\rmR^i
\chi_*\calF^I_{r.\lambda}$ construites ci-dessus sont identiques.
\end{proposition}

\dem Soit $j$ l'immersion ouverte de l'image $\pi(U)$ dans
$(X'-I)^{[r]}$. Puisque $\pi$ est un morphisme fini, on a
$$\pi_*B_i =j_* j^* \pi_* B_i$$
où $j_*$ est le foncteur image directe sur les faisceaux
ordinaires. Il suffit alors de vérifier que les deux actions de
$\frakS_r$ coïncident sur $\pi(U)$ ce qui est une tautologie.
\findem

\bigskip Au-dessus du point $[r\bab x]$ image de $r.\bab x\in {X'}^{r}$
dans ${X'}^{[r]}$ le morphisme fini $\pi:{X'}^r \rta {X'}^{[r]}$
est complètement ramifié c'est-à-dire l'image inverse de $[r\bab
x]$ est un épaississement de $r.\bab x$. Par conséquent la fibre
de $\pi_* B_i$ en $[r\bab x]$ est canoniquement isomorphe à la
fibre de $B_i$ en $r.\bab x$. La compatibilité entre les deux
actions de $\la\tau\ra$ démontrée dans la proposition précédente,
implique l'égalité
$$\Tr(\tau\circ \sigma \circ \Phi_{\beta_{T'}},[B_{r.\bab x}])
=\Tr(\tau\circ \sigma \circ \Phi_{\beta_{T'}},
\rmR\Gamma(\calS^I_{[r\lambda]}([r\bab x]),
\rmR(\pi_\calS)_*\calF^I_{r.\lambda}))$$ où dans le membre droit
l'action de $\tau$ provient de l'action triviale sur l'espace
$\calS^I_{[r\lambda]}([r\bab x])$ et de l'action non triviale de
$\tau$ sur le faisceau pervers
$\rmR(\pi_\calS)_*\calF^I_{r.\lambda}$ provenant de l'action de
$\tau$ sur $\rmR(\pi_{\calQ})_* \calA_{r.\lambda}$.

On est maintenant en position d'appliquer la formule des traces de
Gro\-then\-dieck-Lefschetz
$$
\begin{array}{c}
\Tr(\tau\circ \sigma \circ \Phi_{\beta_{T'}},
\rmR\Gamma(\calS^I_{[r\lambda]}([r\bab x]),\calF^I_{[r\lambda]}))\\
=\sum_{z\in{\rm Fix}(\tau\circ \sigma \circ \Phi_{\beta_{T'}})}
\#{\rm Aut}(z) \Tr(\tau\circ\sigma\circ \Phi_{\beta_{T'}},
(\rmR(\pi_\calS)_*\calF^I_{r.\lambda})_{z}).
\end{array}$$
Si $z=(\calV,t,t',\iota)$, on sait que $\Tr(\tau\circ\sigma\circ
\Phi_{\beta_{T'}},
(\rmR(\pi_\calS)_*\calF^I_{r.\lambda})_{(\calV;t,t')})$ ne dépend
que de l'invariant $\alpha=\inv(t)$ et est égal au nombre
$$\Tr(\tau\circ\sigma,(\rmR(\pi_\calQ)_*\calA_{r.\lambda})_{\alpha(x)})$$
qui est un coefficient de ${\mathbf b}(\psi_{\lambda_{x}})$
exprimé comme combinaison linéaire des $\phi_{\alpha(x)}$ d'après
la proposition 2. On déduit la formule suivante du théorème 1 de
4.7 et de la proposition 2 de cette section.

\begin{corollaire}
On a
$$\begin{array}{c}
\Tr(\tau\circ \sigma \circ \Phi_{\beta_{T'}},[B_{r.\bab x}])\\
=\sum_{\gamma_0} {\rm vol}(G_{\gamma_0}(F)a^\ZZ\backslash
G_{\gamma_0}(\bbA_F)) \prod_{v\not=x} {\mathbf
O}_{\gamma_0}(\phi_{\beta_v}) {\mathbf O}_{\gamma_0}({\mathbf
b}(\psi_{\lambda_{x}}))
\end{array}$$
où la fonction $\psi_{\lambda_{x}}$ est celle définie dans le
corollaire 1 de 5.5.
\end{corollaire}

En appliquant la théorème 1 de 3.3, et en tenant compte du fait
que les fonctions $\phi_{\beta_v}$ peut être prises
arbitrairement, on obtient la formule.

\begin{corollaire}
Soit $\bigotimes_{v\not= x}\phi_v$ une fonction  telle que presque
partout $\phi_v$ est l'unité de l'algèbre de Hecke, alors on a
l'égalité
$$
\begin{array}{c}
\sum_{(\gamma_0,\delta_x)}{\rm vol}(J_{\gamma_0,\delta_x}(F)a^\ZZ
\backslash J_{\gamma_0,\delta_x}(\bbA_F)) \prod_{v\not=x} {\mathbf
O}_{\gamma_0}(\phi_v)
{\mathbf{TO}}_{\delta_x}(\psi_{\lambda_{x}})\\
=\sum_{\gamma_0} {\rm vol}(G_{\gamma_0}(F)a^\ZZ\backslash
G_{\gamma_0}(\bbA_F))\prod_{v\not=x} {\mathbf
O}_{\gamma_0}(\phi_v) {\mathbf O}_{\gamma_0}({\mathbf
b}(\psi_{\lambda_{x}})).
\end{array}
$$
\end{corollaire}

Notons qu'à la différence de la formule globale de Arthur-Clozel
\cite{AC}, en les places $v\not=x$, on n'a que des intégrales
orbitales non-tordues.

\subsection{Application au lemme fondamental} \setcounter{lemme}{0}

On cette section, on va déduire du corollaire 5 de 5.6 le lemme
fondamental pour le changement de base pour des fonctions de Hecke
$\phi_\lambda$ avec $|\lambda|=0$ dans le cas de caractéristique
positive.

\begin{theoreme}
Soient $F_x$ un corps local de caractéristique $p$ et $E_x$ son
extension non ramifiée de degré $r$. Soient $\delta_x$ une classe
de $\sigma$-conjugaison dans $\GL_d(E_x)$ et sa norme $\gamma_x$
une classe de conjugaison de $\GL_d(F_x)$. Supposons que
$\gamma_x$ est une classe semi-simple régulière séparable. Alors
pour toutes $\lambda\in\ZZ^d_+$ avec $|\lambda|=0$, on a
$${\rm TO}_{\delta_x}(\phi_\lambda)={\rm O}_{\gamma_x}({\bf
b}(\phi_\lambda))
$$
où $\phi_\lambda$ est la fonction caractéristique de la
$\GL_d(\calO_{E_x})$-double classe correspondant à $\lambda$.
\end{theoreme}

Ce théorème a été démontré par Arthur-Clozel sans l'hypothèse
$|\lambda|=0$ dans le cas $p$-adique à l'aide de la formule des
traces. Il ne fait pas de doute que leur démonstration puisse être
adaptée en caractéristique $p$. \footnote{Cette démonstration est
connue de Henniart qui ne l'a pas rédigée.}

Avec la méthode utilisée ici, on peut également se passer de
l'hypothèse $|\lambda|=0$ au prix de quelques complications
techniques. On a renoncé à cette généralisation dans la dernière
version de ce texte.

\bigskip
\dem Soit $\FF_q$ le corps résiduel de $F_v$. Soit $X$ une courbe
projective lisse géomé\-tri\-quement connexe sur $\FF_q$ avec un
$\FF_q$-point $x$ de sorte que le complété du corps des fonctions
$F$ de $X$ en $v$ soit isomorphe au corps local $F_v$ de l'énoncé.
L'extension $E_x$ est alors le complété en $x$ de
$F\otimes_{\FF_q}\FF_{q^r}$. Notons $\bab x$ un point géométrique
de $X$ au-dessus de $x$.

Soient $\gamma_x$ un élément semi-simple régulier séparable de
$\GL(F_v)$ et $\delta_x$ un élément de $\GL(E_v)$ dont la norme
est la classe de conjugaison de $\gamma_x$, comme dans l'énoncé.

On va choisir une algèbre à division $D$ de dimension $d^2$ sur
son centre $F$ qui est non ramifié en $x$ mais qui a un nombre de
places totalement ramifiées supérieur ou égal à $rd^2||\lambda||$.
Sous cette hypothèse, on a des morphismes propres
$c^I_\lambda:\calS^I_\lambda\rta (X'-I)$ et
$c^I_{r.\lambda}:\calS^I_{r.\lambda}\rta (X'-I)^r$. On obtient
alors l'égalité entre des sommes d'intégrales orbitales et
intégrales orbitales tordues du corollaire 5 de 5.6.

Puisque les intégrales orbitales (resp. tordues) sont localement
constantes aux voisinages des classes semi-simple régulières
séparables, on peut trouver, par approximation faible, un élément
$\gamma'\in G(F)$ assez proche de $\gamma_x$ pour la topologie
$x$-adique, qui est la norme d'un élément $\delta'_x\in
G(F_x\otimes_{\FF_q}\FF_{q^r})$ et tel que ${\bf
O}_{\gamma_x}({\bf b}(\psi_{\lambda_x})) ={\bf O}_{\gamma'}({\bf
b}(\psi_{\lambda_x}))$ et que ${\bf
TO}_{\delta_x}(\psi_{\lambda_x})={\bf TO}_{\delta'_x}
(\psi_{\lambda_x})$. Il suffit de démontrer que ${\bf
O}_{\gamma'}({\bf b}(\psi_{\lambda_x}))={\bf TO}_{\delta'_x}
(\psi_{\lambda_x})$ car les fonctions $\psi_{\lambda_x}$ et les
fonctions $\phi_{\lambda_x}$ engendrent le même espace vectoriel.

Dans le corollaire 5 de 5.6, on va choisir les fonctions $\phi_v$
en les autres places $v\not=x$ telles que
$${\mathbf O}_{\gamma'}(\phi_v)\not=0.$$
Pour ce choix, on a dans la somme le terme correspondant à
$\gamma_0=\gamma'$ et éventuellement un nombre fini d'autres
classes de conjugaison $\gamma_0$. Quitte à rétrécir d'avantage le
support de $\phi_v$, on peut se débarrasser de ces classes
également tout en préservant la non nullité ${\mathbf
O}_{\gamma'}(\phi_v)\not=0$. Il ne reste alors plus qu'à diviser
des deux côtés de l'égalité par des facteurs ${\mathbf
O}_{\gamma'}(\phi_v)$ pour $v\not=x$ et trouver
$${\bf O}_{\gamma'}({\bf b}(\psi_{\lambda_x}))={\bf TO}_{\delta'_x}
(\psi_{\lambda_x}).$$

\end{document}